\documentclass[final]{siamltex}%
\usepackage[pdftex]{graphicx}
\usepackage{algorithm}
\usepackage{algorithmic}
\usepackage{amsmath}
\usepackage{amsfonts}
\usepackage{amssymb}
\usepackage{stmaryrd}
\usepackage{url}
\usepackage{relsize}
\usepackage{color}
\usepackage{todonotes}
\usepackage{ulem}
\providecommand{\U}[1]{\protect\rule{.1in}{.1in}}
\newtheorem{example}[theorem]{Example}

\begin{document}

\title{Quasi-Newton methods on Grassmannians and multilinear approximations of tensors}
\author{Berkant Savas\thanks{Department of Mathematics, Link\"{o}ping University,
SE-581 83 Link\"{o}ping. \textit{Current address}: Institute for Computational
Engineering and Sciences, The University of Texas at Austin, Austin TX 78712
(\texttt{berkant@cs.utexas.edu})}
\and Lek-Heng Lim\thanks{Department of Mathematics, University of California,
Berkeley, CA 94720\newline(\texttt{lekheng@math.berkeley.edu})}}
\date{\today}
\maketitle

\begin{abstract}
In this paper we proposed quasi-Newton and limited memory quasi-Newton methods
for objective functions defined on Grassmannians or a product of
Grassmannians. Specifically we defined \textsc{bfgs} and \textsc{l-bfgs}
updates in local and global coordinates on Grassmannians or a product of
these. We proved that, when local coordinates are used, our \textsc{bfgs}
updates on Grassmannians share the same optimality property as the usual
\textsc{bfgs} updates on Euclidean spaces. When applied to the best
multilinear rank approximation problem for general and symmetric tensors, our
approach yields fast, robust, and accurate algorithms that exploit the special
Grassmannian structure of the respective problems, and which work on tensors
of large dimensions and arbitrarily high order. Extensive numerical
experiments are included to substantiate our claims.

\end{abstract}



\section{Introduction\label{qng:sec:intro}}

\subsection{Quasi-Newton and limited memory quasi-Newton algorithms on
Grassmannians}

We develop quasi-Newton and limited memory quasi-Newton algorithms for
functions defined on a Grassmannian $\operatorname*{Gr}(n,r)$ as well as
a product of Grassmannians $\operatorname*{Gr}(n_{1},r_{1})\times
\dots\times\operatorname*{Gr}(n_{k},r_{k})$, with \textsc{bfgs} and
\textsc{l-bfgs} updates. These are algorithms along the lines of the class of
algorithms studied by Edelman, Arias, and Smith in \cite{edelm99} and
more recently, the monograph of Absil, Mahony, and Sepulchre in
\cite{absil08}. They are algorithms that respect the Riemannian metric
structure of the manifolds under consideration, and not mere
applications of the usual \textsc{bfgs} and \textsc{l-bfgs} algorithms for
functions on Euclidean space. The actual computations of our \textsc{bfgs} and
\textsc{l-bfgs} algorithms on Grassmannians, like the algorithms in
\cite{edelm99}, require nothing more than standard numerical linear algebra
routines and can therefore take advantage of the many high quality softwares
developed for matrix computations \cite{lapack99, arpack98}. In other words,
manifold operations such as movement along geodesics and parallel transport of
tangent vectors and linear operators do not require actual numerical solutions
of the differential equations defining these operations; instead they are
characterized as matrix operations on local and global coordinate
representations (as matrices) of points on the Grassmannians or points
on an appropriate vector bundle.

A departure and improvement from existing algorithms for manifold optimization
\cite{absil07, absil08, edelm99, gabay82} is that we undertake a \textit{local
coordinates} approach. This allows our computational costs to be reduced to
the order of the \textit{intrinsic dimension} of the manifold as opposed to
the dimension of ambient Euclidean space. For a Grassmannian embedded in
the Euclidean space of $n\times r$ matrices, i.e.$\ \operatorname*{Gr}%
(n,r)\subseteq\mathbb{R}^{n\times r}$, computations in local coordinates have
$r(n-r)$ unit cost whereas computations in global coordinates, like the ones
in \cite{edelm99}, have $nr$ unit cost. This difference becomes more
pronounced when we deal with products of Grassmannians
$\operatorname*{Gr}(n_{1},r_{1})\times\dots\times\operatorname*{Gr}%
(n_{k},r_{k})\subseteq\mathbb{R}^{n_{1}\times r_{1}}\times\dots\times
\mathbb{R}^{n_{k}\times r_{k}}$---for $n_{i}=O(n)$ and $r_{i}=O(r)$, we {have
a computational unit cost of {$O(kr(n-r))$ and $O(krn)$ flops between
the local and global coordinates versions the \textsc{bfgs} algorithms. More
importantly, we will show that our \textsc{bfgs} update in local coordinates
on a product of Grassmannians (and Grassmannian in particular)
shares the same well-known optimality property of its Euclidean counterpart,
namely, it is the best possible update of the current Hessian approximation
{that} satisfies the secant equations and preserves symmetric positive
definiteness (cf.\ Theorem~\ref{thm:Opt}). For completeness and as an
alternative, we also provide the global coordinate version of our
\textsc{bfgs} and \textsc{l-bfgs} algorithms analogous to the algorithms
described in \cite{edelm99}. However the aforementioned optimality is not
possible for \textsc{bfgs} in global coordinates. }}

While we have limited our discussions to \textsc{bfgs} and \textsc{l-bfgs}
updates, it is straightforward to substitute these updates with other
quasi-Newton updates (e.g.\ \textsc{dfp} or more general Broyden class
updates) by applying the same principles in this paper.

\subsection{Multilinear approximations of tensors and symmetric tensors}

In part to illustrate the efficiency of these algorithms, this paper also
addresses the following two related problems about the multilinear
approximations of tensors and symmetric tensors, which are also important
problems in their own right with various applications in analytical chemistry
\cite{smild04}, bioinformatics \cite{ombe07}, computer vision \cite{vasil03},
machine learning \cite{ML, LM}, neuroscience \cite{mo06}, quantum chemistry
\cite{khor07}, signal processing \cite{comon96,latha97,delath04b}, etc. See
also the very comprehensive bibliography of the recent survey \cite{kolda09}.
In data analytic applications, the multilinear approximation of general
tensors is the basis behind the \textit{Tucker model} \cite{tucke66} while the
multilinear approximation of symmetric tensors is used in \textit{independent
components analysis} \cite{delath00a,delath04b} and \textit{principal cumulant
components analysis} \cite{ML,LM}. The algorithms above provide a natural method
to solve these problems that exploits their unique structures.

The first problem is that of finding a best multilinear rank-$(p,q,r)$
approximation to a tensor, i.e.\ approximating a given tensor $\mathcal{A}%
\in\mathbb{R}^{l\times m\times n}$ by another tensor $\mathcal{B}\in
\mathbb{R}^{l\times m\times n}$ of lower multilinear rank \cite{hit:27, dSL},
\[
\min_{\operatorname*{rank}(\mathcal{B})\leq(p,q,r)}\lVert\mathcal{A}%
-\mathcal{B}\rVert.
\]
For concreteness we will assume that the norm in question is the
Frobenius or Hilbert-Schmidt norm $\lVert\,\cdot\,\rVert_{F}$. In
notations that we will soon define, we seek matrices $X,Y,Z$ with orthonormal
columns and a tensor $\mathcal{C}\in\mathbb{R}^{p\times q\times r}$ such that%
\begin{equation}
\operatorname*{argmin}_{X\in\operatorname*{O}(l,p),Y\in\operatorname*{O}%
(m,q),Z\in\operatorname*{O}(n,r),\mathcal{C}\in\mathbb{R}^{l\times m\times n}%
}\lVert\mathcal{A}-(X,Y,Z)\cdot\mathcal{C}\rVert_{F}.\label{eqT}%
\end{equation}

The second problem is that of finding a best multilinear rank-$r$
approximation to a symmetric tensor $\mathcal{S}\in\mathsf{S}^{3}%
(\mathbb{R}^{n})$. In other words, we seek a matrix $Q$ whose columns are
mutually orthonormal, and a symmetric tensor $\mathcal{C}\in\mathsf{S}%
^{3}(\mathbb{R}^{r})$ such that a multilinear transformation of $\mathcal{C}$
by $Q$ approximates $\mathcal{S}$ in the sense of minimizing a sum-of-squares
loss. Using the same notation as in \eqref{eqT}, the problem is
\begin{equation}
\operatorname*{argmin}_{Q\in\operatorname*{O}(n,r),\mathcal{C}\in
\mathsf{S}^{3}(\mathbb{R}^{r})}\lVert\mathcal{S}-(Q,Q,Q)\cdot\mathcal{C}%
\rVert_{F}. \label{eqS}%
\end{equation}
This problem is significant because many important tensors that arise in
applications are symmetric tensors.

We will often refer to the first problem as the \textit{general case} and the
second problem as the \textit{symmetric case}. Most discussions are
presented for the case of $3$-tensors for notational simplicity
but key expressions are given for tensors of arbitrary order
to facilitate structural analysis of the problem and algorithmic
implementation. The \textsc{matlab} codes of all algorithms in this
paper are available for download at \cite{savas08a, savas08b}. All of
our implementations will handle $3$-tensors and in addition,
our implementation of the \textsc{bfgs} with scaled identity as
initial Hessian approximation will handle tensors of arbitrary order.
In fact the reader will find an example of a tensor of order-$10$ in
Section~\ref{qng:sec:comEx}, included to show that our algorithms indeed work
on high-order tensors.

Our approach is summarized as follows. Observe that due to the unitary
invariance of the sum-of-squares norm $\lVert\,\cdot\,\rVert_{F}$, the
orthonormal matrices $U,V,W$ in \eqref{eqT} and the orthonormal matrix $Q$ in
\eqref{eqS} are only determined up to an action of $\operatorname*{O}%
(p)\times\operatorname*{O}(q)\times\operatorname*{O}(r)$ and an action of
$\operatorname*{O}(r)$, respectively. We exploit this to our advantage by
reducing the problems to an optimization problem on a product of Grassmannians and a Grassmannian, respectively. Specifically, we reduce
\eqref{eqT}, a minimization problem over a product of three Stiefel manifolds
and a Euclidean space $\operatorname*{O}(l,p)\times\operatorname*{O}%
(m,q)\times\operatorname*{O}(n,r)\times\mathbb{R}^{p\times q\times r}$, to a
maximization problem over a product of three Grassmannians
$\operatorname*{Gr}(l,p)\times\operatorname*{Gr}(m,q)\times\operatorname*{Gr}%
(n,r)$; and likewise we reduce \eqref{eqS} from $\operatorname*{O}%
(n,r)\times\mathsf{S}^{3}(\mathbb{R}^{r})$ to a maximization problem over
$\operatorname*{Gr}(n,r)$. This reduction of \eqref{eqT} to product of
Grassmannians has been exploited in \cite{elsa09,ishteva09,ishteva09b}. 
The algorithms in \cite{elsa09,ishteva09,ishteva09b}
involve the Hessian, either explicitly or implicitly via its approximation on a tangent.
Whichever the case, the reliance on Hessian in these methods results in them quickly becoming infeasible as the size of the problem increases. With this in mind, we consider the quasi-Newton and limited memory quasi-Newton approaches described in the first paragraph of this section.

An important case not addressed in {\cite{elsa09,ishteva09, ishteva09b} }is
the multilinear approximation of symmetric tensors \eqref{eqS}. Note that the
general \eqref{eqT} and symmetric \eqref{eqS} cases are related but different,
not unlike the way the singular value problem differs from the symmetric eigenvalue
problem for matrices. The problem \eqref{eqT} for general tensors is linear in
the entries of $U,V,W$ (quadratic upon taking norm-squared) whereas the
problem \eqref{eqS} for symmetric tensors is cubic in the entries of $Q$
(sextic upon taking norm-squared). To the best of our knowledge, all existing
solvers for \eqref{eqS} are unsatisfactory because they rely on algorithms for
\eqref{eqT}. A typical heuristic is as follows: find \textit{three}
orthonormal matrices $Q_{1},Q_{2},Q_{3}$ and a nonsymmetric $\mathcal{C}%
^{\prime}\in\mathbb{R}^{r\times r\times r}$ that approximates~$\mathcal{S}$,%
\[
\mathcal{S}\approx(Q_{1},Q_{2},Q_{3})\cdot\mathcal{C}^{\prime},
\]
then artificially set $Q_{1}=Q_{2}=Q_{3}=Q$ by either averaging or choosing
the last iterate and then symmetrize $\mathcal{C}^{\prime}$. This of course is
not ideal. Furthermore, using the framework developed for the general tensor
approximation to solve the symmetric tensor approximation problem will be
computationally much more expensive. In particular, to optimize $\mathcal{S}%
\approx(Q_{1},Q_{2},Q_{3})\cdot\mathcal{C}^{\prime}$ without taking the
symmetry into account incurs a $k$-fold increase in computational cost
relative to $\mathcal{S}\approx(Q,\dots,Q)\cdot\mathcal{C}$. The algorithm
proposed in this paper solves \eqref{eqS} directly. It finds a \textit{single}
$Q\in\operatorname*{O}(n,r)$ and a \textit{symmetric} $\mathcal{C}%
\in\mathsf{S}^{3}(\mathbb{R}^{r})$ with%
\[
\mathcal{S}\approx(Q,Q,Q)\cdot\mathcal{C}.
\]
The symmetric case can often be more important than the general case, not
surprising since symmetric tensors are common in practice, arising as higher
order derivatives of smooth multivariate real-valued functions, higher order
moments and cumulants of a vector-valued random variable, etc.

Like the Grassmann-Newton algorithms in \cite{elsa09,ishteva09} and the
trust-region approach in \cite{ishteva09b,ishteva08}, the quasi-Newton
algorithms proposed in this article have guaranteed convergence to
stationary points and represent an improvement over Gauss-Seidel type
coordinate-cycling heuristics like alternating least squares (\textsc{als}),
higher-order orthogonal iteration (\textsc{hooi}), or higher-order singular
value decomposition (\textsc{hosvd}) \cite{delath00a}. And as far as
accuracy is concerned, our algorithms perform as well as the algorithms in
\cite{elsa09,ishteva08,ishteva09} and outperform Gauss-Seidel type strategies
in many cases. As far as robustness and speed are concerned, our algorithms
work on much larger problems and perform vastly faster than Grassmann-Newton
algorithm. Asymptotically the memory storage requirements of our algorithms
are of the same order-of-magnitude as Gauss-Seidel type strategies. For large
problems, our Grassmann \textsc{l-bfgs} algorithm outperforms even
Gauss-Seidel strategies (in this case \textsc{hooi}), which is not unexpected
since it has the advantage of requiring only a small number of prior iterates.

We will give the reader a rough idea of the performance of our algorithms.
Using \textsc{matlab} on a laptop computer, we attempted to find a solution to
an accuracy within machine precision, i.e.\ $\approx10^{-13}$. For general
$3$-tensors of size $200\times200\times200$, our Grassmann \textsc{l-bfgs}
algorithm took less than $13$ minute{s} while for general $4$-tensors of size
$50\times50\times50\times50$, it took about $7$ minutes. For symmetric
$3$-tensors of size $200\times200\times200$, our Grassmann \textsc{bfgs}
algorithm took about $5$ minutes while for symmetric $4$-tensors of size
$50\times50\times50\times50$, it took less than $2$ minutes. In all cases, we
{seek} a rank-$(5,5,5)$ or rank-$(5,5,5,5)$ approximation. For a general
order-$10$ tensor of dimensions $5\times\dots\times5$, a rank-$(2,\dots,2)$
approximation took about $15$ minutes to reach the same accuracy as above.
More extensive numerical experiments are reported in
Section~\ref{qng:sec:comEx}. The reader is welcomed to try our algorithms,
which have been made publicly available at \cite{savas08a, savas08b}.

\subsection{Outline}

The structure of the article is as follows. In Sections \ref{sec:Tensors} and
\ref{sec:mrank}, we present a more careful discussion of tensors, symmetric
tensors, multilinear rank, and their corresponding multilinear approximation
problems. In Section~\ref{qng:sec:optRM} we will discuss how quasi-Newton
methods in Euclidean space may be extended to Riemannian manifolds and, in
particular, Grassmannians. Section~\ref{qng:sec:geodParTran} contains a
discussion on geodesic curves and transport of vectors on Grassmannians. In
Section~\ref{qng:sec:QN} we present the modifications on quasi-Newton methods
with \textsc{bfgs} updates in order for them to be well-defined on
Grassmannians. Also, the reader will find proof of the optimality properties
of \textsc{bfgs} updates on product{s} of Grassmannians.
Section~\ref{qng:sec:lbfgs} gives the corresponding modifications for limited
memory \textsc{bfgs} updates. Section~\ref{qng:sec:brApp} states the
corresponding expressions for the tensor approximation problem, which are
defined on a product of Grassmannians. The symmetric case is detailed in
section~\ref{qng:sec:symmCase}. Section \ref{sec:ex} contains a few examples
with numerical calculations illustrating the presented concepts. The
implementation and the experimental results are found in
Section~\ref{qng:sec:comEx}. Related work and the conclusions are discussed in Section \ref{sec:relWork} and \ref{qng:sec:conclusion} respectively.

\subsection{Notations}

Tensors will be denoted by calligraphic letters, e.g.\ $\mathcal{A}$,
$\mathcal{B}$, $\mathcal{C}$. Matrices will be denoted in upper case letters,
e.g.\ $X$, $Y$, $Z$. We will also use upper case letters $X$, $X_{k}$ to
denote iterates or elements of a Grassmannian since we represent them as
(equivalence classes of) matrices with orthonormal columns. Vectors and
iterates in vector form are denoted with lower case letters, e.g.\ $x$, $y$,
$x_{k}$, $y_{k}$, where the subscript is the iteration index. To denote
scalars we use lower case Greek letters, e.g.\ $\alpha$, $\beta$, and $t$,
$t_{k}$.

We will use the usual symbol
$\mathbin{\scalebox{.88}{$\displaystyle\otimes$}}$ to denote the \textit{outer
product of tensors} and a large boldfaced version
$\mathbin{\scalebox{1.12}{$\displaystyle\varotimes$}}$ to denote the
\textit{Kronecker product of operators}. For example, if $A\in\mathbb{R}%
^{m\times n}$ and $B\in\mathbb{R}^{p\times q}$ are matrices, then
$A\mathbin{\scalebox{.88}{$\displaystyle\otimes$}}B$ will be a $4$-tensor in
$\mathbb{R}^{m\times n\times p\times q}$ whereas
$A\mathbin{\scalebox{1.12}{$\displaystyle\varotimes$}}B$ will be a matrix in
$\mathbb{R}^{mp\times nq}$. In the former case, we regard $A$ and $B$ as
$2$-tensors while in the latter case, we regard them as matrix representations
of linear operators. The contracted products used in this paper are
defined in Appendix \ref{app:1}.

For $r\leq n$, we will let $\operatorname*{O}(n,r)=\{X\in\mathbb{R}^{n\times
r}\mid X^{\mathsf{T}}X=I\}$ denote the \textit{Stiefel manifold} of $n\times
r$ matrices with orthonormal columns. The special case $r=n$, i.e.\ the
\textit{orthogonal group}, will be denoted $\operatorname*{O}(n)$. For $r\leq
n$, $\operatorname*{O}(r)$ acts on $\operatorname*{O}(n,r)$ via right
multiplication. The set of orbit classes $\operatorname*{O}%
(n,r)/\operatorname*{O}(r)$\ is a manifold called the \textit{Grassmann manifold} or \textit{Grassmannian}
(we adopt the latter name throughout this article)
and will be denoted $\operatorname*{Gr}(n,r)$.

In this paper, we will only cover a minimal number of notions and notations
required to describe our algorithm. Further mathematical details concerning
tensors, tensor ranks, tensor approximations, as well as the counterpart for
symmetric tensors may be found in \cite{comon08, dSL, elsa09}. Specifically we
will use the notational and analytical framework for tensor manipulations
introduced in \cite[Section~2]{elsa09} and assume that these concepts are
familiar to the reader.

\section{General and symmetric tensors\label{sec:Tensors}}

Let $V_{1},\dots,V_{k}$ be real vector spaces of dimensions $n_{1},\dots
,n_{k}$ respectively and let $\mathbf{A}$ be an element of the tensor product
$V_{1}\mathbin{\scalebox{.88}{$\displaystyle\otimes$}}\dots
\mathbin{\scalebox{.88}{$\displaystyle\otimes$}}V_{k}$, i.e.\ $\mathbf{A}$ is
a \textit{tensor} of order $k$ \cite{greub78, L, Y}. Up to a choice of bases
on $V_{1},\dots,V_{k}$, one may represent a tensor $\mathbf{A}$ as a
$k$-dimensional \textit{hypermatrix} $\mathcal{A}=[a_{i_{1}\cdots i_{k}}%
]\in\mathbb{R}^{n_{1}\times\dots\times n_{k}}$. Similarly, let $V$ be a real
vector space of dimension $n$ and $\mathbf{S}\in\mathsf{S}^{k}(V)$ be a
\textit{symmetric tensor} of order $k$ \cite{greub78, L, Y}. Up to a choice of
basis on $V$, $\mathbf{S}$ may be represented as a $k$-dimensional hypermatrix
$\mathcal{S}=[s_{i_{1}\cdots i_{k}}]\in\mathbb{R}^{n\times\dots\times n}$
whose entries are invariant under any permutation of indices, i.e.\
\begin{equation}
s_{i_{\sigma(1)}\cdots i_{\sigma(k)}}=s_{i_{1}\cdots i_{k}}\qquad\text{for
every }\sigma\in\mathfrak{S}_{k}\label{symm}%
\end{equation}
We will write $\mathsf{S}^{k}(\mathbb{R}^{n})$ for the subspace of
$\mathbb{R}^{n\times\dots\times n}$ satisfying \eqref{symm}. Henceforth, we
will assume that there are some predetermined bases and will not distinguish
between a tensor $\mathbf{A}\in V_{1}%
\mathbin{\scalebox{.88}{$\displaystyle\otimes$}}\dots
\mathbin{\scalebox{.88}{$\displaystyle\otimes$}}V_{k}$ and its hypermatrix
representation $\mathcal{A}\in\mathbb{R}^{n_{1}\times\dots\times n_{k}}$ and
likewise for a symmetric tensor $\mathbf{S}\in\mathsf{S}^{k}(V)$ and its
hypermatrix representation $\mathcal{S}\in\mathsf{S}^{k}(\mathbb{R}^{n})$.
Furthermore we will sometimes present our discussions for the case $k=3$ for
notational simplicity. We often call an order-$k$ tensor simply as a
$k$-tensor and an order-$k$ symmetric tensor as a symmetric $k$-tensor.

As we have mentioned in Section~\ref{qng:sec:intro}, symmetric tensors are
common in applications, largely because of the two examples below. The use of
higher-order statistics in signal processing and neuroscience, most notably
the technique of independent component analysis, symmetric tensors often play
a central role. The reader is referred to \cite{comon08} for further
discussion of symmetric tensors.

\begin{example}
Let $m\in\{1,2,\dots,\infty\}$ and $\Omega\subseteq\mathbb{R}^{n}$ be an open
subset. If $f\in C^{m}(\Omega)$, then for $k=1,\dots,m$, the $k$th derivative
of $f$ at $\mathbf{a}\in\Omega$ is a symmetric tensor of order $k$,%
\[
D^{k}f(\mathbf{a})=\left[  \frac{\partial^{k}f}{\partial x_{1}^{i_{1}}%
\dots\partial x_{n}^{i_{n}}}(\mathbf{a})\right]  _{i_{1}+\dots+i_{n}=k}%
\in\mathsf{S}^{k}(\mathbb{R}^{n}).
\]
For $k=1,2$, the vector $D^{1}f(\mathbf{a})$ and the matrix $D^{2}%
f(\mathbf{a})$ are the gradient and Hessian of $f$ at $\mathbf{a}$, respectively.
\end{example}

\begin{example}
\label{eg:Cum}Let $X_{1},\dots,X_{n}$ be random variables with respect to the
same probability distribution $\mu$. The moments and cumulants of the random
vector $\mathbf{X}=(X_{1},\dots,X_{n})$ are symmetric tensors of order $k$
defined by%
\begin{align*}
m_{k}(\mathbf{X)}  &  =\bigl[E(x_{i_{1}}x_{i_{2}}\cdots x_{i_{k}%
})\bigr]_{i_{1},\dots,i_{k}=1}^{n}\\
&  =\left[  \idotsint x_{i_{1}}x_{i_{2}}\cdots x_{i_{k}}~d\mu(x_{i_{1}})\cdots
d\mu(x_{i_{k}})\right]  _{i_{1},\dots,i_{k}=1}^{n}%
\end{align*}
and%
\[
\kappa_{k}(\mathbf{X})=\left[  \sum\nolimits_{\substack{A_{1}\sqcup\dots\sqcup
A_{p}\\=\{i_{1},\dots,i_{k}\}}}(-1)^{p-1}(p-1)!E({\textstyle\prod
\nolimits_{i\in A_{1}}}x_{i})\cdots E({\textstyle\prod\nolimits_{i\in A_{p}}%
}x_{i})\right]  _{i_{1},\dots,i_{k}=1}^{n}%
\]
respectively. The sum above is taken over all possible partitions
$\{i_{1},\dots,i_{k}\}=A_{1}\sqcup\dots\sqcup A_{p}$. It is not hard to show
that both $m_{k}(\mathbf{X)}$ and $\kappa_{k}(\mathbf{X)}\in\mathsf{S}%
^{k}(\mathbb{R}^{n})$. For $n=1$, the quantities $\kappa_{k}(\mathbf{X})$ for
$k=1,2,3,4$ have well-known names: they are the expectation, variance,
skewness, and kurtosis of the random variable $X$, respectively.
\end{example}

\section{Multilinear transformation and multilinear rank\label{sec:mrank}}

Matrices can act on other matrices through two independent multiplication
operations: left-multiplication and right-multiplication. If $C\in
\mathbb{R}^{p\times q}$ and $X\in\mathbb{R}^{m\times p}$, $Y\in\mathbb{R}%
^{n\times q}$, then the matrix~$C$ may be transformed into the matrix
$A\in\mathbb{R}^{m\times n}$, by
\begin{equation}
A=XCY^{\mathsf{T}},\qquad a_{ij}=\sum_{\alpha=1}^{p}\sum_{\beta=1}%
^{q}x_{i\alpha}y_{j\beta}c_{\alpha\beta}\,. \label{eq:lrmm}%
\end{equation}
Matrices act on order-$3$ tensors via \textit{three} different multiplication
operations. As in the matrix case, these can be combined into a single
formula. If $\mathcal{C}\in\mathbb{R}^{p\times q\times r}$ and $X\in
\mathbb{R}^{l\times p}$, $Y\in\mathbb{R}^{m\times q}$, $Z\in\mathbb{R}%
^{n\times r}$, then the $3$-tensor~$\mathcal{C}$ may be transformed into the
$3$-tensor $\mathcal{A}\in\mathbb{R}^{l\times m\times n}$ via
\begin{equation}
a_{ijk}=\sum_{\alpha=1}^{p}\sum_{\beta=1}^{q}\sum_{\gamma=1}^{r}x_{i\alpha
}y_{j\beta}z_{k\gamma}c_{\alpha\beta\gamma}\,. \label{eq:tmm}%
\end{equation}
We call this operation the \textit{trilinear multiplication} of~$\mathcal{A}$
by matrices $X$, $Y$ and $Z$, which we write succinctly as
\begin{equation}
\mathcal{A}=(X,Y,Z)\cdot{\mathcal{C}}. \label{eq:tmm1}%
\end{equation}
This is nothing more than the trilinear equivalent of \eqref{eq:lrmm}, which
in this notation has the form
\[
A=XCY^{\mathsf{T}}=(X,Y)\cdot{C}.
\]
Informally, \eqref{eq:tmm} amounts to multiplying the $3$-tensor $\mathcal{A}$
on its three `sides' or modes by the matrices $X$, $Y$, and $Z$ respectively.

An alternative but equivalent way of writing \eqref{eq:tmm} is as follows.
Define the \textit{outer product} of vectors $x\in\mathbb{R}^{l}$,
$y\in\mathbb{R}^{m}$, $z\in\mathbb{R}^{n}$ by
\[
x\mathbin{\scalebox{.88}{$\displaystyle\otimes$}}
y\mathbin{\scalebox{.88}{$\displaystyle\otimes$}} z=[x_{i}y_{j}z_{k}%
]\in\mathbb{R}^{l\times m\times n},
\]
and call a tensor of the form
$x\mathbin{\scalebox{.88}{$\displaystyle\otimes$}}
y\mathbin{\scalebox{.88}{$\displaystyle\otimes$}} z$ a \textit{decomposable
tensor} or, if non-zero, a \textit{rank-}$1$\textit{ tensor}. One may also
view \eqref{eq:tmm} as a \textit{trilinear combination} (as opposed to a
linear combination) of the decomposable tensors given by%
\begin{equation}
\mathcal{A}=\sum_{\alpha=1}^{p}\sum_{\beta=1}^{q}\sum_{\gamma=1}^{r }%
c_{\alpha\beta\gamma}x_{\alpha}%
\mathbin{\scalebox{.88}{$\displaystyle\otimes$}} y_{\beta}%
\mathbin{\scalebox{.88}{$\displaystyle\otimes$}} z_{\gamma}.
\label{eq:trilinear}%
\end{equation}
The vectors $x_{1},\dots,x_{p}\in\mathbb{R}^{l}$, $y_{1},\dots,y_{q}%
\in\mathbb{R}^{m}$, $z_{1},\dots,z_{r}\in\mathbb{R}^{n}$ are, of course, the
column vectors of the respective matrices $X,Y,Z$ above. In this article, we
find it more natural to present our algorithms in the form \eqref{eq:tmm1} and
therefore we refrain from using \eqref{eq:trilinear}.

More abstractly, given linear transformations of real vector spaces
$T_{u}:U\rightarrow U^{\prime}$, $T_{v}:V\rightarrow V^{\prime}$,
$T_{w}:W\rightarrow W^{\prime}$, the functoriality of tensor product
\cite{greub78, L} (denoted by the usual notation
$\mathbin{\scalebox{.88}{$\displaystyle\otimes$}} $ below) implies that one
has an induced linear transformation between the tensor product of the
respective vector spaces%
\[
T_{u}\mathbin{\scalebox{.88}{$\displaystyle\otimes$}} T_{v}
\mathbin{\scalebox{.88}{$\displaystyle\otimes$}} T_{w} :U
\mathbin{\scalebox{.88}{$\displaystyle\otimes$}} V
\mathbin{\scalebox{.88}{$\displaystyle\otimes$}} W\rightarrow U^{\prime}
\mathbin{\scalebox{.88}{$\displaystyle\otimes$}} V^{\prime}
\mathbin{\scalebox{.88}{$\displaystyle\otimes$}} W^{\prime}.
\]
The trilinear matrix multiplication above is a coordinatized version of this
abstract transformation. In the special case where $U=U^{\prime}$,
$V=V^{\prime}$, $W=W^{\prime}$ and $T_{u}$, $T_{v}$, $T_{w}$ are (invertible)
change-of-basis transformations, the operation in \eqref{eq:tmm} describes the
manner a (contravariant) $3$-tensor transforms under change-of-coordinates of
$U$, $V$ and $W$.

Associated with the trilinear matrix multiplication \eqref{eq:tmm1} is the
following notion of tensor rank that generalizes the row-rank and column-rank
of a matrix. Let $\mathcal{A}\in\mathbb{R}^{l\times m\times n}$. For fixed
values of $j\in\{1,\dots,m\}$ and $k\in\{1,\dots,n\}$, consider the column
vector, written in a \textsc{matlab}-like notation, $\mathcal{A}%
(:,j,k)\in\mathbb{R}^{l}$. Likewise we may consider $\mathcal{A}%
(i,:,k)\in\mathbb{R}^{m}$ for fixed values of $i,k$, and $\mathcal{A}%
(i,j,:)\in\mathbb{R}^{n}$ for fixed values of $i,j$. Define
\begin{align*}
r_{1}(\mathcal{A})  &  :=\dim(\operatorname{span}\{\mathcal{A}(:,j,k)\mid1\leq
j\leq m,1\leq k\leq n\}),\\
r_{2}(\mathcal{A})  &  :=\dim(\operatorname{span}\{\mathcal{A}(i,:,k)\mid1\leq
i\leq l,1\leq k\leq n\}),\\
r_{3}(\mathcal{A})  &  :=\dim(\operatorname{span}\{\mathcal{A}(i,j,:)\mid1\leq
i\leq l,1\leq j\leq m\}).
\end{align*}
Note that $\mathbb{R}^{l\times m\times n}$ may also be viewed as
$\mathbb{R}^{l\times mn}$. Then $r_{1}(\mathcal{A})$ is simply the rank
of~$\mathcal{A}$ regarded as an $l\times mn$ matrix, see the discussion on
tensor matricization in \cite{elsa09} with similar interpretations for
$r_{2}(\mathcal{A})$ and $r_{3}(\mathcal{A})$. The \textit{multilinear rank}
of $\mathcal{A}$ is the $3$-tuple $(r_{1}(\mathcal{A}),r_{2}(\mathcal{A}%
),r_{3}(\mathcal{A}))$ and we write
\[
\operatorname*{rank}(\mathcal{A})=(r_{1}(\mathcal{A}),r_{2}(\mathcal{A}%
),r_{3}(\mathcal{A})).
\]
We need to `store' all three numbers as $r_{1}(\mathcal{A})\neq r_{2}%
(\mathcal{A})\neq r_{3}(\mathcal{A})$ in general---a clear departure from the
case of matrices where row-rank and column-rank are always equal. Note that a
rank-$1$ tensor must necessarily have multilinear rank $(1,1,1)$,
i.e.$\ \operatorname*{rank}(x\mathbin{\scalebox{.88}{$\displaystyle\otimes$}}
y\mathbin{\scalebox{.88}{$\displaystyle\otimes$}} z)=(1,1,1)$ if $x,y,z$ are
non-zero vectors.

{For symmetric} tensors, one would be interested in transformation that
preserves the symmetry. For a symmetric matrix $C\in\mathsf{S}^{2}%
(\mathbb{R}^{n})$, this would be%
\begin{equation}
S=XCX^{\mathsf{T}},\qquad s_{ij}=\sum_{\alpha=1}^{r}\sum_{\beta=1}%
^{r}x_{i\alpha}x_{j\beta}c_{\alpha\beta}\,.\label{eq:smm}%
\end{equation}
Matrices act on symmetric order-$3$ tensors $\mathcal{S}\in\mathsf{S}%
^{3}(\mathbb{R}^{n})$ via the symmetric version of~\eqref{eq:tuck1}
\begin{equation}
s_{ijk}=\sum_{\alpha=1}^{r}\sum_{\beta=1}^{r}\sum_{\gamma=1}^{r}x_{i\alpha
}x_{j\beta}x_{k\gamma}c_{\alpha\beta\gamma}\,.\label{eq:stmm}%
\end{equation}
We call this operation the \textit{symmetric trilinear multiplication}
of~$\mathcal{S}$ by matrix $X$, which in the notation above, is written as
\[
\mathcal{S}=(X,X,X)\cdot{\mathcal{C}}.
\]
This is nothing more than the cubic equivalent of \eqref{eq:smm}, which in
this notation becomes
\[
S=XCX^{\mathsf{T}}=(X,X)\cdot{C}.
\]
Informally, \eqref{eq:stmm} amounts to multiplying the $3$-tensor
$\mathcal{A}$ on its three `sides' or modes by the same matrix $X$. In the
multilinear combination form, this is%
\[
\mathcal{S}=\sum_{\alpha=1}^{r}\sum_{\beta=1}^{r}\sum_{\gamma=1}^{r}%
c_{\alpha\beta\gamma}x_{\alpha}%
\mathbin{\scalebox{.88}{$\displaystyle\otimes$}}x_{\beta}%
\mathbin{\scalebox{.88}{$\displaystyle\otimes$}}x_{\gamma}%
\]
where the vectors $x_{1},\dots,x_{r}\in\mathbb{R}^{n}$ are the column vectors
of the matrix $X$ above.

More abstractly, given a linear transformation of real vector spaces
$T:V\rightarrow V^{\prime}$, the functoriality of symmetric tensor product
\cite{greub78, L} implies that one has an induced linear transformation
between the tensor product of the respective vector spaces%
\[
\mathsf{S}^{k}(T):\mathsf{S}^{k}(V)\rightarrow\mathsf{S}^{k}(V^{\prime}).
\]
The symmetric trilinear matrix multiplication above is a coordinatized version
of this abstract transformation. In the special case where $V=V^{\prime}$ and
$T$ is an (invertible) change-of-basis transformation, the operation in
\eqref{eq:stmm} describes the manner a symmetric $3$-tensor transforms under
change-of-coordinates of $V$.

For a symmetric tensor $\mathcal{S}\in\mathsf{S}^{3}(\mathbb{R}^{n})$, we must
have%
\[
r_{1}(\mathcal{S})=r_{2}(\mathcal{S})=r_{3}(\mathcal{S})
\]
by its symmetry. See Lemma~\ref{lem:symrank} for a short proof. The
\textit{symmetric multilinear rank} of $\mathcal{S}$ is the common value,
denoted $r_{\mathsf{S}}(\mathcal{S})$. When referring to a symmetric tensor in
this article, rank would always mean symmetric multilinear rank; e.g.\ a
rank-$s$ symmetric tensor $\mathcal{S}$ would be one with $r_{\mathsf{S}%
}(\mathcal{S})=s$.

We note that there is a different notion of tensor rank and symmetric tensor
rank, defined as the number of terms in a minimal decomposition of a tensor
(resp.\ symmetric tensor) into rank-$1$ tensors (resp.\ rank-$1$ symmetric
tensor{s}). Associated with this notion of rank are low-rank approximation
problems for tensors and symmetric tensors analogous to the ones discussed in
this paper. Unfortunately, these are ill-posed problems that may not even have
a solution \cite{comon08, dSL}. As such, we will not discuss this other notion
of tensor rank. It is implicitly assumed that whenever we discuss tensor rank
or symmetric tensor rank, it is with the multilinear rank or symmetric
multilinear rank defined above in mind.

\subsection{Multilinear approximation as maximization over
Grassmannians\label{sec:approx-max}}

Let $\mathcal{A}\in\mathbb{R}^{l\times m\times n}$ be a given third order
tensor and consider the problem
\[
\min\{\lVert\mathcal{A}-\mathcal{B}\rVert_{F}\mid\operatorname*{rank}%
(\mathcal{B})\leq(p,q,r)\}.
\]
Under this rank constraint, we can write $\mathcal{B}$ in factorized form
\begin{equation}
\mathcal{B}=(X,Y,Z)\cdot{\mathcal{C}},\quad b_{ijk}=\sum_{\alpha=1}^{p}%
\sum_{\beta=1}^{q}\sum_{\gamma=1}^{r}x_{i\alpha}y_{j\beta}z_{k\gamma}%
c_{\alpha\beta\gamma},\label{eq:tuck1}%
\end{equation}
where $\mathcal{C}\in\mathbb{R}^{p\times q\times r}$ and $X\in\mathbb{R}%
^{l\times p}$, $Y\in\mathbb{R}^{m\times q}$, $Z\in\mathbb{R}^{n\times r}$, are
full rank matrices. In other words, one would like to solve the following
\textit{best multilinear rank approximation problem},
\begin{equation}
\min_{X,Y,Z,\mathcal{C}}\lVert\mathcal{A}-(X,Y,Z)\cdot{\mathcal{C}}\rVert
_{F}.\label{eq:bmlra}%
\end{equation}
This is the optimization problem underlying the \textit{Tucker model}
\cite{tucke66} that originated in psychometrics but has become increasingly
popular in other areas of data analysis. In fact, there is no loss of
generality if we assume $X^{\mathsf{T}}X=I$, $Y^{\mathsf{T}}Y=I$ and
$Z^{\mathsf{T}}Z=I$. This is verified as follows, for any full column-rank
matrices $\widetilde{X}$, $\widetilde{Y}$ and $\widetilde{Z}$ we can compute
their \textsc{qr}-factorizations
\[
\widetilde{X}=XR_{X},\quad\widetilde{Y}=YR_{Y},\quad\widetilde{Z}=ZR_{Z}%
\]
and multiply the right triangular matrices into the core tensor, i.e.\
\[
\lVert\mathcal{A}-(\widetilde{X},\widetilde{Y},\widetilde{Z})\cdot
{\widetilde{\mathcal{C}}\rVert_{F}}=\lVert\mathcal{A}-(X,Y,Z)\cdot
{\mathcal{C}}\rVert_{F}\quad\text{where}\quad\mathcal{C}=(R_{X},R_{Y}%
,R_{Z})\cdot{\widetilde{\mathcal{C}}}\,.
\]
With the orthonormal constraints on $X$ ,$Y$ and $Z$ the tensor approximation
problem can be viewed as an optimization problem on a product of Stiefel
manifolds. Using the identity
\[
(U^{\mathsf{T}},V^{\mathsf{T}},W^{\mathsf{T}})\cdot\mathcal{A}\equiv
{\mathcal{A}}\cdot(U,V,W),
\]
we can rewrite the tensor approximation problem as a maximization problem with
the objective function
\begin{equation}
\Phi(X,Y,Z)=\frac{1}{2}\left\Vert \mathcal{A}\cdot(X,Y,Z)\right\Vert _{F}%
^{2}\quad\text{s.t.}\quad X^{\mathsf{T}}X=I,\quad Y^{\mathsf{T}}Y=I,\quad
Z^{\mathsf{T}}Z=I,\label{eq:GrassProd}%
\end{equation}
in which the small core tensor $\mathcal{C}$ is no longer present.
See references \cite{latha00a, elsa09} for the elimination of $\mathcal{C}$.
The objective function $\Phi(X,Y,Z)$ is invariant under orthogonal
transformation of the variable matrices $X,Y,Z$ from the right. Specifically, for
any orthogonal matrices $Q_{1}\in\operatorname*{O}(p)$, $Q_{2}\in
\operatorname*{O}(q)$ and $Q_{3}\in\operatorname*{O}(r)$ it holds that
$\Phi(X,Y,Z)=\Phi(XQ_{1},YQ_{2},ZQ_{3})$. This homogeneity property implies
that $\Phi(X,Y,Z)$ is in fact defined on a product of three Grassmannians $\operatorname*{Gr}(l,p)\times\operatorname*{Gr}(m,q)\times
\operatorname*{Gr}(n,r)$.

Let $\mathcal{S}\in\mathsf{S}^{3}(\mathbb{R}^{n})$ be a given symmetric
$3$-tensor and consider the problem
\[
\min\{\lVert\mathcal{S}-\mathcal{T}\rVert_{F}\mid\mathcal{T}\in\mathsf{S}%
^{3}(\mathbb{R}^{n}),\;r_{\mathsf{S}}(\mathcal{S})\leq r\}.
\]
Under this rank constraint, we can write $\mathcal{T}$ in factorized form
\[
\mathcal{T}=(X,X,X)\cdot{\mathcal{C}},\quad t_{ijk}=\sum_{\alpha=1}^{r}%
\sum_{\beta=1}^{r}\sum_{\gamma=1}^{r}x_{i\alpha}x_{j\beta}x_{k\gamma}%
c_{\alpha\beta\gamma},
\]
where $\mathcal{C}\in\mathsf{S}^{3}(\mathbb{R}^{r})$ and $X\in\mathbb{R}%
^{n\times r}$ is a full-rank matrix. In other words, one would like to solve
the following \textit{best symmetric multilinear rank approximation problem},
\[
\min\{\lVert\mathcal{A}-(X,X,X)\cdot{\mathcal{C}}\rVert_{F}\mid X\in
\operatorname*{O}(n,r),\;\mathcal{C}\in\mathsf{S}^{3}(\mathbb{R}^{r})\}.
\]
As with the general case, there is no loss of generality if we assume
$X^{\mathsf{T}}X=I$. With the orthonormal constraints on $X$, the tensor
approximation problem can be viewed as an optimization problem on a single
Stiefel manifold (as opposed to a product of Stiefel manifolds in
\eqref{eq:GrassProd}). Using the identity
\[
(V^{\mathsf{T}},V^{\mathsf{T}},V^{\mathsf{T}})\cdot\mathcal{S}\equiv
{\mathcal{S}}\cdot(V,V,V),
\]
we may again rewrite the tensor approximation problem as a maximization
problem with the objective function
\[
\Phi(X)=\frac{1}{2}\left\Vert \mathcal{S}\cdot(X,X,X)\right\Vert _{F}^{2}%
\quad\text{s.t.}\quad X^{\mathsf{T}}X=I,
\]
in which the core-tensor $\mathcal{C}$ is no longer present. As with
the general case, the objective function $\Phi(X)$ also has an invariance
property, namely $\Phi(X)=\Phi(XQ)$ for any orthogonal $Q\in\operatorname*{O}%
(r)$. As before, this homogeneity property implies that $\Phi(X)$ is
well-defined on a single Grassmannian $\operatorname*{Gr}(n,r)$.

These multilinear approximation problems may be viewed as `dimension
reduction' or `rank reduction' for tensors and symmetric tensors respectively.
In general, a matrix requires $O(n^{2})$ storage and an order-$k$ tensor
requires $O(n^{k})$ storage. While it is sometimes important to perform
dimension reduction to a matrix, a dimension reduction is almost always
necessary if one wants to work effectively with a tensor of higher order. A
dimension reduction of a matrix $A\in\mathbb{R}^{n\times n}$ of the form
$A\approx UCV^{\mathsf{T}}$, where $U,V\in\operatorname*{O}(n,r)$ and diagonal
$C\in\mathbb{R}^{r\times r}$ reduces dimension from $O(n^{2})$ to
$O(nr+r)$. A dimension reduction of a tensor $\mathcal{A}\in
\mathbb{R}^{n\times\dots\times n}$ of the form $\mathcal{A}\approx(Q_{1}%
,\dots,Q_{k})\cdot\mathcal{C}$ where $Q_{1},\dots,Q_{k}\in\operatorname*{O}%
(n,r)$ and $\mathcal{C}\in\mathbb{R}^{r\times\dots\times r}$, reduces
dimension from $O(n^{k})$ to $O(knr+r^{k})$. If $r$ is significantly smaller
than $n$, e.g.\ $r=O(n^{1/k})$, then a dimension reduction in the higher order
case could reduce problem size by orders of magnitude.

\section{Optimization in Euclidean space and on Riemannian
manifolds\label{qng:sec:optRM}}

In this section we discuss the necessary modifications for generalizing an
optimization algorithm from Euclidean space to Riemannian manifolds.
Specifically we consider the quasi-Newton methods with \textsc{bfgs} and
\textit{limited memory} \textsc{bfgs} (\textsc{l-bfgs}) updates. First we
state the expressions in Euclidean space and then we point out what needs to
be modified. The convergence properties of quasi-Newton methods defined on
manifolds were established by Gabay~\cite{gabay82}. Numerical treatment of
algorithms on the Grassmannian are given in
\cite{edelm99,absil07,simon07,lunds02,ishteva08}. A recent book on
optimization on manifolds is~\cite{absil08}. In this and the next three
sections, i.e.\ Sections~\ref{qng:sec:optRM} through \ref{qng:sec:lbfgs}, we
will discuss our algorithms in the context of minimization problems, as is
conventional. It will of course be trivial to modify them for maximization
problem. Indeed the tensor approximation problems discussed in Sections
\ref{qng:sec:brApp} through \ref{qng:sec:comEx} will all be solved as
maximization problems on Grassmannians or product of Grassmannians.

\subsection{BFGS updates in Euclidean Space\label{qng:sec:BFGS-EucSpace}}

Assume that we want to minimize a nonlinear real valued function $f(x)$ where
$x\in\mathbb{R}^{n}$. As is well-known, in quasi-Newton methods, one solves%
\begin{equation}
H_{k}p_{k}=-g_{k}, \label{eq:QN}%
\end{equation}
to obtain the direction of descent $p_{k}$ from the current iterate $x_{k}$
and the gradient $g_{k}=\nabla f(x_{k})$ at $x_{k}$. Unlike Newton method,
which uses the exact Hessian for $H_{k}$, in \eqref{eq:QN} $H_{k}$ is only an
approximation of the Hessian at $x_{k}$. After computing the (search)
direction $p_{k}$ one obtains the next iterate as $x_{k+1}=x_{k}+t_{k}p_{k}$
in which the step length $t_{k}$ is usually given by a line search method
satisfying the Wolfe or the Goldstein conditions \cite{nowr:06}. Instead of
recomputing the Hessian at each new iterate $x_{k+1}$, it is updated from the
previous approximation. The \textsc{bfgs} update has the following form,
\begin{equation}
H_{k+1}=H_{k}-\frac{H_{k}s_{k}s_{k}^{\mathsf{T}}H_{k}}{s_{k}^{\mathsf{T}}%
H_{k}s_{k}}+\frac{y_{k}y_{k}^{\mathsf{T}}}{y_{k}^{\mathsf{T}}s_{k}},
\label{eq:BFGS-update}%
\end{equation}
where
\begin{align}
s_{k}  &  =x_{k+1}-x_{k}=t_{k}p_{k},\label{eq:eucsk}\\
y_{k}  &  =g_{k+1}-g_{k}. \label{eq:eucyk}%
\end{align}
Quasi-Newton methods with \textsc{bfgs} updates are considered to be the most
computationally efficient algorithms for minimization of general nonlinear
functions. This efficiency is obtained by computing a new Hessian
approximation as a rank-$2$ modification of the previous Hessian. The
convergence of quasi-Newton methods is super-linear in a vicinity of a local
minimum. In most cases the quadratic convergence of the Newton method is
outperformed by quasi-Newton methods since each iteration is computationally
much cheaper than a Newton iteration. A thorough study of quasi-Newton methods may be found in \cite{nowr:06,denni96}. The reader is reminded that quasi-Newton methods do not necessarily converge to local minima but only to stationary points.
Nevertheless, using the Hessians that we derived in Sections~\ref{sec:genDeriv} and \ref{sec:genDerivSym} (for the general and symmetric cases respectively), the nature of these stationary points can often be determined.

\subsection{Quasi-Newton methods on a Riemannian manifold\label{eq:InfMan}}

We will give a very brief overview of Riemannian geometry tailored
specifically to our needs in this paper. First we sketch the modifications
that are needed in order for an optimization {algorithm} to be well-defined
when the objective function is defined on a manifold. For details and proof,
the reader should refer to standard literature on differential and Riemannian
geometry \cite{conlo01,booth86,absil08}. Informally a manifold, denoted with
$M$, is an object locally homeomorphic to $\mathbb{R}^{n}$. We will regard $M$
as a submanifold of some high-dimensional ambient Euclidean space
$\mathbb{R}^{N}$. Our objective function $f:M\rightarrow\mathbb{R}$ will be
assumed to have (at least) continuous second order partial derivatives. We
will write $f(x)$, $g(x)$, and $H(x)$ for the value of the function, the
gradient, and the Hessian at $x\in M$. These will be reviewed in greater
detail in the following.

Equations \eqref{eq:QN}--\eqref{eq:eucyk} are the basis of any algorithmic
implementation involving \textsc{bfgs} or \textsc{l-bfgs} updates. The key
operations are (1) computation of the gradient, (2) computation of the Hessian
or its approximation, (3) subtraction of iterates, e.g.\ to get $s_{k}$ or
$x_{k+1}$ and (4) subtraction of gradients. Each of these points needs to be
modified in order for these operations to be well-defined on manifolds.

\paragraph{Computation of the gradient}

The gradient $g(x)=\operatorname{grad}f(x)=\nabla f(x)$ at $x\in M$ of a
real-valued function defined on a manifold, $f:M\rightarrow\mathbb{R}$, $M\ni
x\mapsto f(x)\in\mathbb{R}$, is a vector in the tangent space $\mathbf{T}%
_{x}=\mathbf{T}_{x}(M)$ of the manifold at the given point $x$. We write
$g\in\mathbf{T}_{x}$.

To facilitate computations, we will often embed our manifold $M$ in some
ambient Euclidean space $\mathbb{R}^{N}$ and in turn endow $M$ with a system
of global coordinates $(x_{1},\dots,x_{N})$ where $N$ is usually larger than
$d:=\dim(M)$, the intrinsic dimension of $M$. The function $f$ then inherits
an expression in terms of $x_{1},\dots,x_{N}$, say, $f(x)=\tilde{f}%
(x_{1},\dots,x_{N})$. We would like to caution our readers that computing
$\nabla f$ on $M$ is not a matter of simply taking partial derivatives of
$\tilde{f}$ with respect to $x_{1},\dots,x_{N}$. An easy way to observe this
is that $\nabla f$ is a $d$-tuple when expressed in local coordinates whereas
$\left(  \partial\tilde{f}/\partial x_{1},\dots,\partial\tilde{f}/\partial
x_{N}\right)  $ is an $N$-tuple.

\paragraph{Computation of the Hessian or its approximation}

The Hessian $H(x)=\operatorname*{Hess}f(x)=\nabla^{2}f(x)$ at $x\in M$ of a
function $f:M\rightarrow\mathbb{R}$ is a linear transformation of the tangent
space $\mathbf{T}_{x}$ to itself, i.e.\
\begin{equation}
H(x):\mathbf{T}_{x}\rightarrow\mathbf{T}_{x}. \label{eq:H}%
\end{equation}
As in the case of gradient, when $f$ is expressed in terms of global
coordinates, differentiating the expression twice will in general not give the
correct Hessian.

\paragraph{Updating the current iterate}

Given an iterate $x_{k}\in M$, a step length $t_{k}$ and a search direction
$p_{k}\in\mathbf{T}_{x_{k}}$ the update $x_{k+1}=x_{k}+t_{k}p_{k}$ will in
general not be a point on the manifold. The corresponding operation on a
manifold is to move along the geodesic curve of the manifold given by the
direction $p_{k}$. Geodesics on manifolds correspond to straight lines in
Euclidean spaces. The operation $s_{k}=x_{k+1}-x_{k}$ is undefined in general
when the points on the right hand side belong to a general manifold.

\paragraph{Updating vectors and operators}

The quantities $s_{k}$ and $y_{k}$ are in fact tangent vectors and the
Hessians (or Hessian approximations) are linear operators all defined at a
specific point of the manifold. A given Hessian $H(x)$ is only defined at a
point $x\in M$ and correspondingly only acts on vectors in the tangent space
$\mathbf{T}_{x}$. In the right hand side of the \textsc{bfgs} update
\eqref{eq:BFGS-update} all quantities need to be defined at the same point
$x_{k}$ in order to have well-defined operations between the terms. In
addition the resulting sum will define an operator at a new point $x_{k+1}$.
The notion of parallel transporting vectors along geodesic curves resolves all
of these issues. The operations to the Hessian are similar and involves
parallel transport operation back and forth between two different points.

\section{Grassmann geodesics and parallel transport of
vectors\label{qng:sec:geodParTran}}

In this paper, the Riemannian manifold $M$ of most interest to us is
$\operatorname*{Gr}(n,r)$, the Grassmannian or Grassmannian of
$r$-planes in $\mathbb{R}^{n}$. Our discussion will proceed with
$M=\operatorname*{Gr}(n,r)$.

\subsection{Algorithmic considerations}

Representing points on Grassmannians as (equivalence classes of)
matrices allow{s} us to take advantage of matrix arithmetic and matrix
algorithms, as well as the readily available libraries of highly optimized and
robust matrix computational softwares developed over the last five decades
\cite{lapack99, arpack98}. A major observation of \cite{edelm99} is the
realization that common differential geometric operations on points of
Grassmann and Stiefel manifolds can all be represented in terms of matrix
operations. For our purposes, the two most important operations are (1) the
determination of a geodesic at a point along a tangent vector and (2) the
parallel transport of a tangent vector along a geodesic. On a product of
Grassmannians, these operations may likewise be represented in terms of
matrix operations \cite{elsa09}. We will give explicit expressions for
geodesic curves on Grassmannians and two different ways of parallel
transporting tangent vectors.

\subsection{Grassmannians in terms of matrices}

First we will review some preliminary materials from \cite{edelm99}. A point
$\langle X \rangle$ on the Grassmannian $\operatorname*{Gr}(n,r)$ is an
equivalence class of orthonormal matrices whose columns form an orthonormal
basis for an $r$-dimensional subspace of $\mathbb{R}^{n}$. Explicitly, we
write
\[
\langle X \rangle=\{XQ\in\operatorname*{O}(n,r)\mid Q\in\operatorname*{O}%
(r)\},
\]
where $X\in\operatorname*{O}(n,r)$, i.e.\ $X$ is an $n\times r$ matrix and
$X^{\mathsf{T}}X=I_{r}$. The set $\operatorname*{O}(n,r)=\{X\in\mathbb{R}%
^{n\times r}\mid X^{\mathsf{T}}X=I_{r}\}$ is also a manifold, often called the
Stiefel manifold. When $n=r$, $\operatorname*{O}(r,r)=\operatorname*{O}(r)$ is
the orthogonal group. It is easy to see that the dimensions of these manifolds
are%
\[
\dim(\operatorname*{O}(r))=\frac{1}{2}r(r-1),\quad\dim(\operatorname*{O}%
(n,r))=nr-\frac{1}{2}r(r+1),\quad\dim(\operatorname*{Gr}(n,r))=r(n-r).
\]
In order to use standard linear algebra in our computations, we will not be
able to work with a whole equivalence class of matrices. So by a point on a
Grassmannian, we will always mean some $X\in\langle X\rangle$ that represents
the equivalence class. The functions $f$ that we optimize in this paper will
take orthonormal matrices in $\operatorname*{O}(n,r)$ as arguments but will
always be well-defined on Grassmannians: given two representatives
$X_{1},X_{2}\in\langle X\rangle$, we will have $f(X_{1})=f(X_{2})$, i.e.%
\[
f(XQ)=f(X)\quad\text{for every }Q\in\operatorname*{O}(r).
\]
Abusing notations slightly, we will sometimes write $f:\operatorname*{Gr}%
(n,r)\rightarrow\mathbb{R}$.

The tangent space $\mathbf{T}_{X}$, where $X\in\operatorname*{Gr}(n,r)$, is an
affine vector space with elements in $\mathbb{R}^{n\times r}$. It can be shown
that any element $\Delta\in\mathbf{T}_{X}$ satisfies
\[
X^{\mathsf{T}}\Delta=0.
\]
The projection on the tangent space is
\[
\Pi_{X}=I-XX^{\mathsf{T}}=X_{\perp}X_{\perp}^{\mathsf{T}},
\]
where $X_{\perp}$ is an orthogonal complement of $X$, i.e.\ the square matrix
$[X\,X_{\perp}]$ is an $n\times n$ orthogonal matrix. Since by definition
$X^{\mathsf{T}}X_{\perp}=0$, any tangent vector can also be written as
$\Delta=X_{\perp}D$ where $D$ is an $(n-r)\times r$ matrix. This shows that
the columns of $X_{\perp}$ may be interpreted as a basis for $\mathbf{T}_{X}$.
We say that $\Delta$ is a \textit{global coordinate representation} and $D$ is
a \textit{local coordinate representation} of the same tangent. Note that the
number of degrees of freedom in $D$ equals the dimension of the tangent space
$\mathbf{T}_{X}$, which is $r(n-r)$. It follows that for a given tangent in
global coordinates $\Delta$, its local coordinate representation is given by
$D=X_{\perp}^{\mathsf{T}}\Delta$. Observe that to a given local representation
$D$ of a tangent there is an associated basis matrix $X_{\perp}$. Tangent
vectors are also embedded in $\mathbb{R}^{nr}$ since in global coordinates
they are given by $n\times r$ matrices. We will define algorithms using both
global coordinates as well as intrinsic local coordinates. When using global
coordinates, the Grassmannian $\operatorname*{Gr}(n,r)$ is
(isometrically) embed{ded} in the Euclidean space $\mathbb{R}^{n\times r}$ and
a product of Grassmannians in a corresponding product of Euclidean
spaces. The use of Pl\"{u}cker coordinates to represent points on Grassmannian is not useful for our purpose.

\subsection{Geodesics}

Let $X\in\operatorname*{Gr}(n,r)$ and $\Delta$ be a tangent vector at $X$,
i.e.\ $\Delta\in\mathbf{T}_{X}$. The geodesic path from $X$ in the direction
$\Delta$ is given by
\begin{equation}
X(t)=[XV\,\,\,\,U]%
\begin{bmatrix}
\cos\Sigma t\\
\sin\Sigma t
\end{bmatrix}
V^{\mathsf{T}}, \label{eq:geodesic}%
\end{equation}
where $\Delta=U\Sigma V^{\mathsf{T}}$ is the thin \textsc{svd} and we identify
$X(0)\equiv X$. Observe that omitting the last $V$ in \eqref{eq:geodesic} will
give the same path on the manifold but with a different\footnote{A given
matrix representation of a point on a Grassmannian can be
postmultiplied by any orthogonal matrix, giving a new representation of the same point.} representation. This information is
useful because some algorithms require a consistency in the matrix
representations along a path but other algorithms do not. For example, in a
Newton-Grassmann algorithm we may omit the second $V$ \cite{elsa09} but in
quasi-Newton-Grassmann algorithms $V$ is necessary.

\subsection{{Parallel} transport in global and local
coordinates\label{qng:sec:glCoord}}

Let $X$ be a point on a Grassmannian and consider the geodesic given by
the tangent vector $\Delta\in\mathbf{T}_{X}$. The matrix expression for the
parallel transport of an arbitrary tangent vector $\Delta_{2}\in\mathbf{T}%
_{X}$ is given by
\begin{equation}
\mathbf{T}_{X(t)}\ni\Delta_{2}(t)=\left(  [XV\,\,\,\,U]%
\begin{bmatrix}
-\sin\Sigma t\\
\cos\Sigma t
\end{bmatrix}
U^{\mathsf{T}}+(I-UU^{\mathsf{T}})\right)  \Delta_{2}\equiv T_{X,\Delta
}(t)\Delta_{2}, \label{eq:parVecTran}%
\end{equation}
where $\Delta=U\Sigma V^{T}$ is the thin \textsc{svd} and we define
$T_{X,\Delta}(t)$ to be the \textit{parallel transport matrix}\footnote{We
will often omit subscripts $X$ and $\Delta$ and just write $T(t)$ when there
is no risk for confusion.} from the point $X$ in the direction $\Delta$. If
$\Delta_{2}=\Delta$ expression \eqref{eq:parVecTran} can be simplified.

Let $X\in\operatorname*{Gr}(n,r)$, $\Delta\in\mathbf{T}_{X}$, and $X_{\perp}$
be an orthogonal complement of $X$ so that $[X\,\,X_{\perp}]$ is orthogonal.
Recall that we may write $\Delta=X_{\perp}D$, where we view $X_{\perp}$ as a
basis for $\mathbf{T}_{X}$ and $D$ as a local coordinate representation of the
tangent vector~$\Delta$. Assuming that $X(t)$ is the geodesic curve given in
\eqref{eq:geodesic}, the parallel transport of the corresponding basis
$X_{\perp}(t)$ for $\mathbf{T}_{X(t)}$ is given by
\begin{equation}
X_{\bot}(t)=T_{X,\Delta}(t)X_{\perp},\label{eq:geoPerp}%
\end{equation}
where $T_{X,\Delta}(t)$ is the transport matrix defined in
\eqref{eq:parVecTran}. It is straightforward to show that the matrix
$[X(t)\,\,X_{\perp}(t)]$ is orthogonal for all $t$, i.e.\
\[
X_{\bot}^{\mathsf{T}}(t)X_{\bot}(t)=I_{n-r}\quad\text{ and }\quad X_{\bot
}^{\mathsf{T}}(t)X(t)=0\quad\text{for every }t.
\]
Using \eqref{eq:geoPerp} we can write the parallel transport of a tangent
vector $\Delta_{2}$ as
\begin{equation}
\Delta_{2}(t)=T(t)\Delta_{2}=T(t)X_{\perp}D_{2}=X_{\perp}(t)D_{2}%
.\label{eq:parTran}%
\end{equation}
Equation \eqref{eq:parTran} shows that the local coordinate representation of
the tangent vector is constant at all points of the geodesic path $X(t)$ when
the basis for $\mathbf{T}_{X(t)}$ is given by $X_{\perp}(t)$. The global
coordinate representation, on the other hand, varies with $t$.
This is an important observation since explicit parallel transport of tangents and Hessians (cf.\ Section~\ref{sec:bfgs-opt}) can be avoided if
the algorithm is implemented using local coordinates.
The computational complexity for these two operations\footnote{Here we assume that the parallel transport operator has been computed and stored.} are $O(n^{2}r)$ and $O(n^{3}r^{2})$ respectively. The cost saved in avoiding parallel transports of tangents and Hessians is paid instead in the parallel transport of the basis $X_{\perp}$. This matrix is computed in the first iteration at a cost of at most $O(n^{3})$ operations, and in each of the consecutive iterations, it is parallel transported at a cost of $O(n^{2}(n-r))$ operations. There are also differences in memory requirements: In global coordinates tangents are stored as $n\times r$ matrices and Hessians as  $nr \times nr$ matrices, whereas in local coordinates tangents and Hessians are stored as $(n-r) \times r$ and $(n-r)r \times (n-r)r$ matrices respectively. Local coordinate implementation also requires the additional storage of $X_{\perp}$ as an $n \times (n-r)$ matrix. In most cases, the local coordinate implementation provides greater computational and memory savings, as we observed in our numerical experiments.

By introducing the thin \textsc{svd} of $D=\bar{U}\bar{\Sigma}\bar
{V}^{\mathsf{T}}$, we can also write \eqref{eq:geoPerp} as
\[
X_{\bot}(t)=T(t)X_{\perp}=%
\begin{bmatrix}
X\bar{V} & X_{\perp}\bar{U}%
\end{bmatrix}%
\begin{bmatrix}
-\sin\bar{\Sigma}t\\
\cos\bar{\Sigma}t
\end{bmatrix}
\bar{U}^{\mathsf{T}}+X_{\bot}\left(  I-\bar{U}\bar{U}^{\mathsf{T}}\right)  .
\]
This follows from the identities
\[
\bar{U}=X_{\bot}^{\mathsf{T}}U,\quad\bar{\Sigma}=\Sigma,\quad\bar{V}=V,
\]
which are obtained from $\Delta=U\Sigma V^{\mathsf{T}}=X_{\perp}D$.
Using this, we will derive a general property of
inner products for our later use.

\begin{theorem}
\label{thm:vecTranProd} Let $X\in\operatorname*{Gr}(n,r)$ and $\Delta
,\Delta_{1},\Delta_{2}\in\mathbf{T}_{X}$. Define the transport matrix in the
direction $\Delta$
\[
T_{X,\Delta}(t)=T(t)=%
\begin{bmatrix}
XV & U
\end{bmatrix}%
\begin{bmatrix}
-\sin\Sigma t\\
\cos\Sigma t
\end{bmatrix}
U^{\mathsf{T}}+(I-UU^{\mathsf{T}}),
\]
where $\Delta=U\Sigma V^{\mathsf{T}}$ is the thin \textsc{svd}. Then
\[
\langle\Delta_{1},\Delta_{2}\rangle=\langle\Delta_{1}(t),\Delta_{2}%
(t)\rangle\qquad\text{for every }t,
\]
where $\Delta_{1}(t)=T_{X,\Delta}(t)\Delta_{1}$ and $\Delta_{2}(t)=T_{X,\Delta
}(t)\Delta_{2}$ are parallel transported tangents.
\end{theorem}

\begin{proof}
The proof is a direct consequence of the Levi-Civita connection used in the definition of the parallel transport of tangents. Or we can use the canonical inner product on the Grassmannian $\langle\Delta_{1},\Delta_{2}\rangle=\operatorname*{tr}(\Delta_{1}%
^{\mathsf{T}}\Delta_{2})$. Then, inserting the parallel transported tangents we
obtain
\begin{align*}
\langle\Delta_{1}(t),\Delta_{2}(t)\rangle &  =\operatorname*{tr}%
\bigl({\Delta_{1}(t)}^{\mathsf{T}}\Delta_{2}(t)\bigr)
 =\operatorname*{tr}\bigl(\Delta_{1}^{\mathsf{T}}{T(t)}^{\mathsf{T}%
}T(t)\Delta_{2}\bigr)\\
&  =\operatorname*{tr}\left(  \Delta_{1}^{\mathsf{T}}\left(  I-U\sin(\Sigma
t)V^{\mathsf{T}}X^{\mathsf{T}}-XV\sin(\Sigma t)U^{\mathsf{T}}\right)
\Delta_{2}\right) \\
&  =\operatorname*{tr}\left(  \Delta_{1}^{\mathsf{T}}\Delta_{2}\right)
-\operatorname*{tr}\left(  \Delta_{1}^{\mathsf{T}}U\sin(\Sigma t)V^{\mathsf{T}%
}X^{\mathsf{T}}\Delta_{2}\right)  -\operatorname*{tr}\left(  \Delta
_{1}^{\mathsf{T}}XV\sin(\Sigma t)U^{\mathsf{T}}\Delta_{2}\right)
\hspace{-2pt}.
\end{align*}
The proof is concluded by observing that the second and third terms
after the last equality are zero because $X^{\mathsf{T}}\Delta_{2}=0$ and $\Delta
_{1}^{\mathsf{T}}X=0$.
%\[
%\langle\Delta_{1}(t),\Delta_{2}(t)\rangle= \langle U_{\perp}(t) \Delta_{1},
%U_{\perp}(t) \Delta_{2} \rangle= \langle U_{\perp}(0) \Delta_{1}, U_{\perp}(0)
%\Delta_{2} \rangle= \langle\Delta_{1}(0), \Delta_{2}(0) \rangle.
%\]
\end{proof}

\paragraph{Remark} We would like to point out that it is the tangents that are parallel transported. In global coordinates tangents are represented by $n \times r$ matrices and their parallel transport is given by \eqref{eq:parVecTran}. On the other hand, in local coordinates, tangents are represented by $(n-r)\times r$ matrices and this representation does not change when the basis for the tangent space  is parallel transported according to \eqref{eq:geoPerp}. In other words, in local coordinates, parallel transported tangents are represented by the \textit{same matrix} at every point along a geodesic. This is to be contrasted with the global coordinate representation of points on the manifold $\operatorname*{Gr}(n,r)$, which are $n\times r$ matrices that differ from point to point on a geodesic.

\section{Quasi-Newton methods with BFGS updates on a
Grassmannian\label{qng:sec:QN}}

In this section we will present the necessary modifications in order for
\textsc{bfgs} updates to be well-defined on a Grassmannian. We will write
$f(X)$ instead of $f(x)$ since the argument to the function is a point
on a Grassmannian and represented by a matrix $X=[x_{ij}%
]_{i,j=1}^{n,r}\in\mathbb{R}^{n\times r}$. Similarly the quantities $s_{k}$
and $y_{k}$ from equations \eqref{eq:eucsk} and \eqref{eq:eucyk} will be
written as matrices $S_{k}$ and $Y_{k}$, respectively.

\subsection{Computations in global coordinates}

We describe here the expressions of various quantities required for defining
\textsc{bfgs} updates in global coordinates.
The corresponding expressions in local coordinates are in the next section.

\paragraph{Gradient}

The Grassmann gradient of the objective function $f(X)$ is given by
\begin{equation}
\label{eq:df}\nabla f(X)=\Pi_{X}\frac{\partial f}{\partial X},\quad
\frac{\partial f}{\partial X}:=\left[  \frac{\partial f}{\partial x_{ij}%
}\right]  _{i,j=1}^{n,r},
\end{equation}
where $\Pi_{X}=I-XX^{\mathsf{T}}$ is the projection on the tangent space
$\mathbf{T}_{X}$.

\paragraph{Computing $S_{k}$}

We will now modify the operations in equation \eqref{eq:eucsk}, i.e.\
\[
s_{k}=x_{k+1}-x_{k}=t_{k}p_{k},
\]
so that it is valid on a Grassmannian. Let $X_{k+1}$ be given by
$X_{k+1}=X_{k}(t_{k})$ where the geodesic path originating from $X_{k}$ is
defined by the tangent (or search direction) $\Delta\in\mathbf{T}_{X_{k}}$.
The step size is given by $t_{k}$. We will later assume that $S_{k}%
\in\mathbf{T}_{X_{k+1}}$ and with the tangent $\Delta\in\mathbf{T}_{X_{k}}$,
corresponding to $p_{k}$, we conclude that
\begin{equation}
S_{k}=t_{k}\Delta(t_{k})=t_{k}T(t_{k})\Delta, \label{eq:grSk}%
\end{equation}
where $T(t_{k})$ is the transport matrix defined in \eqref{eq:parVecTran}.

\paragraph{Computing $Y_{k}$}

Similarly, we will translate
\[
y_{k}=g_{k+1}-g_{k}=\nabla f(x_{k+1})-\nabla f(x_{k})
\]
from equation \eqref{eq:eucyk}. Computing the Grassmann gradient at $X_{k+1}$
we get $\nabla f(X_{k+1})\in\mathbf{T}_{X_{k+1}}$. Parallel transporting
$\nabla f(X_{k})\in\mathbf{T}_{X_{k}}$ along the direction $\Delta$ and
subtracting the two gradients as in equation \eqref{eq:eucyk} we get
\begin{equation}
\mathbf{T}_{X_{k+1}}\ni Y_{k}=\nabla f(X_{k+1})-T(t_{k})\nabla f(X_{k}),
\label{eq:grYk}%
\end{equation}
where we again use the transport matrix \eqref{eq:parVecTran}. Recall that
$Y_{k}$ corresponds to~$y_{k}$.

The expressions for $\nabla f$, $S_{k}$ and $Y_{k}$ are given in matrix form,
i.e.\ they have the same dimensions as the variable matrix $X$. It is
straightforward to obtain the corresponding vectorized expressions. For
example, with $\partial f/\partial X\in\mathbb{R}^{n\times r}$, the vector
form of the Grassmann gradient is given by
\[
\operatorname*{vec}(\nabla f)=\left(  I_{r}%
\mathbin{\scalebox{1.12}{$\displaystyle\varotimes$}} \Pi_{X}\right)
\operatorname*{vec}\left(  \frac{\partial f}{\partial X}\right)  \in
\mathbb{R}^{nr}.
\]
where $\operatorname*{vec}(\cdot)$ is the ordinary column-wise vectorization
of a matrix. For simplicity we switch to this presentation when working with
the Hessian.

\paragraph{Updating the Hessian (approximation)}

Identify the tangents (matrices) $\Delta\in\mathbb{R}^{n\times r}$ with
vectors in $\mathbb{R}^{nr}$ and assume that the Grassmann Hessian
\[
\mathbb{R}^{nr\times nr}\ni H_{k}=H(X_{k}):\mathbf{T}_{X_{k}}\rightarrow
\mathbf{T}_{X_{k}}%
\]
at the iterate $X_{k}$ is given. Then
\begin{equation}
\bar{H_{k}}=(I_{r}\mathbin{\scalebox{1.12}{$\displaystyle\varotimes$}}
T(t_{k}))H_{k}(I_{r}\mathbin{\scalebox{1.12}{$\displaystyle\varotimes$}}
\widetilde{T}(t_{k})):\mathbf{T}_{X_{k+1}}\rightarrow\mathbf{T}_{X_{k+1}},
\label{eq:tranHess}%
\end{equation}
is the transported Hessian defined at iterate $X_{k+1}$. As previously
$T(t_{k})$ is the transport matrix from $X_{k}$ to $X_{k+1}$ given in
\eqref{eq:parVecTran} and $\widetilde{T}(t_{k})$ is the transport matrix from
$X_{k+1}$ to $X_{k}$ along the same geodesic path. Informally we can describe
the operations in \eqref{eq:tranHess} as follows. Tangent vectors from
$\mathbf{T}_{X_{k+1}}$ are transported with $\widetilde{T}(t_{k})$ to
$\mathbf{T}_{X_{k}}$ on which $H_{k}$ is defined. The Hessian $H_{k}$
transforms the transported vectors on $\mathbf{T}_{X_{k}}$ and the result is
then forwarded with $T(t_{k})$ to $\mathbf{T}_{X_{k+1}}$.

Since all vectors and matrices are now defined at $X_{k+1}$, the \textsc{bfgs}
update is computed using equation \eqref{eq:BFGS-update} in which we replace
$H_{k}$ with $(I_{r}\mathbin{\scalebox{1.12}{$\displaystyle\varotimes$}}
T(t_{k}))H_{k}(I_{r}\mathbin{\scalebox{1.12}{$\displaystyle\varotimes$}}
\widetilde{T}(t_{k}))$ and use $s_{k}=\operatorname*{vec}(S_{k})$ and
$y_{k}=\operatorname*{vec}(Y_{k})$ from equations \eqref{eq:grSk} and
\eqref{eq:grYk} respectively.

\subsection{Computations in local coordinates}

Using local coordinates we obtain several simplifications. First given the
current iterate $X$ we need the orthogonal complement $X_{\perp}$. When it is
obvious we will omit the iteration subscript $k$.

\paragraph{Grassmann gradient}

In local coordinates the Grassmann gradient is given~by
\begin{equation}
\nabla\widehat{f}=X_{\perp}^{\mathsf{T}}\nabla f=X_{\perp}^{\mathsf{T}}\Pi
_{X}\frac{\partial f}{\partial X}=X_{\perp}^{\mathsf{T}}\frac{\partial
f}{\partial X},\label{eq:grad-lc}%
\end{equation}
where we have used the global coordinate representation for the Grassmann
gradient \eqref{eq:df}. We denote quantities in local coordinates
with a hat to distinguish them from those in global coordinates.

\paragraph{Parallel transporting the basis $X_{\perp}$}

It is necessary to parallel transport the basis matrix $X_{\perp}$ from the
current iterate $X_{k}$ to the next iterate $X_{k+1}$. Only in this basis will
the local coordinates of parallel transported tangents be constant. The
parallel transport of the basis matrix is given by equation \eqref{eq:geoPerp}.

\paragraph{Computing $\widehat{S}_{k}$ and $\widehat{Y}_{k}$}

According to the discussion in Section \ref{qng:sec:glCoord}, in the
transported tangent basis $X_{\perp}(t)$, the local coordinate representation
of any tangent is constant. Specifically this is true for $\widehat{S}_{k}$
and $\widehat{Y}_{k}$. The two quantities are obtained with the same
expressions as in the Euclidean space.

\paragraph{Updating the Hessian (approximation)}

Since explicit {parallel} transport is not required in local coordinates, the
Hessian remains constant as well. The local coordinate representations for
$H_{k}$ in the basis $X_{\perp}(t)$ for points on the geodesic path $X(t)$ are
the same. This statement is proven in Theorem \ref{thm:indlo}.

The effect of using local coordinates on the Grassmannian is only in the
geodesic transport of the current point $X_{k}$ and its orthogonal complement
$X_{k\perp}$. The transported orthogonal complement $X_{k\perp}(t_{k})$ is
used to compute the Grassmann gradient $\nabla\widehat{f}(X_{k+1})$ in local
coordinates at the new iterate $X_{k+1}=X_{k}(t_{k})$. Assuming tangents are
in local coordinates at $X_{k}$ in the basis $X_{k\perp}$ and tangents at
$X_{k+1}$ are given in the basis $X_{k\perp}(t_{k})$, the \textsc{bfgs} update
is give{n} by \eqref{eq:BFGS-update}, i.e.\ exactly the same update as in the
Euclidean space. This is a major advantage compared with the global coordinate
update of $H_{k}$. In global coordinates $H_{k}$ is multiplied by matrices
from the left and from the right \eqref{eq:tranHess}. This is relatively
expensive since the \textsc{bfgs} update itself is just a rank-$2$ update, see
equation \eqref{eq:BFGS-update}.

\subsection{BFGS update in tensor form}

It is not difficult to see that if the gradient $\nabla f(X)$ is written as an
$n\times r$ matrix, then the second derivative will take the form of a
4-tensor $\mathcal{H}_{k}\in\mathbb{R}^{n\times r\times n\times r}$. The
\textsc{bfgs} update \eqref{eq:BFGS-update} can be written in a different form
using the tensor structure of the Hessian. The action of this operator will
map matrices to matrices. Assuming $s_{k}=\operatorname*{vec}(S_{k})$ and
$H_{k}$ is a matricized form of $\mathcal{H}_{k}$, the matrix-vector
contraction $H_{k}s_{k}$ can be written as $\langle\mathcal{H}_{k}%
,S_{k}\rangle_{1,2}$. Obviously the result of the first operation is a vector
whereas the result of the second operation is a matrix, and of course
$H_{k}s_{k}=\operatorname*{vec}\left(  \langle\mathcal{H}_{k},S_{k}%
\rangle_{1,2}\right)  $.

Furthermore keeping the tangents, e.g.\ $S_{k}$ or $Y_{k}$, in matrix form,
the parallel transport of the Hessian in equation \eqref{eq:tranHess} can be
written as a multilinear product between $\mathcal{H}_{k}\in\mathbb{R}%
^{n\times r\times n\times r}$ and the two transport matrices $T(t_{k})$ and
$\widetilde{T}(t_{k})$, both in $\mathbb{R}^{n\times n}$, along the first and
third modes,
\[
\mathbb{R}^{n\times r\times n\times r}\ni\bar{\mathcal{H}}=\mathcal{H}%
_{k}\cdot(T(t_{k}),I,\widetilde{T}(t_{k}),I).
\]

Finally, noting that the outer product between vectors corresponds to tensor
products between matrices the \textsc{bfgs} update becomes{\footnote{{The
contractions denoted by $\langle\cdot,\cdot\rangle_{\ast}$ are defined in
Appendix \ref{app:1}.}}}
\begin{equation}
\mathcal{H}_{k+1}=\bar{\mathcal{H}}_{k}+\frac{\langle\bar{\mathcal{H}}%
_{k},S_{k}\rangle_{1,2}\mathbin{\scalebox{.88}{$\displaystyle\otimes$}}\langle
\bar{\mathcal{H}}_{k},S_{k}\rangle_{1,2}}{\langle\langle\bar{\mathcal{H}}%
_{k},S_{k}\rangle_{1,2},S_{k}\rangle}+\frac{Y_{k}%
\mathbin{\scalebox{.88}{$\displaystyle\otimes$}}Y_{k}}{\langle S_{k}%
,Y_{k}\rangle},\label{eq:bfgs_tfg}%
\end{equation}
where the matrices $S_{k},Y_{k}\in\mathbf{T}_{X_{k+1}}$ are given by
\eqref{eq:grSk} and \eqref{eq:grYk} respectively.

In local coordinates the update is even simpler since we do not have to
parallel transport the Hessian operator,
\begin{equation}
\widehat{\mathcal{H}}_{k+1}=\widehat{\mathcal{H}}_{k}+\frac{\langle
\widehat{\mathcal{H}}_{k},\widehat{S}_{k}\rangle_{1,2}%
\mathbin{\scalebox{.88}{$\displaystyle\otimes$}}\langle\widehat{\mathcal{H}%
}_{k},\widehat{S}_{k}\rangle_{1,2}}{\langle\langle\widehat{\mathcal{H}}%
_{k},\widehat{S}_{k}\rangle_{1,2},\widehat{S}_{k}\rangle}+\frac{\widehat
{Y}_{k}\mathbin{\scalebox{.88}{$\displaystyle\otimes$}}\widehat{Y}_{k}%
}{\langle\widehat{S}_{k},\widehat{Y}_{k}\rangle},\label{eq:bfgs_tfl}%
\end{equation}
where $\widehat{\mathcal{H}}_{k}\in\mathbb{R}^{(n-r)\times r\times(n-r)\times
r}$ and $\widehat{S}_{k},\widehat{Y}_{k}\in\mathbb{R}^{(n-r)\times r}$. The
\textit{hat} indicates that the corresponding variables are in local coordinates.

\subsection{BFGS update on a product of Grassmannians}

Assume now that the objective function $f$ is defined on a product of
three\footnote{We assume $k=3$ for notational simplicity; generalization of
these discussions to arbitrary $k$ is straightforward.} Grassmannians, i.e.\
\[
f:\operatorname*{Gr}(l,p)\times\operatorname*{Gr}(m,q)\times\operatorname*{Gr}%
(n,r)\rightarrow\mathbb{R},
\]
and is twice continuously differentiable. We write $f(X,Y,Z)$ where
$X\in\operatorname*{Gr}(l,p)$, $Y\in\operatorname*{Gr}(m,q)$ and
$Z\in\operatorname*{Gr}(n,r)$. The Hessian of the objective function will have
a `block tensor' structure but the blocks will not have conforming dimensions.
The action of a (approximate) Hessian operator on tangents $\Delta_{X}%
\in\mathbf{T}_{X}$, $\Delta_{Y}\in\mathbf{T}_{Y}$ and $\Delta_{Z}\in
\mathbf{T}_{Z}$ may be written symbolically as
\begin{align}
&  \mathrel{\phantom{=}}%
\begin{bmatrix}
\mathcal{H}_{XX} & \mathcal{H}_{XY} & \mathcal{H}_{XZ}\\
\mathcal{H}_{YX} & \mathcal{H}_{YY} & \mathcal{H}_{YZ}\\
\mathcal{H}_{ZX} & \mathcal{H}_{ZY} & \mathcal{H}_{ZZ}%
\end{bmatrix}%
\begin{bmatrix}
\Delta_{X}\\
\Delta_{Y}\\
\Delta_{Z}%
\end{bmatrix}
\label{eq:prod-op}\\
&  =%
\begin{bmatrix}
\langle\mathcal{H}_{XX},\Delta_{X}\rangle_{3,4;1,2}+\langle\mathcal{H}%
_{XY},\Delta_{Y}\rangle_{3,4;1,2}+\langle\mathcal{H}_{XZ},\Delta_{Z}%
\rangle_{3,4;1,2}\\
\langle\mathcal{H}_{YX},\Delta_{X}\rangle_{3,4;1,2}+\langle\mathcal{H}%
_{YY},\Delta_{Y}\rangle_{3,4;1,2}+\langle\mathcal{H}_{YZ},\Delta_{Z}%
\rangle_{3,4;1,2}\\
\langle\mathcal{H}_{ZX},\Delta_{X}\rangle_{3,4;1,2}+\langle\mathcal{H}%
_{ZY},\Delta_{Y}\rangle_{3,4;1,2}+\langle\mathcal{H}_{ZZ},\Delta_{Z}%
\rangle_{3,4;1,2}%
\end{bmatrix}
\nonumber\\
&  =%
\begin{bmatrix}
\langle\mathcal{H}_{XX},\Delta_{X}\rangle_{1,2}+\langle\mathcal{H}_{YX}%
,\Delta_{Y}\rangle_{1,2}+\langle\mathcal{H}_{ZX},\Delta_{Z}\rangle_{1,2}\\
\langle\mathcal{H}_{XY},\Delta_{X}\rangle_{1,2}+\langle\mathcal{H}_{YY}%
,\Delta_{Y}\rangle_{1,2}+\langle\mathcal{H}_{ZY},\Delta_{Z}\rangle_{1,2}\\
\langle\mathcal{H}_{XZ},\Delta_{X}\rangle_{1,2}+\langle\mathcal{H}_{YZ}%
,\Delta_{Y}\rangle_{1,2}+\langle\mathcal{H}_{ZZ},\Delta_{Z}\rangle_{1,2}%
\end{bmatrix}
.\nonumber
\end{align}
The blocks of the Hessian are $4$-tensors and elements of the tangent spaces
are matrices. The result of the operation\footnotemark[3] is a triplet where
each element is in the corresponding tangent space. For example $\mathcal{H}%
_{XX}$ is an $l\times p\times l\times p$ tensor which acts on the tangent
matrix $\Delta_{X}$ of size $l\times p$ with the result $\langle
\mathcal{H}_{XX},\Delta_{X}\rangle_{1,2}\in\mathbf{T}_{X}$. Off diagonal
example may look as follows, $\mathcal{H}_{YZ}$ is an $m\times q\times n\times
r$ tensor which acts on the tangent matrix $\Delta_{Z}$ of size $n\times r$
with the result $\langle\mathcal{H}_{YZ},\Delta_{Z}\rangle_{3,4;1,2}%
=\langle\mathcal{H}_{ZY},\Delta_{Z}\rangle_{1,2}\in\mathbf{T}_{Y}$. The
equality in the last step follows from the fact that the $n\times r\times
m\times q$ tensor $\mathcal{H}_{ZY}$ is a permutation of the $m\times q\times
n\times r$ tensor $\mathcal{H}_{YZ}$. This is expected since for twice
continuously differentiable functions $f_{xy}=f_{yx}$. But in our case they
have different `shapes'. The three tangent spaces $\mathbf{T}_{X}$,
$\mathbf{T}_{Y}$ and $\mathbf{T}_{Z}$ are interconnected through the Hessian
of $f(X,Y,Z)$ in the sense that every block in \eqref{eq:prod-op} is a linear
operator mapping matrices from one tangent space to another tangent space. For
example $\mathcal{H}_{YX}:\mathbf{T}_{X}\rightarrow\mathbf{T}_{Y}$ and
$\mathcal{H}_{ZX}:\mathbf{T}_{X}\rightarrow\mathbf{T}_{Z}$.

The corresponding \textsc{bfgs} in the product manifold case has basically the
same form as equations \eqref{eq:bfgs_tfg} and \eqref{eq:bfgs_tfl} where the
action of the Hessian on $S_{k}$, which will be a triplet with an element on
each tangent space, is replaced with formulas as in \eqref{eq:prod-op}. Also
the tensor/outer product needs to be modified in the obvious way, i.e.\ if
$\Delta=(\Delta_{X},\Delta_{Y},\Delta_{Z})$ and $\Gamma=(\Gamma_{X},\Gamma
_{Y},\Gamma_{Z})$ then we let
\begin{align}
\Delta\mathbin{\scalebox{.88}{$\displaystyle\hat{\otimes}$}} \Gamma &
=(\Delta_{X},\Delta_{Y},\Delta_{Z}%
)\mathbin{\scalebox{.88}{$\displaystyle\hat{\otimes}$}} (\Gamma_{X},\Gamma
_{Y},\Gamma_{Z})\label{eq:modop}\\
&  :=%
\begin{bmatrix}
\Delta_{X}\mathbin{\scalebox{.88}{$\displaystyle\otimes$}} \Gamma_{X} &
\Delta_{X}\mathbin{\scalebox{.88}{$\displaystyle\otimes$}} \Gamma_{Y} &
\Delta_{X}\mathbin{\scalebox{.88}{$\displaystyle\otimes$}} \Gamma_{Z}\\
\Delta_{Y}\mathbin{\scalebox{.88}{$\displaystyle\otimes$}} \Gamma_{X} &
\Delta_{Y}\mathbin{\scalebox{.88}{$\displaystyle\otimes$}} \Gamma_{Y} &
\Delta_{Y}\mathbin{\scalebox{.88}{$\displaystyle\otimes$}} \Gamma_{Z}\\
\Delta_{Z}\mathbin{\scalebox{.88}{$\displaystyle\otimes$}} \Gamma_{X} &
\Delta_{Z}\mathbin{\scalebox{.88}{$\displaystyle\otimes$}} \Gamma_{Y} &
\Delta_{Z}\mathbin{\scalebox{.88}{$\displaystyle\otimes$}} \Gamma_{Z}%
\end{bmatrix}
,\nonumber
\end{align}
where the results are conveniently stored in a `block matrix' whose blocks are
tensors of different dimensions (possibly nonconforming).

\subsection{Optimality of BFGS on Grassmannians\label{sec:bfgs-opt}}

The \textsc{bfgs} update in quasi-Newton methods is optimal because it is the
solution to
\[
\min_{H\in\mathbb{R}^{n\times n}}\lVert H-H_{k}\rVert_{F}\qquad\text{subject
to}\qquad H=H^{\mathsf{T}},\qquad Hs_{k}=y_{k},
\]
where $s_{k}$ and $y_{k}$ are given by \eqref{eq:eucsk} and \eqref{eq:eucyk}
respectively \cite{nowr:06}. For the Euclidean case it is immaterial
whether $H$ is considered as an abstract operator or explicitly
represented as a matrix. The final conclusion with respect to
optimality is the same---it amounts to a rank-$2$ change of $H_{k}%
$. The situation is different when considering the corresponding optimality
problem on Grassmannians. In particular, a given Hessian (or approximate
Hessian) matrix $H_{k}$ considered in a global coordinate representation and
defined at $X_{k}\in\operatorname*{Gr}(n,r)$ has the following form when
parallel transported along a geodesic,
\[
\bar{H_{k}}=(I_{r}\mathbin{\scalebox{1.12}{$\displaystyle\varotimes$}}T(t_{k}%
))H_{k}(I_{r}\mathbin{\scalebox{1.12}{$\displaystyle\varotimes$}}\widetilde
{T}(t_{k})).
\]
This is the same expression as equation \eqref{eq:tranHess}.
While the Hessian operator should not change by a
parallel transport to a new point on the manifold, its representation
evidently changes. This has important numerical and
computational ramifications. In fact, the global coordinate
representation of the Hessian at the previous point is usually very different
from the global coordinate representation of the transported Hessian at the
current point.

Assume now the Hessian matrix (or its approximation) is given in local
coordinates $\widehat{H}_{k}$ at $X_{k}\in\operatorname*{Gr}(n,r)$ and let
$X_{k\,\perp}$ be the associated basis matrix for the tangent space.
Representation of the parallel transported Hessian will not change
if the associated basis matrix $X_{k\,\perp}(t)$ is transported
according to \eqref{eq:geoPerp}. The updated Hessian at the current
point is a rank-$2$ modification of the Hessian from the previous
point given by the \textsc{bfgs} update. The optimality of \textsc{bfgs} update on
Euclidean spaces is with respect to a change in successive Hessian
matrices; we will prove that in the correct tangent space basis and
in local coordinates, the \textsc{bfgs} update is also optimal on
Grassmannians.

We now give a self contained proof for this statement. First we will state the
optimality results for the Euclidean case. The proofs of
Theorem~\ref{thm:rank1update}, Lemma~\ref{lem:SPDexist}, and
Theorem~\ref{thm:eucl-bfgs-opt} are based on \cite{denmor77,densch81}. We will
then use these to deduce the corresponding optimality result on a product of
Grassmannians in Theorem~\ref{thm:Opt}.

\begin{theorem}
\label{thm:rank1update}Let $B\in\mathbb{R}^{n\times n}$, $y\in\mathbb{R}^{n}$,
$0\neq s\in\mathbb{R}^{n}$. The solution to%
\[
\min\{\lVert A-B\rVert_{F}\mid A\in\mathbb{R}^{n\times n},\;As=y\}
\]
is given by%
\[
\bar{B}=B+\frac{(y-Bs)s^{\mathsf{T}}}{s^{\mathsf{T}}s}.
\]

\end{theorem}

\begin{proof}
Note that while the set $Q(y,s):=\{A\in\mathbb{R}^{n\times n}\mid As=y\}$ is
non-compact (closed but unbounded), for a fixed $B$, the function
$f:Q(y,s)\rightarrow\mathbb{R}$, $f(A)=\lVert A-B\rVert_{F}$ is coercive and
therefore a minimizer $A_{\ast}\in Q(y,s)$ is attained. This demonstrates
existence. The minimizer is also unique since $Q(y,s)$ is convex while $f$ is
strictly convex. We claim that $A_{\ast}=\bar{B}$: Observe that $\bar{B}s=y$
and so $\bar{B}\in Q(y,s)$; for any $A\in Q(y,s)$,
\[
\lVert\bar{B}-B\rVert_{F}=\left\|  \frac{ys^{\mathsf{T}}}{s^{\mathsf{T}}%
s}-\frac{Bss^{\mathsf{T}}}{s^{\mathsf{T}}s}\right\|  _{F}=\left\|
(A-B)\frac{ss^{\mathsf{T}}}{s^{\mathsf{T}}s}\right\|  _{F}\leq\lVert
A-B\rVert_{F}\left\|  \frac{ss^{\mathsf{T} }}{s^{\mathsf{T}}s}\right\|
_{F}=\lVert A-B\rVert_{F}.
\]

\end{proof}

\begin{lemma}
\label{lem:SPDexist}Let $y\in\mathbb{R}^{n}$, $0\neq s\in\mathbb{R}^{n}$. Then
the set $Q(y,s)=\{A\in\mathbb{R}^{n\times n}\mid As=y\}$ contains a symmetric
positive definite matrix iff $y=Lv$ and $v=L^{\mathsf{T}}s$ for some $0\neq
v\in\mathbb{R}^{n}$ and $L\in\operatorname*{GL}(n)$.
\end{lemma}

\begin{proof}
If such $v$ and $L$ exist, then $y=Lv=LL^{\mathsf{T}}s$ and so $LL^{\mathsf{T}%
}$ is a symmetric positive definite matrix in $Q(y,s)$. On the other hand, if
$A\in Q(y,s)$ is symmetric positive definite, its Cholesky factorization
$A=LL^{\mathsf{T}}$ yields an $L\in\operatorname*{GL}(n)$. If we let
$v=L^{\mathsf{T}}s$, then $Lv=As=y$, as required.
\end{proof}

\begin{theorem}
\label{thm:eucl-bfgs-opt} Let $y\in\mathbb{R}^{n}$, $0\neq s\in\mathbb{R}^{n}%
$. Let $L\in\operatorname*{GL}(n)$ and $H=LL^{\mathsf{T}}$. There is a
symmetric positive definite matrix $H_{+}\in Q(y,s)$ iff $y^{\mathsf{T}}s>0$.
In this case, the \textsc{bfgs} update $H_{+}=L_{+}L_{+}^{\mathsf{T}}$ is one
where%
\begin{equation}
L_{+}=L+\frac{(y-\alpha Hs)(L^{\mathsf{T}}s)^{\mathsf{T}}}{\alpha
s^{\mathsf{T}}Hs} \quad\text{ with } \quad\alpha=\pm\sqrt{\frac{y^{\mathsf{T}%
}s}{s^{\mathsf{T}}Hs}}. \label{update}%
\end{equation}
%where%
%\[
%\]

\end{theorem}

\begin{proof}
In order for the update \eqref{update} to exist it is necessary that there
exists $0\neq v\in\mathbb{R}^{n}$ and $L_{+}\in\operatorname*{GL}(n)$ such
that $y=L_{+}v$ and $v=L_{+}^{\mathsf{T}}s$. Hence%
\[
0<v^{\mathsf{T}}v=(L_{+}^{\mathsf{T}}s)^{\mathsf{T}}(L_{+}^{-1}%
y)=s^{\mathsf{T}}y
\]
as required.

If $v$ is known, then the nearest matrix to $L$ that takes $v$ to $y$ would be
the update given in Theorem~\ref{thm:rank1update}, i.e.%
\[
L_{+}=L+\frac{(y-Lv)v^{\mathsf{T}}}{v^{\mathsf{T}}v}.
\]
Hence we need to find the vector $v$. By Lemma~\ref{lem:SPDexist},%
\begin{equation}
v=L_{+}^{\mathsf{T}}s=L^{\mathsf{T}}s+\frac{y^{\mathsf{T}}s-v^{\mathsf{T}%
}L^{\mathsf{T}}s}{v^{\mathsf{T}}v}v \label{eq1}%
\end{equation}
and so
\begin{equation}
v=\alpha L^{\mathsf{T}}s \label{eq2}%
\end{equation}
for some $\alpha\in\mathbb{R}$. Now it remains to find the scalar $\alpha$.
Plugging \eqref{eq2} into \eqref{eq1} and using $H=LL^{\mathsf{T}}$, we get%
\[
\alpha=1+\frac{y^{\mathsf{T}}s-\alpha s^{\mathsf{T}}Hs}{\alpha^{2}%
s^{\mathsf{T}}Hs}\cdot\alpha\quad\Rightarrow\quad\alpha^{2}=\frac
{y^{\mathsf{T}}s}{s^{\mathsf{T}}Hs}.
\]
If $y^{\mathsf{T}}s>0$, this defines an update in $Q(y,s)$ that is symmetric
positive definite. It is straightforward to verify that $H_{+}=L_{+}%
L_{+}^{\mathsf{T}}$ yields the \textsc{bfgs} update%
\[
H_{+}=H-\frac{Hss^{\mathsf{T}}H}{s^{\mathsf{T}}Hs}+\frac{yy^{\mathsf{T}}%
}{y^{\mathsf{T}}s}.
\]

\end{proof}

\begin{theorem}
\label{thm:indlo} Let $X\in\operatorname*{Gr}(n,r)$ and $X_{\perp}$ be the
orthogonal complement to $X$, i.e.\ $[X\,X_{\perp}]$ is orthogonal. Let
$\Delta\in\mathbf{T}_{X}$ and $X_{\Delta}(t)$ be a geodesic with $T_{X,\Delta
}(t)$ the corresponding transport matrix, defined according to equations
\eqref{eq:geodesic} and \eqref{eq:parVecTran}. Identify $\mathbf{T}_{X}$ with
$\mathbb{R}^{(n-r)r}$ and consider a linear operator in local coordinates
$\hat{A}:\mathbb{R}^{(n-r)r}\rightarrow\mathbb{R}^{(n-r)r}$. Consider the
corresponding linear operator in global coordinates $A:\mathbb{R}%
^{nr}\rightarrow\mathbb{R}^{nr}$, in which tangents in $\mathbf{T}_{X}$ are
embedded. The relation between the two operators is given by
\begin{align}
A  &  =(I\mathbin{\scalebox{1.12}{$\displaystyle\varotimes$}} X_{\perp}%
)\hat{A}(I\mathbin{\scalebox{1.12}{$\displaystyle\varotimes$}} X_{\perp
}^{\mathsf{T}}),\label{eq:ident1}\\
\hat{A}  &  =(I\mathbin{\scalebox{1.12}{$\displaystyle\varotimes$}} X_{\perp
}^{\mathsf{T}})A(I\mathbin{\scalebox{1.12}{$\displaystyle\varotimes$}}
X_{\perp}). \label{eq:ident2}%
\end{align}
Furthermore, the parallel transported operator $\hat{A}$ has the same
representation for all $t$ along the geodesic $X(t)$, i.e.\ $\hat{A}%
(t)\equiv\hat{A}$.
\end{theorem}

\begin{proof}
Let $d_{1}\in\mathbb{R}^{(n-r)r}$ be a tangent vector with corresponding
global coordinate matrix representation $\Delta_{1}=X_{\perp}D_{1}%
\in\mathbf{T}_{X}$. Obviously $d_{1}=\operatorname*{vec}(D_{1})$. We may write
$\operatorname*{vec}(\Delta_{1}%
)=(I\mathbin{\scalebox{1.12}{$\displaystyle\varotimes$}} X_{\perp})d_{1}$. Set
$d_{2}=\hat{A}d_{1}$ and it follows that $\operatorname*{vec}(\Delta
_{2})=(I\mathbin{\scalebox{1.12}{$\displaystyle\varotimes$}} X_{\perp})d_{2}.$
The corresponding operation in global coordinates are
\begin{align*}
\operatorname*{vec}(\Delta_{2})=A\operatorname*{vec}(\Delta_{1})\quad &
\Leftrightarrow\quad(I\mathbin{\scalebox{1.12}{$\displaystyle\varotimes$}}
X_{\perp})d_{2}=A(I\mathbin{\scalebox{1.12}{$\displaystyle\varotimes$}}
X_{\perp})d_{1}\\
&  \Leftrightarrow\quad d_{2}%
=(I\mathbin{\scalebox{1.12}{$\displaystyle\varotimes$}} X_{\perp}^{\mathsf{T}%
})A(I\mathbin{\scalebox{1.12}{$\displaystyle\varotimes$}} X_{\perp})d_{1}%
\end{align*}
and it follows that $\hat{A}%
=(I\mathbin{\scalebox{1.12}{$\displaystyle\varotimes$}} X_{\perp}^{\mathsf{T}%
})A(I\mathbin{\scalebox{1.12}{$\displaystyle\varotimes$}} X_{\perp})$, which
proves \eqref{eq:ident2}.

For any tangent $\Delta_{\ast}\in\mathbf{T}_{X}$ it holds that $\Delta_{\ast
}=\Pi_{X}\Delta_{\ast}=X_{\perp}X_{\perp}^{\mathsf{T}}\Delta_{\ast}$, where
$\Pi_{X}$ is a projection onto $\mathbf{T}_{X}$, and consequently
$\operatorname*{vec}(\Delta_{\ast}%
)=(I\mathbin{\scalebox{1.12}{$\displaystyle\varotimes$}} X_{\perp}X_{\perp
}^{\mathsf{T}})\operatorname*{vec}(\Delta_{\ast})$. Thus the operations in
global coordinates also satisfy
\begin{align*}
\operatorname*{vec}(\Delta_{2})  &
=(I\mathbin{\scalebox{1.12}{$\displaystyle\varotimes$}} X_{\perp}X_{\perp
}^{\mathsf{T}})A(I\mathbin{\scalebox{1.12}{$\displaystyle\varotimes$}}
X_{\perp}X_{\perp}^{\mathsf{T}})\operatorname*{vec}(\Delta_{1})\\
&  =(I\mathbin{\scalebox{1.12}{$\displaystyle\varotimes$}} X_{\perp})\hat
{A}(I\mathbin{\scalebox{1.12}{$\displaystyle\varotimes$}} X_{\perp
}^{\mathsf{T}})\operatorname*{vec}(\Delta_{1})\\
&  \equiv A\operatorname*{vec}(\Delta_{1}).
\end{align*}
This proves equation \eqref{eq:ident1}.

For the third part we have $A:\mathbf{T}_{X}\rightarrow\mathbf{T}_{X}$ and
$A(t):\mathbf{T}_{X(t)}\rightarrow\mathbf{T}_{X(t)}$ with $A(0)\equiv A$. We
want to prove that $\hat{A}(t)\equiv\hat{A}$ for all $t$. The operator $A(t)$
is defined in the following sense: a tangent $\Delta_{1}(t)\in\mathbf{T}%
_{X(t)}$ is parallel transported with $\widetilde{T}_{X(t),-\Delta(t)}(t)$ to
$X(0)$ along $X(t)$, the operator transformations is performed in
$\mathbf{T}_{X}$, thus $\Delta_{2}=A(\Delta_{1}(0))\in\mathbf{T}_{X}$ and the
result is forwarded to $\mathbf{T}_{X(t)}$, i.e.\ $\Delta_{2}(t)=T_{X,\Delta
}(t)\Delta_{2}$. The parallel transported operator in global coordinates takes
the form
\begin{equation}
A(t)=(I\mathbin{\scalebox{1.12}{$\displaystyle\varotimes$}} T_{X,\Delta
}(t))A(I\mathbin{\scalebox{1.12}{$\displaystyle\varotimes$}} \widetilde
{T}_{X(t),-\Delta(t)}(t)). \label{eq:op}%
\end{equation}
Then, in the basis $X_{\perp}(t)$, the local coordinate representation of the
operator is
\begin{equation}
\label{eq:opt}\hat{A}%
(t)=(I\mathbin{\scalebox{1.12}{$\displaystyle\varotimes$}} X_{\perp
}^{\mathsf{T}}(t))A(t)(I\mathbin{\scalebox{1.12}{$\displaystyle\varotimes$}}
X_{\perp}(t)).
\end{equation}
Substituting \eqref{eq:op} into \eqref{eq:opt}, we obtain
\[
\hat{A}(t)=(I\mathbin{\scalebox{1.12}{$\displaystyle\varotimes$}} X_{\perp
}^{\mathsf{T}}(t)T_{X,\Delta}(t)X_{\perp})\hat{A}%
(I\mathbin{\scalebox{1.12}{$\displaystyle\varotimes$}} X_{\perp}^{\mathsf{T}%
}\widetilde{T}_{X(t),-\Delta(t)}(t)X_{\perp}(t)).
\]
Recall that $T_{X,\Delta}(t)X_{\perp}=X_{\perp}(t)$ and thus $X_{\perp
}^{\mathsf{T}}(t)T_{X,\Delta}(t)X_{\perp}=I$. Similarly one can show that
$X_{\perp}^{\mathsf{T}}\widetilde{T}_{X(t),-\Delta(t)}(t)X_{\perp}(t)=I$ and
we get $\hat{A}(t)=\hat{A}$ for all $t$.
\end{proof}

A different proof of essentially the same statement may be found in
\cite{simon07}.

\begin{lemma}
\label{lem:tind-prod-grass} Let $X_{i}\in\operatorname*{Gr}(n_{i},r_{i})$,
$i=1,\dots,k,$ with corresponding tangent spaces $\mathbf{T}_{X_{i}}$. Let
$Y_{i}$ be given such that $[X_{i}\,\,Y_{i}]$ is orthogonal. On each
Grassmannian $\operatorname*{Gr}(n_{i},r_{i})$, let $X_{i}(t)$ be a geodesic
and $Y_{i}(t)$ be its orthogonal complement corresponding to the tangent
$\Delta_{i}\in\mathbf{T}_{X_{i}}$. Then a local coordinate representation of
the linear operator
\[
\hat{A}:\mathbf{T}_{X_{1}}\times\dots\times\mathbf{T}_{X_{k}}\rightarrow
\mathbf{T}_{X_{1}}\times\dots\times\mathbf{T}_{X_{k}}%
\]
is independent of $t$ when parallel transported along the geodesics $X_{i}(t)$
and in the tangent basis $Y_{i}(t)$.
\end{lemma}

\begin{proof}
First we observe that the operator $\hat{A}$ must necessarily have the
structure
\[
\hat{A} =
\begin{bmatrix}
\hat{A}_{11} & \cdots & \hat{A}_{1k}\\
\vdots & \ddots & \vdots\\
\hat{A}_{k1} & \cdots & \hat{A}_{kk}%
\end{bmatrix}
\]
where each $\hat{A}_{ij}$, $1 \leq i,j \leq k$ is such that $\hat{A}_{ij} :
\mathbf{T}_{X_{j}} \rightarrow\mathbf{T}_{X_{i}}$. Now, applying a similar
procedure as in Theorem~\ref{thm:indlo} on each block $\hat{A}_{ij}$ proves
that local coordinate representation of $\hat{A}_{ij}$ and thus $\hat{A}$ is
independent of $t$ along the geodesics $X_{i}(t)$ in the tangent space basis
$Y_{i}(t)$.
\end{proof}

Now we will give an explicit expression for the general \textsc{bfgs} update
in tensor form and in local coordinates. We omit the hat and iteration index
below for clarity. For a function defined on a product of $k$ Grassmannians
$f:\operatorname*{Gr}(n_{1},r_{1})\times\dots\times\operatorname*{Gr}%
(n_{k},r_{k})\rightarrow\mathbb{R}$, we write $f(X_{1},\dots,X_{k})$ and
$S=(S_{1},\dots,S_{k})$, $Y=(Y_{1},\dots,Y_{k})$ where $X_{i}\in
\operatorname*{Gr}(n_{i},r_{i})$ and $S_{i},Y_{i}\in\mathbf{T}_{X_{i}}$ for
$i=1,\dots,k$. The Hessian or its approximation has the symbolic form
\[
\mathcal{H}=%
\begin{bmatrix}
\mathcal{H}_{11} & \cdots & \mathcal{H}_{1k}\\
\vdots & \ddots & \vdots\\
\mathcal{H}_{k1} & \cdots & \mathcal{H}_{kk}%
\end{bmatrix}
\]
where each block is a $4$-tensor. The \textsc{bfgs} update takes the form,
\begin{equation}
\mathcal{H}_{+}=\mathcal{H}+\frac{\langle\mathcal{H},S\rangle_{1,2}%
\mathbin{\scalebox{.88}{$\displaystyle\hat{\otimes}$}} \langle\mathcal{H}%
,S\rangle_{1,2}}{\langle\langle\mathcal{H},S\rangle_{1,2},S\rangle}%
+\frac{Y\mathbin{\scalebox{.88}{$\displaystyle\hat{\otimes}$}} Y}{\langle
S,Y\rangle}, \label{eq:BFGSprod}%
\end{equation}
where the $\langle\mathcal{H},S\rangle_{1,2}$ is given by a formula similar to
\eqref{eq:prod-op} with the result being a $k$-tuple, the tensor product
between $k$-tuples of tangents is an obvious generalization of
\eqref{eq:modop}, and of course $\langle S,Y\rangle=\sum_{i=1}^{d}\langle
S_{i},Y_{i}\rangle$.

Finally, we have all the ingredients required to prove the optimality of the
\textsc{bfgs} update on a product of Grassmannians.

\begin{theorem}
[Optimality of \textsc{bfgs} update on product of Grassmannians]%
\label{thm:Opt}Consider a function $f(X_{1},\dots,X_{k})$ in the variables
$X_{i}\in\operatorname*{Gr}(n_{i},r_{i})$, $i=1,\dots,k,$ that we want to
minimize. Let $X_{i}(t)$ be geodesic defined by $\Delta_{i}\in\mathbf{T}%
_{X_{i}}$ with the corresponding tangent space basis matrices $Y_{i}(t)$. In
these basis for the tangent spaces, the \textsc{bfgs} updates in
\eqref{eq:BFGSprod} on the product Grassmannians have the same optimality
properties as a function with variables in a Euclidean space, i.e.\ it is the
least change update of the current Hessian approximation that satisfies the
secant equations.
\end{theorem}

\begin{proof}
First we observe that the Grassmann Hessian of $f(X_{1},\dots,X_{k})$ (or its
approximation) is a linear operator
\[
H_{f}:\mathbf{T}_{X_{1}}\times\dots\times\mathbf{T}_{X_{k}}\rightarrow
\mathbf{T}_{X_{1}}\times\dots\times\mathbf{T}_{X_{k}}%
\]
and according to Lemma~\ref{lem:tind-prod-grass} its local coordinate
representation is constant along the geodesics $X_{i}(t)$. Given this, the
\textsc{bfgs} optimality result on product Grassmannians is a consequence from
Theorem~\ref{thm:eucl-bfgs-opt}---the optimality of \textsc{bfgs} in Euclidean space.
\end{proof}

\paragraph{Remark}

An important difference on (product) Grassmannians is that we need to keep
track of the basis for the tangent spaces---$Y_{i}(t)$ from equation
\eqref{eq:geoPerp}. Only then will the local coordinate representation of an
operator be independent of $t$ when transported along geodesics.

Note that Theorem~\ref{thm:Opt} is a coordinate dependent result. If
we regard Hessians as abstract operators, there will no longer be any
difference between the global and the local scenario. But the corresponding
optimality as the least amount of change in successive Hessians cannot be
obtained in global coordinate representation and is thus not true if the
Hessians are regarded as abstract operators.

\subsection{Other alternatives}

Movement along geodesics and parallel transport of tangents are the
most straightforward and natural generalizations to the key operations from
Euclidean spaces to manifolds. There are also methods for dealing
with the manifold structure in optimization algorithms based on
different principles. For example, instead of moving along geodesics from one
point to another on the manifold one could use the notion of
\textit{retractions}, which is a smooth mapping from the tangent bundle of the
manifold onto the manifold. Another example is the notion of \textit{vector
transport} that generalizes the parallel translation/transport of tangents
used in this paper. All these notions are defined and
described in \cite{absil08}. It is not clear how the use of the more  general vector transport would effect the convergence properties of the resulting \textsc{bfgs} methods.

\section{Limited memory BFGS\label{qng:sec:lbfgs}}

We give a brief summary of the limited memory quasi-Newton method with
\textsc{l-bfgs} updates on Euclidean spaces \cite{byrd94} that we need later
for our Grassmann variant. See also the discussion in \cite[Chapter~7]%
{simon07}. In Euclidean space the \textsc{bfgs} update can be represented in
the following compact form
\begin{equation}
H_{k}=H_{0}+%
\begin{bmatrix}
S_{k} & H_{0}Y_{k}%
\end{bmatrix}%
\begin{bmatrix}
R_{k}^{-\mathsf{T}}(D_{k}+Y_{k}^{\mathsf{T}}H_{0}Y_{k})R_{k}^{-1} &
-R_{k}^{-\mathsf{T}}\\
-R_{k}^{-1} & 0
\end{bmatrix}%
\begin{bmatrix}
S_{k}^{\mathsf{T}}\\
Y_{k}^{\mathsf{T}}H_{0}%
\end{bmatrix}
, \label{qng:eg:cf_bfgs}%
\end{equation}
where $S_{k}=\left[  s_{0},\dots,s_{k-1}\right]  $, $Y_{k}=\left[  y_{0}%
,\dots,y_{k-1}\right]  $, $D_{k}=\operatorname*{diag}\left[  s_{0}%
^{\mathsf{T}}y_{0},\dots,s_{k-1}^{\mathsf{T}}y_{k-1}\right]  $ and
\[
R_{k}=%
\begin{bmatrix}
s_{0}^{\mathsf{T}}y_{0} & s_{0}^{\mathsf{T}}y_{1} & \cdots & s_{0}%
^{\mathsf{T}}y_{k-1}\\
0 & s_{1}^{\mathsf{T}}y_{1} & \cdots & s_{1}^{\mathsf{T}}y_{k-1}\\
\vdots &  & \ddots & \vdots\\
0 & \cdots & 0 & s_{k-1}^{\mathsf{T}}y_{k-1}%
\end{bmatrix}
,
\]
are obtained using equations \eqref{eq:eucsk} and \eqref{eq:eucyk}. Observe
that in this section $S_{k}$ and $Y_{k}$ are not the same as in
\eqref{eq:grSk} and \eqref{eq:grYk} respectively. The limited memory version
of the algorithm is obtained when replacing the initial Hessian $H_{0}$ by a
sparse matrix, usually this is a suitably scaled identity matrix $\gamma_{k}%
I$, and only keep the $m$ most resent $s_{j}$ and $y_{j}$ in the update
\eqref{qng:eg:cf_bfgs}. Since $m\ll n$ the amount of storage and computations
in each iteration is only a small fraction compared to the regular
\textsc{bfgs}. According to \cite{nowr:06} satisfactory results are often
achieved with $5\leq m\leq20$, even for large problems. Our experiments
confirm this heuristic. Thus for the limited memory \textsc{bfgs} we have
\begin{equation}
H_{k}=\gamma_{k}I+%
\begin{bmatrix}
S_{k} & \gamma_{k}Y_{k}%
\end{bmatrix}%
\begin{bmatrix}
R_{k}^{-\mathsf{T}}(D_{k}+\gamma_{k}Y_{k}^{\mathsf{T}}Y_{k})R_{k}^{-1} &
-R_{k}^{-\mathsf{T}}\\
-R_{k}^{-1} & 0
\end{bmatrix}%
\begin{bmatrix}
S_{k}^{\mathsf{T}}\\
\gamma_{k}Y_{k}^{\mathsf{T}}%
\end{bmatrix}
, \label{eq:L-BFGS}%
\end{equation}
where now
\[
S_{k}=\left[  s_{k-m},\dots,s_{k-1}\right]  ,\quad Y_{k}=\left[  y_{k-m}%
,\dots,y_{k-1}\right]  ,\quad D_{k}=\operatorname*{diag}\left[  s_{k-m}%
^{\mathsf{T}}y_{k-m},\dots,s_{k-1}^{\mathsf{T}}y_{k-1}\right]
\]
and
\[
R_{k}=%
\begin{bmatrix}
s_{k-m}^{\mathsf{T}}y_{k-m} & s_{k-m}^{\mathsf{T}}y_{k-m+1} & \cdots &
s_{k-m}^{\mathsf{T}}y_{k-1}\\
0 & s_{k-m+1}^{\mathsf{T}}y_{k-m+1} & \cdots & s_{k-m+1}^{\mathsf{T}}y_{k-1}\\
\vdots &  & \ddots & \vdots\\
0 & \cdots & 0 & s_{k-1}^{\mathsf{T}}y_{k-1}%
\end{bmatrix}
.
\]

\subsection{Limited memory BFGS on Grassmannians}

Analyzing the \textsc{l-bfgs} update above with the intent of modifying it to
be applicable on Grassmannians, we observe the following:

\begin{enumerate}
\item The columns in the matrices $S_{k}$ and $Y_{k}$ represent tangents, and
as such, they are defined on a specific point of the manifold. In each
iteration we need to parallel transport these vectors to the next tangent
space. Assuming $s_{k}$ and $y_{k}$ are vectorized forms of \eqref{eq:grSk}
and \eqref{eq:grYk} the transport amounts to computing $\bar{S_{k}}%
=(I_{r}\mathbin{\scalebox{1.12}{$\displaystyle\varotimes$}} T(t_{k}))S_{k}$
and $\bar{Y_{k}}=(I_{r}\mathbin{\scalebox{1.12}{$\displaystyle\varotimes$}}
T(t_{k}))Y_{k}$ where $T(t_{k})$ is the Grassmann transport matrix.

\item The matrices $R_{k}$ and $D_{k}$ contain inner products between
tangents. Fortunately, the inner products are invariant with respect to
parallel transporting. Given vectors $u,v\in\mathbf{T}_{X_{k}}$ and a
transport matrix $T$ from $\mathbf{T}_{X_{k}}$ to $\mathbf{T}_{X_{k+1}}$,
i.e.\ $Tu,Tv\in\mathbf{T}_{X_{k+1}}$, we have that $\langle Tu,Tv\rangle
=\langle u,v\rangle$. This is a direct result from
Theorem~\ref{thm:vecTranProd}, showing that there is no need for modifying
$R_{k}$ or $D_{k}$. Because of this property one may wonder whether the
transport matrix $T$ is orthogonal, but this is not the case, $T^{\mathsf{T}%
}T\neq I$.

\item Recalling the relation from equation \eqref{eq:ident1} between local and
global coordinate representation of an operator, we conclude that the global
representation is necessarily a singular matrix, simply because the local
coordinate representation of the operator is a smaller matrix. The same is
true for the Hessian using global coordinates. But by construction, the
\textsc{l-bfgs} update $H_{k}$ in \eqref{eq:L-BFGS} is positive definite and
thus nonsingular. This causes no problem since $\mathbf{T}_{X_{k}}$ is an
invariant subspace of $H_{k}$, i.e.\ if $v\in\mathbf{T}_{X_{k}}$ then
$H_{k}v\in\mathbf{T}_{X_{k}}$, see Lemma~\ref{lem:invSS}. Similarly for the
solution of the (quasi-)Newton equations \eqref{eq:QN} since $y_{k}%
\in\mathbf{T}_{X_{k}}$ and $H_{k}:\mathbf{T}_{X_{k}}\rightarrow\mathbf{T}%
_{X_{k}}$, then obviously $p_{k}\in\mathbf{T}_{X_{k}}$. This is valid for
$H_{k}$ from both \eqref{qng:eg:cf_bfgs} and \eqref{eq:L-BFGS}.
\end{enumerate}

\begin{lemma}
\label{lem:invSS} The tangent space $\mathbf{T}_{X_{k}}$ is an invariant
subspace of the operator obtained by the \textsc{l-bfgs} update.
\end{lemma}

\begin{proof}
This is straightforward. Simply observe that for a vector $v_{k}\in
\mathbf{T}_{X_{k}}$ we have that $H_{k}v_{k}$ is a linear combination of
vectors, and all of them belong to $\mathbf{T}_{X_{k}}$.
\end{proof}

\textsc{l-bfgs} algorithms are intended for large scale problems where the
storage of the full Hessian may not be possible. With this in mind we realize
that the computation and storage of the orthogonal complement $X_{\perp}$,
which is used in local coordinate implementations, may not be practical. For
large and sparse problems it is more economical to do the {parallel}
transports explicitly than to update a basis for the tangent space. The
computational time is reasonable since only $2(m-1)$ vectors are {parallel}
transported each step and $m$ is usually very small compared to the dimensions
of the Hessian.

\section{Quasi-Newton methods for the best multilinear rank approximation of a
tensor\label{qng:sec:brApp}}

In this section we apply the algorithms developed in the last three sections
to the tensor approximation problem described earlier. Recall from
Section~\ref{sec:approx-max} that the best multilinear rank-$(p,q,r)$
approximation of a general tensor is equivalent to the maximization of
\[
\Phi(X,Y,Z)=\frac{1}{2}\lVert\mathcal{A}\cdot(X,Y,Z)\rVert_{F}^{2}\quad\text{
s.t. }\quad X^{\mathsf{T}}X=I,\quad Y^{\mathsf{T}}Y=I,\quad Z^{\mathsf{T}}Z=I,
\]
where $\mathcal{A}\in\mathbb{R}^{l\times m\times n}$ and $X\in\mathbb{R}%
^{l\times p}$, $Y\in\mathbb{R}^{m\times q}$, $Z\in\mathbb{R}^{n\times r}$.
Recall also that $X,Y,Z$ may be regarded as elements of $\operatorname*{Gr}%
(l,p)$, $\operatorname*{Gr}(m,q)$, and $\operatorname*{Gr}(n,r)$ respectively
and $\Phi$ may be regarded as a function defined on a product of the three
Grassmannians. The Grassmann gradient of $\Phi$ will consist of three
parts. Setting $\mathcal{F}=\mathcal{A}\cdot(X,Y,Z)$, one can show that in
global coordinates the gradient is the triplet $\nabla\Phi=(\Pi_{X}\Phi
_{x},\Pi_{Y}\Phi_{y},\Pi_{Z}\Phi_{z})$, where
\begin{align}
\Pi_{X}\Phi_{x} &  =\langle\mathcal{A}\cdot(\Pi_{X},Y,Z),\mathcal{F}%
\rangle_{-1}\quad\in\mathbb{R}^{l\times p}, &  &  \Pi_{X}=I-XX^{\mathsf{T}%
},\label{eq:globGradX}\\
\Pi_{Y}\Phi_{y} &  =\langle\mathcal{A}\cdot(X,\Pi_{Y},Z),\mathcal{F}%
\rangle_{-2}\quad\in\mathbb{R}^{m\times q}, &  &  \Pi_{Y}=I-YY^{\mathsf{T}%
},\label{eq:globGradY}\\
\Pi_{Z}\Phi_{z} &  =\langle\mathcal{A}\cdot(X,Y,\Pi_{Z}),\mathcal{F}%
\rangle_{-3}\quad\in\mathbb{R}^{n\times r}, &  &  \Pi_{Z}=I-ZZ^{\mathsf{T}%
},\label{eq:globGradZ}%
\end{align}
and $\Phi_{x}=\partial\Phi/\partial X$, $\Phi_{y}=\partial\Phi/\partial Y$ and
$\Phi_{z}=\partial\Phi/\partial Z$, see equation \eqref{eq:df}. For derivation
of these formulas\footnotemark[3] see \cite{elsa09}.

To obtain the corresponding expressions in local coordinates we observe that a
projection matrix can also be written as $\Pi_{X}=X_{\perp}X_{\perp
}^{\mathsf{T}}$. Then for tangent vectors $\Delta_{x}\in\mathbf{T}_{X}$, we
have
\[
\Delta_{x}=\Pi_{X}\Delta_{x}=X_{\perp}X_{\perp}^{\mathsf{T}}\Delta_{x}\equiv
X_{\perp}D_{x},
\]
which gives the local coordinates of $\Delta_{x}$ as $X_{\perp}^{\mathsf{T}%
}\Delta_{x}=D_{x}$. The practical implication of these manipulations is that
in local coordinates we simply replace the projection matrices $\Pi_{X}%
,\Pi_{Y},\Pi_{Z}$ with $X_{\perp}^{\mathsf{T}},Y_{\perp}^{\mathsf{T}}%
,Z_{\perp}^{\mathsf{T}}$. We get $\nabla\widehat{\Phi}=(X_{\perp}^{\mathsf{T}%
}\Phi_{x},Y_{\perp}^{\mathsf{T}}\Phi_{y},Z_{\perp}^{\mathsf{T}}\Phi_{z})$,
where
\begin{align}
X_{\perp}^{\mathsf{T}}\Phi_{x}  &  =\langle\mathcal{A}\cdot(X_{\perp
},Y,Z),\mathcal{F}\rangle_{-1}\quad\in\mathbb{R}^{(l-p)\times p}%
,\label{eq:locGradX}\\
Y_{\perp}^{\mathsf{T}}\Phi_{y}  &  =\langle\mathcal{A}\cdot(X,Y_{\perp
},Z),\mathcal{F}\rangle_{-2}\quad\in\mathbb{R}^{(m-q)\times q}%
,\label{eq:locGradY}\\
Z_{\perp}^{\mathsf{T}}\Phi_{z}  &  =\langle\mathcal{A}\cdot(X,Y,Z_{\perp
}),\mathcal{F}\rangle_{-3}\quad\in\mathbb{R}^{(n-r)\times r}.
\label{eq:locGradZ}%
\end{align}
Note that the expressions of the gradient in global and local coordinates are
different. In order to distinguish between them we put a hat on the gradient,
i.e.\ $\nabla\widehat{\Phi}$, when it is expressed in local coordinates.

\subsection{General expression for Grassmann gradients and
Hessians\label{sec:genDeriv}}

In the general case we will have an order-$k$ tensor $\mathcal{A}\in
\mathbb{R}^{n_{1}\times\dots\times n_{k}}$ and the objective function takes
the form
\[
\Phi(X_{1},\dots,X_{k})=\frac{1}{2}\lVert\mathcal{A}\cdot(X_{1},\dots
,X_{k})\rVert_{F}^{2}.
\]
The low rank approximation problem becomes
\[
\max\Phi(X_{1},\dots,X_{k})\quad\text{ s.t. }\quad X_{i}^{\mathsf{T}}%
X_{i}=I,\quad i=1,\dots,k.
\]
The same procedure used to derive the gradients for the order-3 case can be
used for the general case. The results are obvious modifications of what we
have for 3-tensors. First we introduce matrices $X_{i\perp}$, $i=1,\dots,k,$
such that each $[X_{i}\,\,X_{i\perp}]$ forms an orthogonal matrix and we
define the tensors
\begin{align}
\mathcal{F}  &  =\mathcal{A}\cdot(X_{1},X_{2},X_{3},\dots,X_{k}),\nonumber\\
\mathcal{B}_{1}  &  =\mathcal{A}\cdot(X_{1\perp},X_{2},X_{3},\dots
,X_{k}),\nonumber\\
\mathcal{B}_{2}  &  =\mathcal{A}\cdot(X_{1},X_{2\perp},X_{3},\dots
,X_{k}),\nonumber\\
&  \vdots\label{eq:intermediate}\\
\mathcal{B}_{k-1}  &  =\mathcal{A}\cdot(X_{1},\dots,X_{k-2},X_{k-1,\perp
},X_{k}),\nonumber\\
\mathcal{B}_{k}  &  =\mathcal{A}\cdot(X_{1},\dots,X_{k-2},X_{k-1},X_{k\perp
}).\nonumber
\end{align}
The Grassmann gradient of the objective function in local coordinates is given
by the $k$-tuple
\[
\nabla\widehat{\Phi}=\left(  \Phi_{1},\Phi_{2},\dots,\Phi_{k}\right)
,\qquad\Phi_{i}=\langle\mathcal{B}_{i},\mathcal{F}\rangle_{-i},\quad
i=1,2,\dots,k.
\]
Each $\Phi_{i}$ is an $(n_{i}-r_{i})\times r_{i}$ matrix representing a
tangent in $\mathbf{T}_{X_{i}}$. To obtain the corresponding global coordinate
representation, simply replace each $X_{i\perp}$ with the projection
$\Pi_{X_{i}}=I-X_{i}X_{i}^{\mathsf{T}}$.

We will also give the expression of the Hessian since we may wish to
initialize our approximate Hessian with the exact Hessian. Furthermore, in our
numerical experiments in Section~\ref{qng:sec:comEx}, the expression for the
Hessian will be useful for checking whether our algorithms have indeed arrived
at a local maximum. In order to express the Hessian, we will need to introduce
the additional variables
\begin{equation}%
\begin{matrix}
\mathcal{C}_{12} &  &  & \\
\mathcal{C}_{13} & \mathcal{C}_{23} &  & \\
\vdots & \vdots & \ddots & \\
\mathcal{C}_{1,k} & \mathcal{C}_{2,k} & \cdots & \mathcal{C}_{k-1,k}%
\end{matrix}
\label{eq:Ctensors}%
\end{equation}
where each term is a multilinear tensor-matrix product involving the tensor
$\mathcal{A}$ and a subset of the matrices in $\{X_{1},\dots,X_{k},X_{1\perp
},\dots,X_{k\perp}\}$. The subscripts $i$ and $j$ in $\mathcal{C}_{ij}$
indicate that $X_{i\perp}$ and $X_{j\perp}$ are multiplied in the $i$th and
$j$th mode of $\mathcal{A}$, respectively. All other modes are multiplied with
the corresponding $X_{d}$, $d\neq i$ and $d\neq j$. For example we have
\begin{align*}
\mathcal{C}_{12}  &  =\mathcal{A}\cdot(X_{1\perp},X_{2\perp},X_{3},\dots
,X_{k}),\\
\mathcal{C}_{24}  &  =\mathcal{A}\cdot(X_{1},X_{2\perp},X_{3},X_{4\perp}%
,X_{5},\dots,X_{k}),\\
\mathcal{C}_{k-1,k}  &  =\mathcal{A}\cdot(X_{1},\dots,X_{k-2},X_{k-1,\perp
},X_{k\perp}).
\end{align*}
Together with $\mathcal{B}_{1},\dots,\mathcal{B}_{k}$, introduced earlier, one
can express the complete Grassmann Hessian of the objective function
$\Phi(X_{1},\dots,X_{k})$. The derivation of the Hessian is somewhat tricky.
The interested reader should refer to \cite{elsa09} for details. In this paper
we only state the final result in a form that can be directly implemented in a solver.

The diagonal blocks of the Hessian are Sylvester operators and have the form
\begin{align*}
\mathcal{H}_{ii}(D_{i})  &  = \langle\mathcal{B}_{i}, \mathcal{B}_{i}
\rangle_{-i} D_{i} - D_{i} \langle\mathcal{F},\mathcal{F} \rangle_{-i} , \quad
i = 1,2,\dots, k.
\end{align*}
The off-diagonal block operators are
\begin{align*}
\mathcal{H}_{12}(D_{2})  &  = \langle\langle\mathcal{C}_{12},\mathcal{F}
\rangle_{-(1,2)} ,D_{2} \rangle_{2,4;1,2} + \langle\langle\mathcal{B}%
_{1},\mathcal{B}_{2} \rangle_{-(1,2)} ,D_{2} \rangle_{4,2;1,2},\\
&  \vdots\\
\mathcal{H}_{ij}(D_{j})  &  = \langle\langle\mathcal{C}_{ij},\mathcal{F}
\rangle_{-(i,j)} ,D_{j} \rangle_{2,4;1,2} + \langle\langle\mathcal{B}%
_{i},\mathcal{B}_{j} \rangle_{-(i,j)} ,D_{j} \rangle_{4,2;1,2},
\end{align*}
where $i \neq j$, $i < j$ and $i, j = 1,2,\dots,k.$ {See Appendix \ref{app:1}
for definition of the contracted products $\langle\ \cdot,\cdot\rangle
_{-(i,j)}$.}

\section{Best multilinear rank approximation of a symmetric
tensor\label{qng:sec:symmCase}}

Recall from Section~\ref{sec:Tensors} that an order-$k$ tensor $\mathcal{S}%
\in\mathbb{R}^{n\times\dots\times n}$ is called \textit{symmetric} if
\[
s_{i_{\sigma(1)}\cdots i_{\sigma(k)}}=s_{i_{1}\cdots i_{k}},\qquad i_{1}%
,\dots,i_{k}\in\{1,\dots,n\},
\]
where $\sigma\in\mathfrak{S}_{k}$, the set of all permutations with $k$
integers. For example, a third order cubical tensor $\mathcal{S}\in
\mathbb{R}^{n\times n\times n}$ is symmetric iff
\[
s_{ijk}=s_{ikj}=s_{jik}=s_{jki}=s_{kij}=s_{kji}%
\]
for all $i,j,k\in\{1,\dots,n\}$. The definition given above is equivalent to
the usual definition given in, say \cite{greub78}; see \cite{comon08} for a
proof of this simple equivalence. Recall also that the set of all order-$k$
dimension-$n$ symmetric tensors is denoted $\mathsf{S}^{k}(\mathbb{R}^{n})$.
This is a subspace of $\mathbb{R}^{n\times\dots\times n}$ and%
\[
\dim\mathsf{S}^{k}(\mathbb{R}^{n})=\dbinom{n+k-1}{k}.
\]

\begin{lemma}
\label{lem:symrank}If $\mathcal{S}\in\mathsf{S}^{k}(\mathbb{C}^{n})$ and
$\operatorname*{rank}(\mathcal{S})=(r_{1},\dots,r_{k})$, then%
\[
r_{1}=\dots=r_{k}.
\]
In other words, the multilinear rank of a symmetric tensor is always of the
form $(r,\dots,r)$ for some $r$. We will write $r_{\mathsf{S}}(\mathcal{S})$
for this common value. Furthermore, we have a multilinear decomposition of the
following form
\begin{equation}
\mathcal{S}=(X,X,\dots,X)\cdot{\mathcal{C}}, \label{eq:symdecomp}%
\end{equation}
where $\mathcal{C}\in\mathsf{S}^{k}(\mathbb{R}^{r})$ and $X\in
\operatorname*{O}(n,r)$.
\end{lemma}

\begin{proof}
The ranks $r_{i}$ being equal follows from observing that the matricizations
$S^{(1)},\dots,S^{(k)}$ of $\mathcal{S}$ are, due to symmetry, all equal. The
factorization \eqref{eq:symdecomp} is a consequence of the higher order
singular value decomposition (\textsc{hosvd})~\cite{latha00}.
\end{proof}

In application where noise is an inevitable factor, we would like to study
instead the approximation problem
\[
\mathcal{S}\approx(X,X,\dots,X)\cdot{\mathcal{C}},
\]
instead of the exact decomposition in \eqref{eq:symdecomp}. More precisely, we
want to solve%
\begin{equation}
\min\{\lVert\mathcal{S}-\mathcal{T}\rVert_{F}\mid\text{$\mathcal{T}\in$%
}\mathsf{S}^{k}(\mathbb{R}^{n}),\;r_{\mathsf{S}}(\mathcal{T})\leq r\}.
\label{eq:symProb}%
\end{equation}

Similar analysis as in the general case shows that the minimization problem
\eqref{eq:symProb} can be reformulated as a maximization of $\lVert
\mathcal{S}\cdot(X,\dots,X)\rVert_{F}$, with the constraint $X^{\mathsf{T}%
}X=I$. The objective function becomes $\Phi(X)=\frac{1}{2}\langle
\mathcal{F},\mathcal{F}\rangle$ where now $\mathcal{F}=\mathcal{S}%
\cdot(X,\dots,X)$. Observe that the symmetric tensor approximation problem is
defined on one Grassmannian only, regardless of the order of the tensor.
These problems require much less storage and computations compared to a
general problem of the same dimensions. Applications involving symmetric
tensors are found in signal processing, independent component analysis, and
the analysis of multivariate cumulants in statistics
\cite{comon96,comon08,kofid02, delath00a,delath04b,delath04c,comon06,ML,LM}. We refer
interested readers to \cite{comon08} for discussion of a different notion of
rank for symmetric tensors.

\subsection{The symmetric Grassmann gradient}

The same procedure for deriving the gradient for the general case can be used
to obtain the gradient for the symmetric case. In particular it involves the
very same terms as the nonsymmetric gradient with obvious modifications. It is
straightforward to show that, due to symmetry of $\mathcal{S}$,
\[
\langle\mathcal{S}\cdot(\Pi_{X},X,X),\mathcal{F}\rangle_{-1}=\langle
\mathcal{S}\cdot(X,\Pi_{X},X),\mathcal{F}\rangle_{-2}=\langle\mathcal{S}%
\cdot(X,X,\Pi_{X}),\mathcal{F}\rangle_{-3}.
\]
We will use the first expression without loss of generality. In which case,
the Grassmann gradient in global coordinates becomes
\begin{equation}
\nabla\Phi=\Pi_{X}\Phi_{x}=3\langle\mathcal{S}\cdot(\Pi_{X},X,X),\mathcal{F}%
\rangle_{-1}, \label{eq:symGradGlob}%
\end{equation}
where $\Pi_{X}=I-XX^{\mathsf{T}}$; and in local coordinate it is
\begin{equation}
\nabla\widehat{\Phi}=X_{\perp}\Phi_{x}=3\langle\mathcal{S}\cdot(X_{\perp
},X,X),\mathcal{F}\rangle_{-1}, \label{eq:symGradLoc}%
\end{equation}
where $X_{\perp}$ is the orthogonal complement of $X$. Compare these with
equations \eqref{eq:globGradX}--\eqref{eq:globGradZ} for the general case.

\subsection{The symmetric Grassmann Hessian}

As for the general case discussed in \cite{elsa09}, we may identify the second
order terms in the Taylor expansion of $\Phi(X_{\Delta}(t))$. There are 15
second order terms and all have the form
\[
\langle\Delta,\mathcal{H}_{\ast}(\Delta)\rangle,\quad\Delta\in\mathbf{T}%
_{X}\text{ and }X\in\operatorname*{Gr}(n,r),
\]
for some linear operator $\mathcal{H}_{\ast}$. Two specific examples are
\begin{align*}
\langle\Delta,\mathcal{H}_{11}(\Delta)\rangle &  =\left\langle \Delta
,\langle\mathcal{B}_{1},\mathcal{B}_{1}\rangle_{-1}\Delta-\Delta
\langle\mathcal{F},\mathcal{F}\rangle_{-1}\right\rangle ,\\
\langle\Delta,\mathcal{H}_{12}(\Delta)\rangle &  =\left\langle \Delta
,\langle\langle\mathcal{C}_{12},\mathcal{F}\rangle_{-(1,2)},\Delta
\rangle_{2,4;1,2}+\langle\langle\mathcal{B}_{1},\mathcal{B}_{2}\rangle
_{-(1,2)},\Delta\rangle_{4,2;1,2}\right\rangle ,
\end{align*}
where $\mathcal{B}_{1}=\mathcal{S}\cdot(\Pi_{X},X,X)$, $\mathcal{B}%
_{2}=\mathcal{S}\cdot(X,\Pi_{X},X)$ and $\mathcal{C}_{12}=\mathcal{S}\cdot
(\Pi_{X},\Pi_{X},X)$. The subscripts $1$ and $2$ indicate that the projection
matrix $\Pi_{X}$ is multiplied with $\mathcal{S}$ in the first and second mode
respectively. Not surprisingly, analysis of these terms reveals equality among
the second order terms due to the symmetry of $\mathcal{S}$. Gathering like
terms and summing up the expressions, we see that the Hessian is a sum of
three different terms,
\begin{align}
\langle\Delta,\mathcal{H}_{1}(\Delta)\rangle=  &  \langle\Delta,3\langle
\mathcal{B}_{1},\mathcal{B}_{1}\rangle_{-1}\Delta-3\Pi_{X}\Delta
\langle\mathcal{F},\mathcal{F}\rangle_{-1}\rangle,\label{eq:hess1}\\
\langle\Delta,\mathcal{H}_{2}(\Delta)\rangle=  &  \langle\Delta,6\langle
\langle\mathcal{C}_{12},\mathcal{F}\rangle_{-(1,2)},\Delta\rangle
_{2,4;1,2}\rangle,\label{eq:hess2}\\
\langle\Delta,\mathcal{H}_{3}(\Delta)\rangle=  &  \langle\Delta,6\langle
\langle\mathcal{B}_{1},\mathcal{B}_{2}\rangle_{-(1,2)},\Delta\rangle
_{4,2;1,2}\rangle. \label{eq:hess3}%
\end{align}
So the action of the Hessian on a tangent is simply%
\[
\mathcal{H}(\Delta)=\mathcal{H}_{1}(\Delta)+\mathcal{H}_{2}(\Delta
)+\mathcal{H}_{3}(\Delta).
\]
Observe that the second term in \eqref{eq:hess1} arises from the fact that the
objective function is defined on a Grassmannian, see \cite{edelm99} for details.

\subsection{General expression for Grassmann gradients and Hessians for a
symmetric tensor\label{sec:genDerivSym}}

With the analysis and expressions for symmetric $3$-tensors at hand,
generalization to symmetric $k$-tensors is straightforward. We will
only state the final results and in local coordinates. Assume we have an
order-$k$ symmetric tensor $\mathcal{S}\in\mathsf{S}^{k}(\mathbb{R}^{n})$. The
corresponding symmetric low rank tensor approximation problem is written as
\[
\max\Phi(X)=\max\frac{1}{2}\lVert\mathcal{S}\cdot(X,\dots,X)\rVert_{F}%
^{2}\quad\text{ s.t. }\quad X\in\operatorname*{Gr}(n,r).
\]
Using the tensor products
\begin{align*}
\mathcal{F}  &  =\mathcal{S}\cdot(X,X,\dots,X)\qquad\text{$X$ appears $k$
times},\\
\mathcal{B}_{1}  &  =\mathcal{S}\cdot(X_{\perp},X,\dots,X)\qquad\text{$X$
appears $k-1$ times},
\end{align*}
where $X_{\perp}$ is such that $[X\,\,X_{\perp}]$ forms an orthogonal matrix,
the Grassmann gradient becomes
\[
\nabla\Phi=k\langle\mathcal{B}_{1},F\rangle_{-1}.
\]
Observe that the symmetric case involves the very same tensor products
$\mathcal{B}_{i}$ as in the general case (given in Section~\ref{sec:genDeriv})
but due to the symmetry of the problem all terms are equal.

We also introduce tensor-matrix multilinear products $\mathcal{C}_{ij}$
similar to those in equation \eqref{eq:Ctensors}. Two specific examples are
\begin{align*}
\mathcal{C}_{12}  &  =\mathcal{S}\cdot(X_{\perp},X_{\perp},X,\dots,X),\\
\mathcal{C}_{24}  &  =\mathcal{S}\cdot(X,X_{\perp},X,X_{\perp},X,\dots,X).
\end{align*}
In general $\mathcal{C}_{ij}$, where $i\neq j$, $i<j$ and $i,j=1,\dots,k$, is
a multilinear product of two $X_{\perp}$'s that are multiplied in the $i$th
and $j$th mode of $\mathcal{S}$. All other modes are multiplied with $X$.

The second order terms of the Taylor expansion of $\Phi(X)$ contain the
following diagonal block operators
\[
H_{ii}(D)=\langle\mathcal{B}_{i},\mathcal{B}_{i}\rangle_{-i}D-\Delta
\langle\mathcal{F},\mathcal{F}\rangle_{-i},\quad i=1,2,\dots,k.
\]
Again, due to symmetry all these are identical and summing them up we get
\[
H_{\text{diag}}(D)=k\left(  \langle\mathcal{B}_{1},\mathcal{B}_{1}\rangle
_{-1}D-D\langle\mathcal{F},\mathcal{F}\rangle_{-1}\right)  .
\]
The off-diagonal block operators have the form
\begin{align*}
\mathcal{H}_{12}(D)  &  =\langle\langle\mathcal{C}_{12},\mathcal{F}%
\rangle_{-(1,2)},D\rangle_{2,4;1,2}+\langle\langle\mathcal{B}_{1}%
,\mathcal{B}_{2}\rangle_{-(1,2)},D\rangle_{4,2;1,2},\\
&  \vdots\\
\mathcal{H}_{ij}(D)  &  =\langle\langle\mathcal{C}_{ij},\mathcal{F}%
\rangle_{-(i,j)},D\rangle_{2,4;1,2}+\langle\langle\mathcal{B}_{i}%
,\mathcal{B}_{j}\rangle_{-(i,j)},D\rangle_{4,2;1,2},
\end{align*}
where $i\neq j$, $i<j$ and $i,j=1,\dots,k.$ Similarly, due to symmetry all of
them are identical. We have
\[
H_{\text{off-diag}}(D)=k(k-1)\left(  \langle\langle\mathcal{C}_{12}%
,\mathcal{F}\rangle_{-(1,2)},D\rangle_{2,4;1,2}+\langle\langle\mathcal{B}%
_{1},\mathcal{B}_{2}\rangle_{-(1,2)},D\rangle_{4,2;1,2}\right)  .
\]
The complete Grassmann Hessian operator is simply
\[
H=H_{\text{diag}}+H_{\text{off-diag}}.
\]

\subsection{Matricizing the Hessian operator}

The second order terms are described using the canonical inner product on
Grassmannians and contracted tensor products. Next we will derive the
expression of the Hessian as a matrix acting on the vector
$d=\operatorname*{vec}(\Delta)$.

The terms in \eqref{eq:hess1} involve only matrix operations and vectorizing
the second argument in the inner product yields
\begin{align*}
\operatorname*{vec}\left(  (\mathcal{H}_{1}(\Delta)\right)   &  =
\operatorname*{vec}(3\langle\mathcal{B}_{1},\mathcal{B}_{1}\rangle_{-1}%
\Delta-3\Pi_{X}\Delta\langle\mathcal{F},\mathcal{F}\rangle_{-1})\\
&  =3\left(  I\mathbin{\scalebox{1.12}{$\displaystyle\varotimes$}}
\langle\mathcal{B}_{1},\mathcal{B}_{1}\rangle_{-1}-\langle\mathcal{F}%
,\mathcal{F}\rangle_{-1}\mathbin{\scalebox{1.12}{$\displaystyle\varotimes$}}
\Pi_{X}\right)  \operatorname*{vec}(\Delta)\\
&  \equiv H_{1}d.
\end{align*}
The vectorization of the terms from \eqref{eq:hess2} and \eqref{eq:hess3}
involve the 4-tensors
\begin{align*}
\mathcal{H}_{2}  &  =\langle\mathcal{C}_{12},\mathcal{F} \rangle_{-(1,2)}
\in\mathbb{R}^{n \times n \times r \times r},\\
\mathcal{H}_{3}  &  =\langle\mathcal{B}_{1},\mathcal{B} _{2}\rangle_{-(1,2)}
\in\mathbb{R}^{n \times r \times r \times n},
\end{align*}
and is done using the tensor matricization described in
\cite{elsa09}. We get
\begin{align}
&  \operatorname*{vec}(\langle\mathcal{H}_{2},\Delta\rangle_{2,4;1,2}%
)=H_{2}^{(3,1;4,2)}\operatorname*{vec}(\Delta)\equiv H_{2}d,\\
&  \operatorname*{vec}(\langle\mathcal{H}_{3},\Delta\rangle_{4,2;1,2}%
)=H_{3}^{(3,1;2,4)}\operatorname*{vec}(\Delta)\equiv H_{3}d.
\end{align}
In $\mathcal{H}_{2}$ we map indices of the first and third mode to row-indices
and indices of the second and fourth mode to column-indices obtaining the
matrix $H_{2}^{(3,1;4,2)}$. In this way the contractions in the matrix-vector
product coincide with the tensor-matrix contractions. Similarly for $H_{3}$.
The matrix form of the Hessian becomes
\[
H=H_{1}+H_{2}+H_{3}.
\]
To obtain the Hessian in local coordinates we replace $\Pi_{X}$ with
$X_{\perp}$ in the computations of the factors involved and thereafter perform
the same matricization procedure.

\section{Examples\label{sec:ex}}

We will now give two small explicit examples to illustrate the computations
involved in the algorithms for tensor approximation described before.

\begin{example}
In this example we will compute the gradient of the objective function, both
in global and in local coordinates. Let the $3\times3\times3$ tensor
$\mathcal{A}$ be given by {\small
\[
\mathcal{A}(:,:,1)=%
\begin{bmatrix}
9 & -3 & 8\\
2 & 7 & 0\\
7 & 0 & -1
\end{bmatrix}
,\,\mathcal{A}(:,:,2)=%
\begin{bmatrix}
2 & 7 & 0\\
-7 & 5 & -3\\
0 & -3 & 1
\end{bmatrix}
,\,\mathcal{A}(:,:,3)=%
\begin{bmatrix}
3 & 0 & -2\\
0 & 4 & -1\\
0 & -2 & 1
\end{bmatrix}
.
\]
} Let the current point of the product manifold be given by $(X,Y,Z)$ where
\[
X=Y=Z=%
\begin{bmatrix}
1\\
0\\
0
\end{bmatrix}
,\quad\text{ and }\quad\Pi_{X}=\Pi_{Y}=\Pi_{Z}=%
\begin{bmatrix}
0 & 0 & 0\\
0 & 1 & 0\\
0 & 0 & 1
\end{bmatrix}
\]
are the corresponding projection matrices onto the three tangent spaces. The
expression for the Grassmann gradient at the current iterate is given by
\eqref{eq:globGradX}--\eqref{eq:globGradZ}. The intermediate quantities,
cf.\ equation \eqref{eq:intermediate}, needed in the calculations of the
Grassmann gradient are%
\begin{align}
\mathcal{F}  &  =\mathcal{A}\cdot(X,Y,Z)=9,\nonumber\\
\mathcal{B}_{x}  &  =\mathcal{A}\cdot(\Pi_{X}%
,Y,Z)=(0\,\,\,2\,\,\,7)^{\mathsf{T}},\label{eq:inter1}\\
\mathcal{B}_{y}  &  =\mathcal{A}\cdot(X,\Pi_{Y}%
,Z)=(0\,\,\,-3\,\,\,8)^{\mathsf{T}},\nonumber\\
\mathcal{B}_{z}  &  =\mathcal{A}\cdot(X,Y,\Pi_{Z}%
)=(0\,\,\,2\,\,\,3)^{\mathsf{T}},\nonumber
\end{align}
and the Grassmann gradient in global coordinates is given by
\[
\nabla\Phi=\left(  \langle\mathcal{B}_{x},\mathcal{F}\rangle_{-1}%
,\langle\mathcal{B}_{y},\mathcal{F}\rangle_{-2},\langle\mathcal{B}%
_{z},\mathcal{F}\rangle_{-3}\right)  =\left(
\begin{bmatrix}
0\\
18\\
63
\end{bmatrix}
,%
\begin{bmatrix}
0\\
-27\\
72
\end{bmatrix}
,%
\begin{bmatrix}
0\\
18\\
27
\end{bmatrix}
\right)  .
\]
To compute the Grassmann gradient in local coordinates we need a basis for the
tangent spaces. For the current iterate we choose
\[
{X}_{\perp}=Y_{\perp}=Z_{\perp}=%
\begin{bmatrix}
0 & 0\\
1 & 0\\
0 & 1
\end{bmatrix}
,
\]
as the corresponding basis matrices for the tangent spaces at $X$, $Y$ and
$Z$. Obviously $[X\,\,X_{\perp}]$, $[Y\,\,Y_{\perp}]$ and $[Z\,\,Z_{\perp}]$
are orthogonal and $X^{\mathsf{T}}X_{\perp}=Y^{\mathsf{T}}Y_{\perp
}=Z^{\mathsf{T}}Z_{\perp}=0$. Replacing the projection matrices $\Pi_{X}$,
$\Pi_{Y}$ and $\Pi_{Z}$ by the orthogonal complements $X_{\perp}$, $Y_{\perp}$
and $Z_{\perp}$ in \eqref{eq:inter1}, we obtain $\widehat{\mathcal{B}}%
_{x},\widehat{\mathcal{B}}_{y},\widehat{\mathcal{B}}_{z}$, and thus the local
coordinate representation of the Grassmann gradient is given by
\[
\nabla\widehat{\Phi}=(\langle\widehat{\mathcal{B}}_{x},\mathcal{F}\rangle
_{-1},\langle\widehat{\mathcal{B}}_{y},\mathcal{F}\rangle_{-2},\langle
\widehat{\mathcal{B}}_{z},\mathcal{F}\rangle_{-3})=\left(
\begin{bmatrix}
18\\
63
\end{bmatrix}
,%
\begin{bmatrix}
-27\\
72
\end{bmatrix}
,%
\begin{bmatrix}
18\\
27
\end{bmatrix}
\right)  .
\]
Recall that we use a hat to distinguish local coordinate representation from
global coordinate representation. The local coordinate representation is
depending on the choice of basis matrices for the tangent spaces. A different
choice of $X_{\perp}$, $Y_{\perp}$ and $Z_{\perp}$ would yield a different
representation of $\nabla\widehat{\Phi}$.
\end{example}

\begin{figure}[t]
\centering
\includegraphics[width=0.8\textwidth]{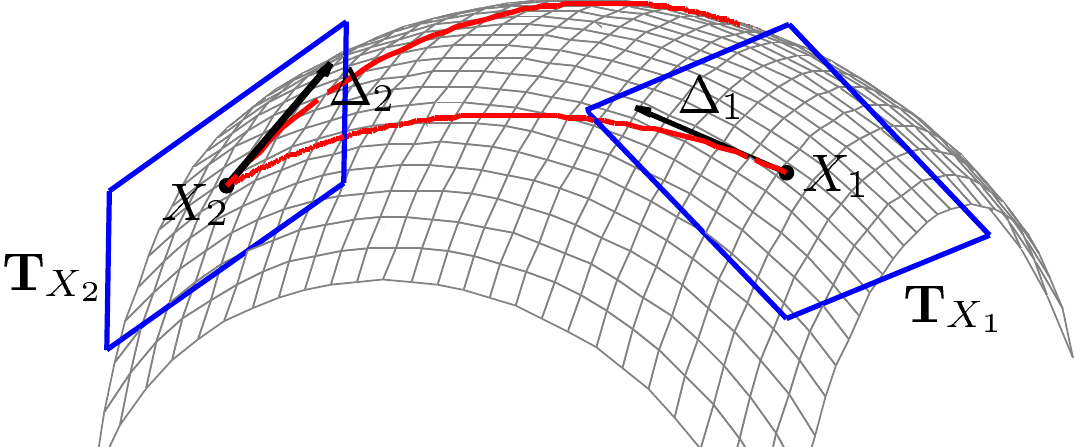}\caption{Pictorial
depiction of the main algorithmic procedure for the Grassmannian
$\operatorname*{Gr}(3,1)$, which is simply the sphere $S^{2}$.}%
\label{qng:man}%
\end{figure}

\begin{example}
\label{sec:exp2}Next we will illustrate the parallel transport of tangent
vectors along geodesics on a product of Grassmannians. Let the tensor
$\mathcal{A}$, the current iterate, and the corresponding gradient be the same
as in the previous example. Introduce tangent vectors
\[
\Delta=(\Delta_{x},\,\Delta_{y},\,\Delta_{z})=\left(
\begin{bmatrix}
0\\
-1\\
0
\end{bmatrix}
,%
\begin{bmatrix}
0\\
0\\
1
\end{bmatrix}
,%
\begin{bmatrix}
0\\
1\\
0
\end{bmatrix}
\right)  .
\]
Clearly we have $X^{\mathsf{T}}\Delta_{x}=Y^{\mathsf{T}}\Delta_{y}%
=Z^{\mathsf{T}}\Delta_{z}=0$. The tangent $\Delta$ will determine the geodesic
path from the current point and in turn the transport of the Grassmann
gradient (see Figure~\ref{qng:man}). We may also verify that $\nabla\Phi$ is
indeed a tangent of the product Grassmannian at the current iterate.

The thin or compact \textsc{svd}s, written $\Delta_{\ast}=U_{\ast}\cdot
\Sigma_{\ast}\cdot{V_{\ast}^{\mathsf{T}}},$ of the tangents are%
\[
\Delta_{x}=%
\begin{bmatrix}
0\\
-1\\
0
\end{bmatrix}
\cdot1\cdot1,\quad\Delta_{y}=%
\begin{bmatrix}
0\\
0\\
1
\end{bmatrix}
\cdot1\cdot1,\quad\Delta_{z}=%
\begin{bmatrix}
0\\
1\\
0
\end{bmatrix}
\cdot1\cdot1.
\]

The transport matrix, cf.\ equation \eqref{eq:parVecTran}, in the direction
$\Delta_{x}$ at $X$ with a step size $t=\pi/4$ is given by
\[
\left.  T_{X,\Delta_{x}}(t)\right\vert _{t=\pi/4}=%
\begin{bmatrix}
XV_{x} & U_{x}%
\end{bmatrix}%
\begin{bmatrix}
-\sin\Sigma_{x}(\pi/4)\\
\cos\Sigma_{x}(\pi/4)
\end{bmatrix}
{U_{x}^{\mathsf{T}}}+(I-U_{x}{U}_{x}^{\mathsf{T}})=%
\begin{bmatrix}
1 & 1/\sqrt{2} & 0\\
0 & 1/\sqrt{2} & 0\\
0 & 0 & 1
\end{bmatrix}
.
\]
Similarly, it is straightforward to calculate
\[
\left.  T_{Y,\Delta_{y}}(t)\right\vert _{t=\pi/4}=%
\begin{bmatrix}
1 & 0 & -1/\sqrt{2}\\
0 & 1 & 0\\
0 & 0 & 1/\sqrt{2}%
\end{bmatrix}
\quad\text{ and }\quad\left.  T_{Z,\Delta_{z}}(t)\right\vert _{t=\pi/4}=%
\begin{bmatrix}
1 & -1/\sqrt{2} & 0\\
0 & 1/\sqrt{2} & 0\\
0 & 0 & 1
\end{bmatrix}
.
\]
Parallel transporting one tangent we get $T_{X,\Delta_{x}}(\pi/4)\Delta
_{x}=(-1/\sqrt{2},\,\,-1/\sqrt{2},\,\,0)^{\mathsf{T}}$. For all tangents in
$\Delta$ and $\nabla\Phi$ we get
\begin{align*}
\Delta(t)|_{t=\pi/4}  &  =\left(
\begin{bmatrix}
-1/\sqrt{2}\\
-1/\sqrt{2}\\
0
\end{bmatrix}
,%
\begin{bmatrix}
-1/\sqrt{2}\\
0\\
1/\sqrt{2}%
\end{bmatrix}
,%
\begin{bmatrix}
-1/\sqrt{2}\\
1/\sqrt{2}\\
0
\end{bmatrix}
\right)  ,\\
\nabla\Phi(t)|_{t=\pi/4}  &  =\left(
\begin{bmatrix}
18/\sqrt{2}\\
18/\sqrt{2}\\
63
\end{bmatrix}
,%
\begin{bmatrix}
-72/\sqrt{2}\\
-27\\
72/\sqrt{2}%
\end{bmatrix}
,%
\begin{bmatrix}
-18/\sqrt{2}\\
18/\sqrt{2}\\
27
\end{bmatrix}
\right)  .
\end{align*}

The above are calculations in global coordinates. In local coordinates we
parallel transport the basis matrices $X_{\perp}$, $Y_{\perp}$ and $Z_{\perp}$
so that the local coordinate representation of a tangent is the same as in the
previous point. The computations are given by equation \eqref{eq:geoPerp} and
in this example we get
\[
X_{\perp}(\pi/4)=%
\begin{bmatrix}
1/\sqrt{2} & 0\\
1/\sqrt{2} & 0\\
0 & 1
\end{bmatrix}
,\,\,Y_{\perp}(\pi/4)=%
\begin{bmatrix}
0 & -1/\sqrt{2}\\
1 & 0\\
0 & 1/\sqrt{2}%
\end{bmatrix}
,\,\,Z_{\perp}(\pi/4)=%
\begin{bmatrix}
-1/\sqrt{2} & 0\\
1/\sqrt{2} & 0\\
0 & 1
\end{bmatrix}
,
\]
i.e.\ the second and third columns of each transport matrix due to the
specific choice of $X_{\perp}$, $Y_{\perp}$ and $Z_{\perp}$.

Taking a step of size $t=\pi/4$ from $X$, $Y$ and $Z$ along the specified
geodesic we arrive~at
\[
X(\pi/4)=%
\begin{bmatrix}
1/\sqrt{2}\\
-1/\sqrt{2}\\
0
\end{bmatrix}
,\quad Y(\pi/4)=%
\begin{bmatrix}
1/\sqrt{2}\\
0\\
1/\sqrt{2}%
\end{bmatrix}
,\quad Z(\pi/4)=%
\begin{bmatrix}
1/\sqrt{2}\\
1/\sqrt{2}\\
0
\end{bmatrix}
.
\]
The value of the objective function at the starting point is $\Phi
(X,Y,Z)=40.5$ and at the new point is $\Phi\left(  X(\pi/4),Y(\pi
/4),X(\pi/4)\right)  \approx45.5625$, an increment as expected.
\end{example}

Figure~\ref{qng:man} illustrates the procedures involved in the algorithms on
the Grassmannian $\operatorname*{Gr}(3,1)$, which we may regard as the
$2$-sphere $S^{2}$ (unit sphere in $\mathbb{R}^{3}$). For the best rank-$1$
tensor approximation of a $3\times3\times3$ tensor, the optimization take{s}
place on a product of three spheres $S^{2}\times S^{2}\times S^{2}$, one for
each vector that needs to be determined. The procedure starts at a point
$X_{1}$ and a direction of ascent\footnote{Recall that we are maximizing
$\Phi$, therefore `ascent' as opposed to `descent'.}, the tangent $\Delta
_{1}\in\mathbf{T}_{X_{1}}$, is obtained through some method. Next we perform a
movement of the point $X_{1}$ along the geodesic defined by $\Delta_{1}$.
Geodesics on spheres are just great circles. At the new point $X_{2}%
\in\operatorname*{Gr}(3,1)$ we repeat the procedure, i.e.\ determine a new
direction of ascent $\Delta_{2}$ and take a geodesic step in this direction.

\section{Numerical experiments and computational
complexity\label{qng:sec:comEx}}

All algorithms described here and the object oriented Grassmann classes
required for them are available for download as two \textsc{matlab} packages
\cite{savas08a} and \cite{savas08b}. We encourage our readers to try them out.

\subsection{Initialization and stopping condition\label{qng:sec:exp}}

We will now test the actual performance of our algorithms with a few large
numerical examples. All algorithms in a given test are started with the same
initial points on a Grassmannian, represented as truncated singular
matrices from the \textsc{hosvd} and a number of additional \textit{higher
order orthogonal iterations}---\textsc{hooi} iterations \cite{latha00,
latha00a}, which are introduced to make the initial Hessian of $\Phi$ negative
definite. The number of initial \textsc{hooi} iterations ranges between $5$
and $50$ depending on the size of the problem. The \textsc{bfgs} algorithm is
either started with (possibly a modification of) the exact Hessian or a scaled
identity matrix according to \cite[pp.~143]{nowr:06}. The \textsc{l-bfgs}
algorithm is always started with a scaled identity matrix but one can modify
the number of columns $m$ in the matrices representing the Hessian
approximation, see equation \eqref{eq:L-BFGS}. This number is between $5$ and
$30$. {Although we use the \textsc{hosvd} to initialize our
algorithms, any other reasonable initialization procedure would work} as long
as the initial Hessian approximate is negative definite. The quasi-Newton
methods can be used as stand-alone algorithms for solving the tensor
approximation problem as well as other problems defined on Grassmannians.

In the following figures, the $y$-axis measures the norm of the relative
gradient, i.e.\ $\lVert\nabla\Phi(X)\rVert/\lVert\Phi(X)\rVert$, and the
$x$-axis shows iterations. This ratio is also used as our stopping condition,
which typically requires that $\lVert\nabla\Phi(X)\rVert/\lVert\Phi
(X)\rVert\approx10^{-13}$, the machine precision of our computer. At a true
local maximizer the gradient of the objective function is zero and its Hessian
is negative definite. In the various figures we present convergence results
for four principally different algorithms. These are (1)
quasi-Newton-Grassmann with \textsc{bfgs}, (2) quasi-Newton-Grassmann with
\textsc{l-bfgs}, (3) Newton-Grassmann, denoted with \textsc{ng} and (4)
\textsc{hooi} which is an alternating least squares approach. In addition, the
tags for \textsc{bfgs} methods may be accompanied by \textsc{i} or \textsc{h}
indicating whether the initial Hessian was a scaled identity matrix or the
exact Hessian, respectively.

\subsection{Experimental results}

We run all our numerical experiments in \textsc{matlab} on a MacBook with a
$2.4$-GHz Intel Core 2 Duo processor and $4$~GB of physical memory.

\begin{figure}[t]
\centering
\includegraphics[width=0.49\textwidth]{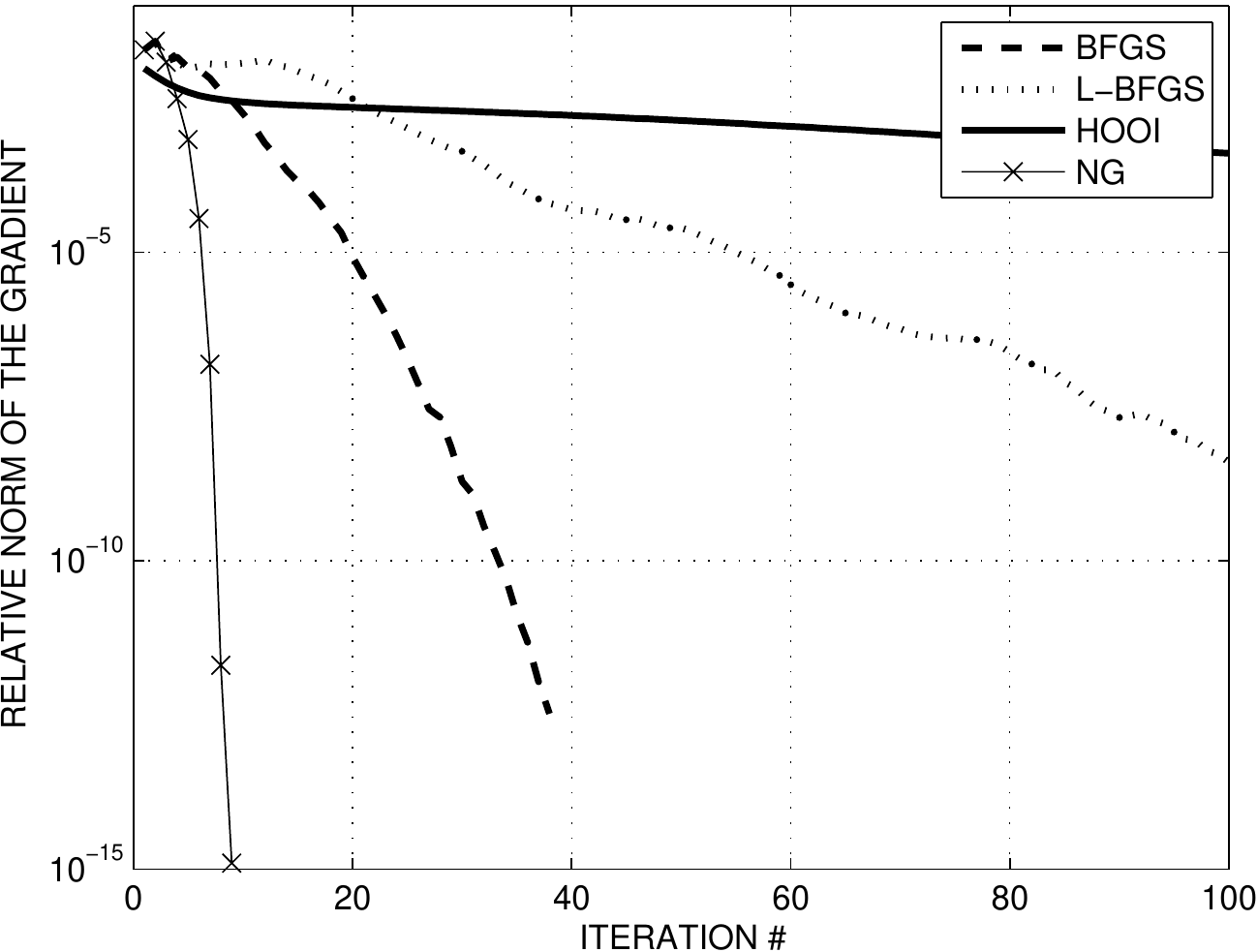}
\includegraphics[width=0.49\textwidth]{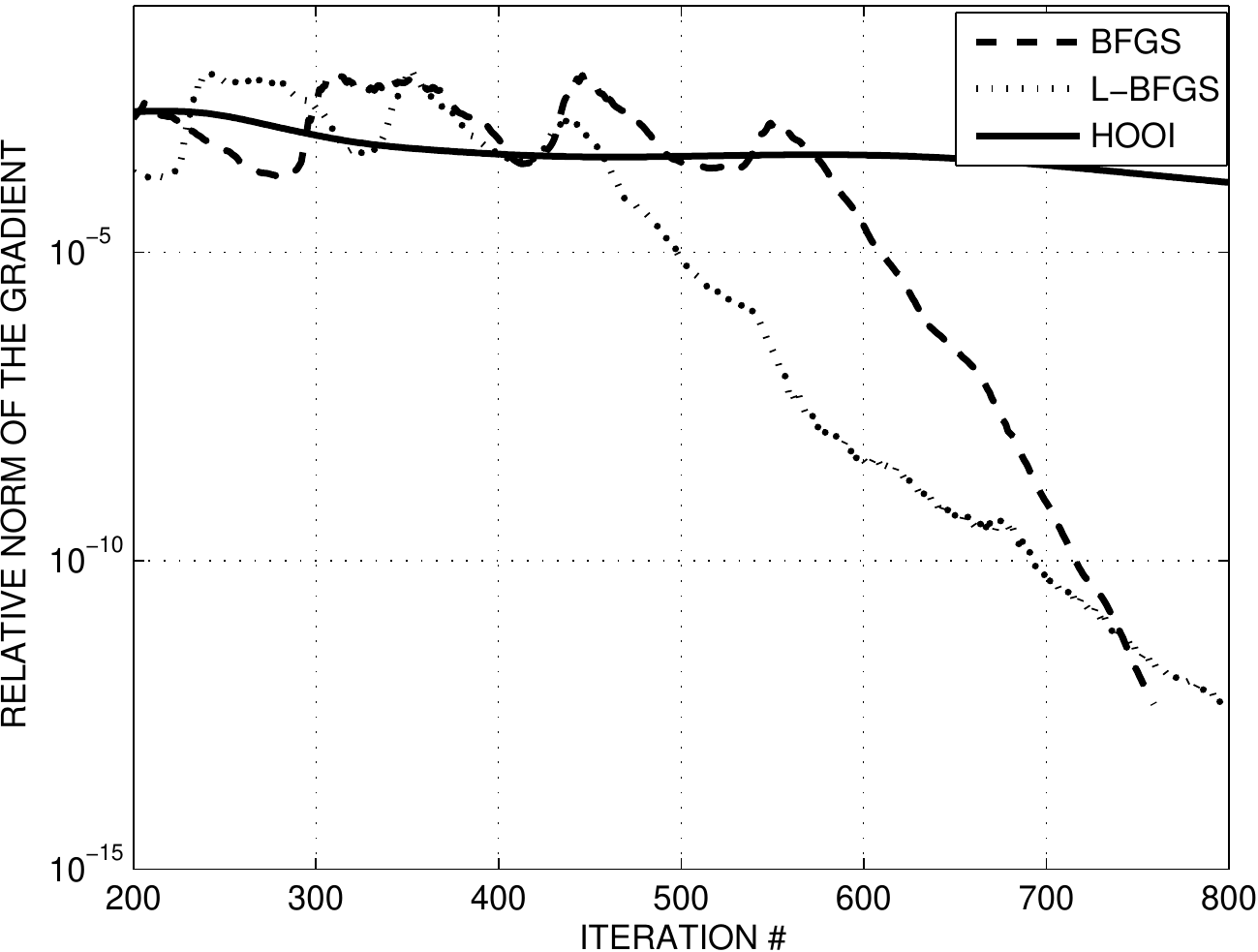}\caption{\textit{Left:}
A $20\times20\times20$ tensor is approximated by a rank-$(5,5,5)$ tensor.
\textsc{bfgs} initiated with the exact Hessian; In \textsc{l-bfgs} $m=5$.
\textit{Right:} A $100\times100\times100$ tensor is approximated by a
rank-$(5,10,20)$ tensor. In this case the initial Hessian is a scaled identity
and $m=10$.}%
\label{qng:fig:1}%
\end{figure}

Figure~\ref{qng:fig:1} shows convergence results for two tests with tensors
generated with $N(0,1)$-distributed values. In the left plot a $20\times20\times20$ tensor is approximated
with a rank-$(5,5,5)$ tensor. One can observe superlinear convergence in the
\textsc{bfgs} method. The right plot shows convergence results of a
$100\times100\times100$ tensor approximated with a rank-$(5,10,20)$ tensor.
Both \textsc{bfgs} and \textsc{l-bfgs} methods exhibit rapid convergence in
the vicinity of a stationary point.

\begin{figure}[t]
\centering
\includegraphics[width=0.49\textwidth]{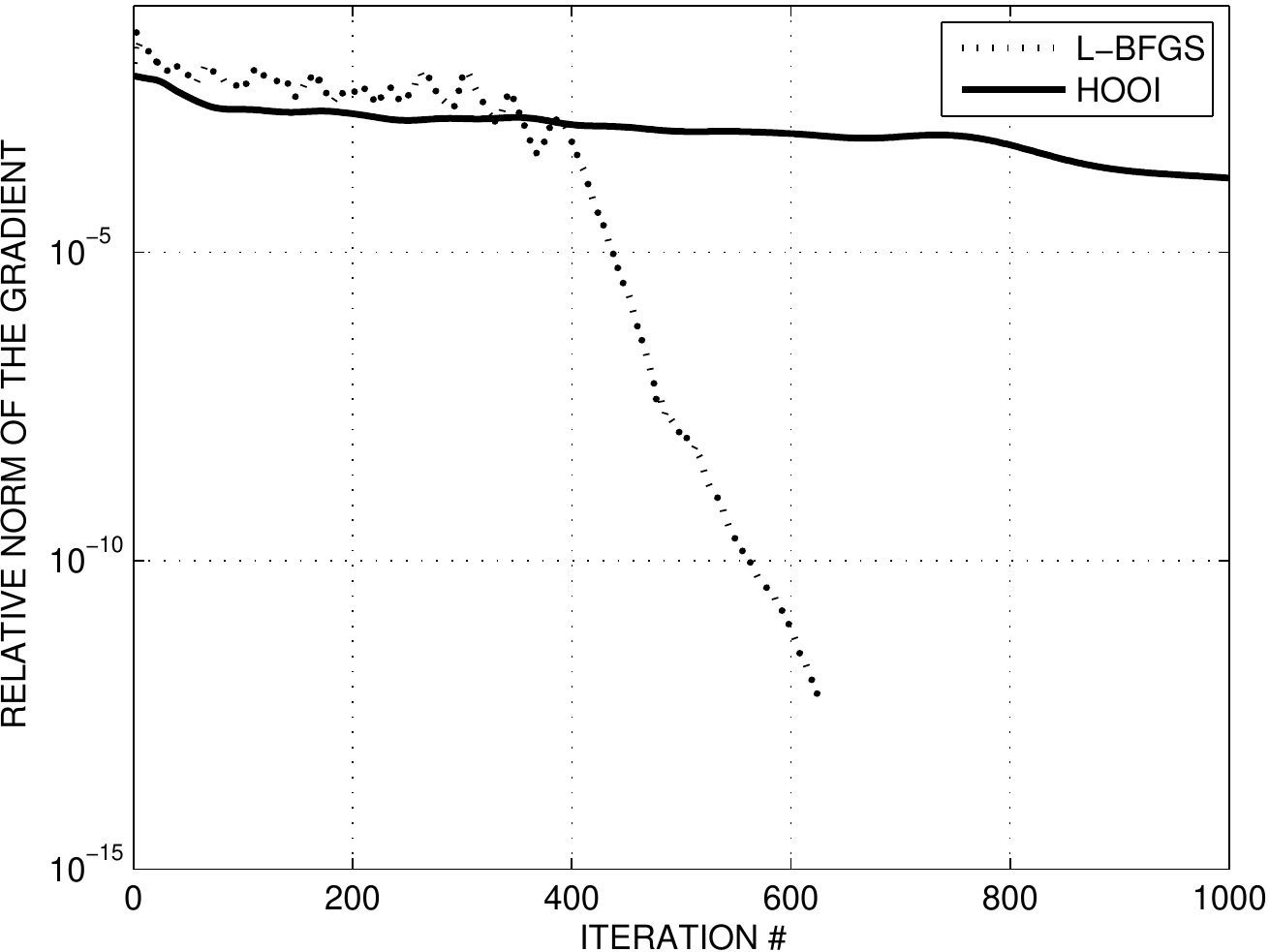}
\includegraphics[width=0.49\textwidth]{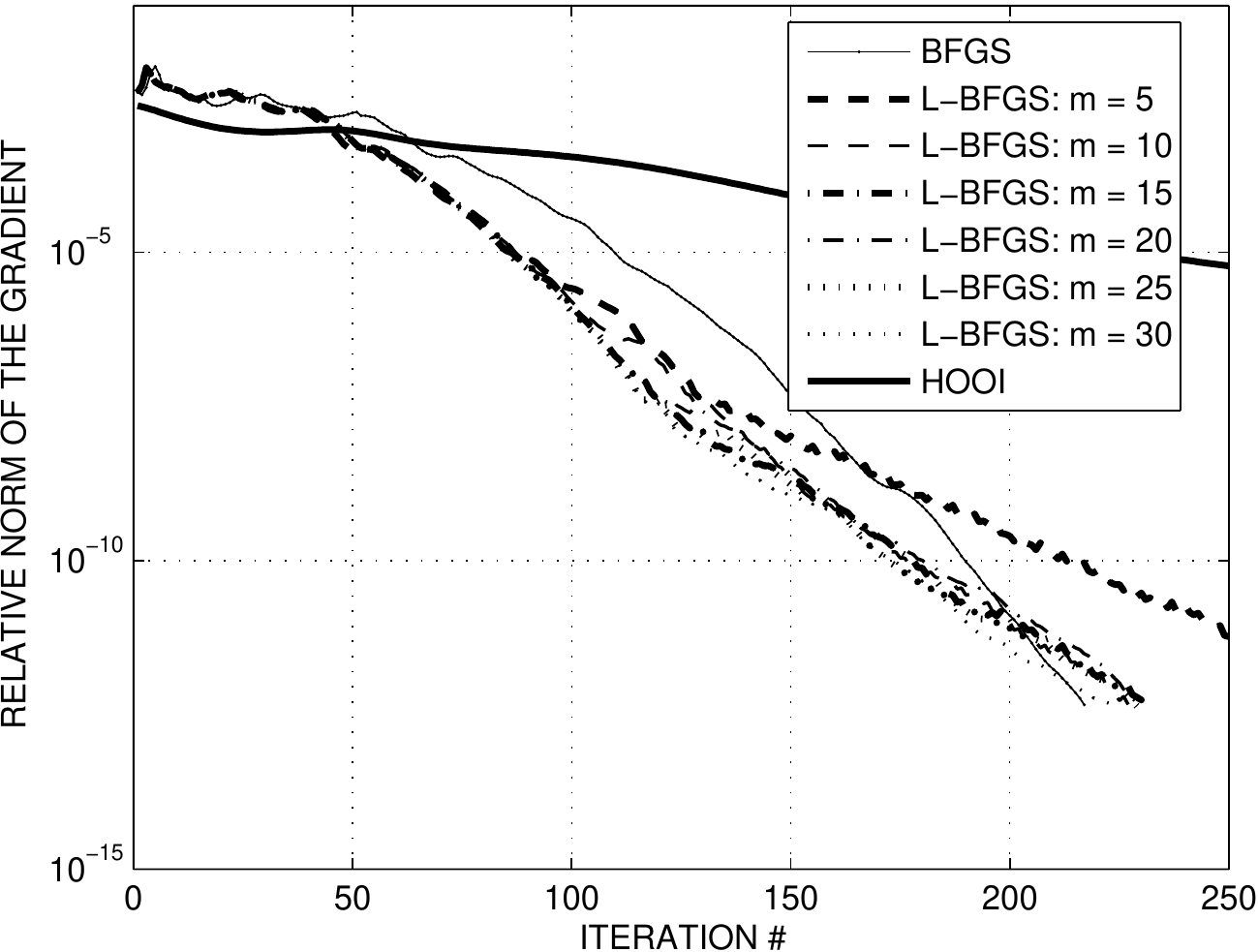}\caption{\textit{Left:}
Convergence plots of a $200\times200\times200$ tensor approximated by a
rank-$(10,10,10)$ tensor. \textit{Right:} Effect of varying $m$ in
\textsc{l-bfgs}. A $50\times50\times50$ tensor approximated by a
rank-$(20,20,20)$ tensor with $m=5,10,15,20,25,30$.}%
\label{qng:fig:2}%
\end{figure}

Figure~\ref{qng:fig:2} (\textit{left}) shows convergence for an even larger
$200\times200\times200$ tensor approximated by a tensor of rank-$(10,10,10)$
using \textsc{l-bfgs} with $m=20$. In the right plot we approximate a $50\times50\times50$ tensor by
a rank-$(20,20,20)$ tensor where we vary over a range of values of $m$ in the
\textsc{l-bfgs} algorithm, namely, $m=5,10,15,20,25,30$. $m=5$ gives (in
general) slightly poorer performance, otherwise the different runs cannot be
distinguished. In other words, our Grassmann \textsc{l-bfgs} algorithm can in
practice work as well as our Grassmann \textsc{bfgs} algorithm, just as one
would expect (from the numerical experiments performed) in the Euclidean case.

\begin{figure}[t]
\centering
\includegraphics[width=0.49\textwidth]{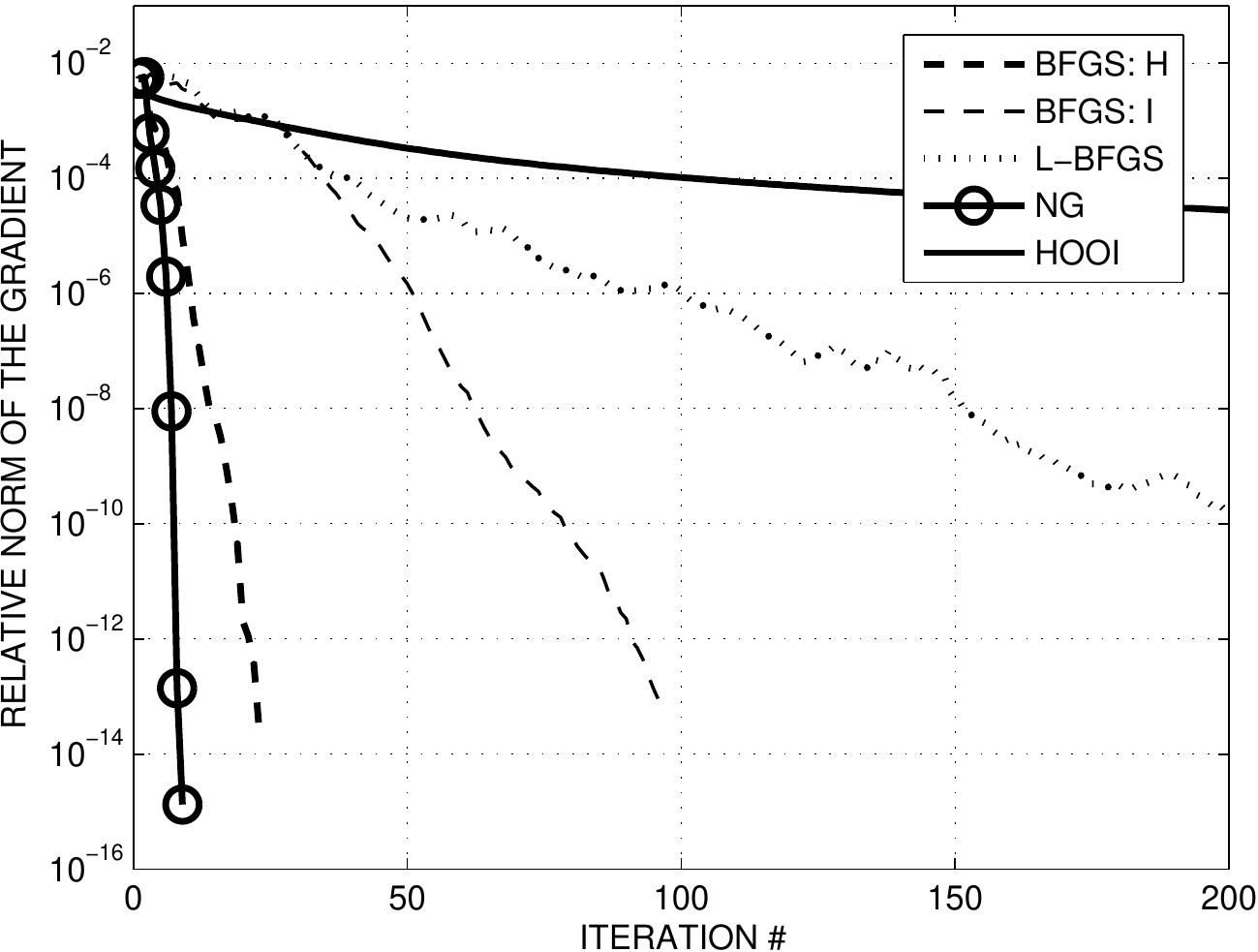}
\includegraphics[width=0.49\textwidth]{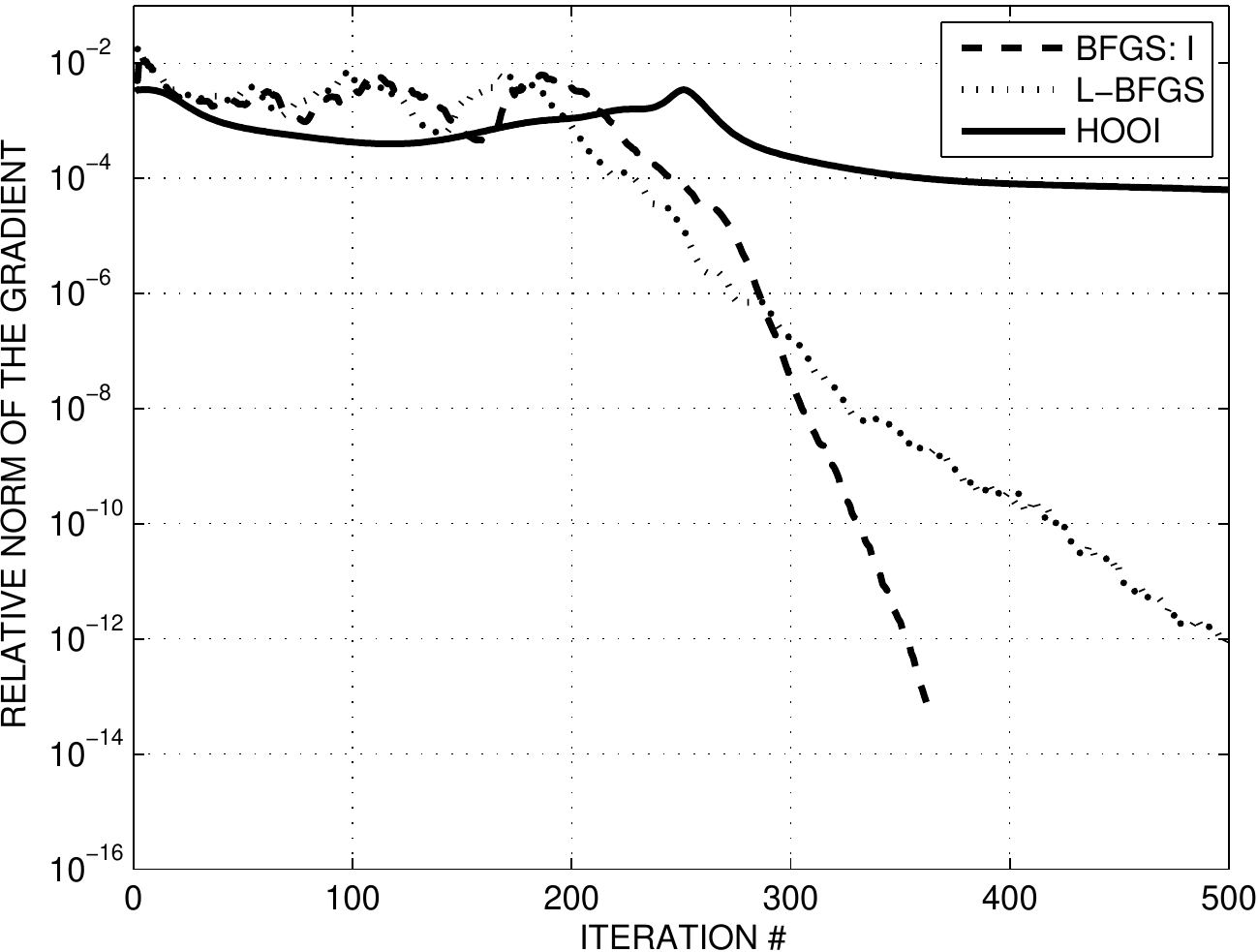}\caption{\textit{Left:}
A $50\times50\times50$ symmetric tensor is approximated by a rank-$5$
symmetric tensor; $m=10$. \textit{Right:} Here we have a $100\times
100\times100$ symmetric tensor approximated by a rank-$20$ symmetric tensor;
$m=10$.}%
\label{qng:fig:3}%
\end{figure}

Figure~\ref{qng:fig:3} shows convergence plots for two symmetric tensor
approximation problems. In the left plot we approximate a symmetric $50\times50\times50$
tensor by a rank-$5$ symmetric tensor. We observe that \textsc{bfgs}
initialized with the exact Hessian (\textsc{bfgs:h} tag) converges much more
rapidly, almost as fast as the Newton-Grassmann method, than when initialized
with a scaled identity matrix (\textsc{bfgs:i} tag). In the right plot we give
convergence results for a $100\times100\times100$ symmetric tensor
approximated by a rank-$20$ symmetric tensor. In both cases $m=10$.

\begin{figure}[t]
\centering
\includegraphics[width=0.49\textwidth]{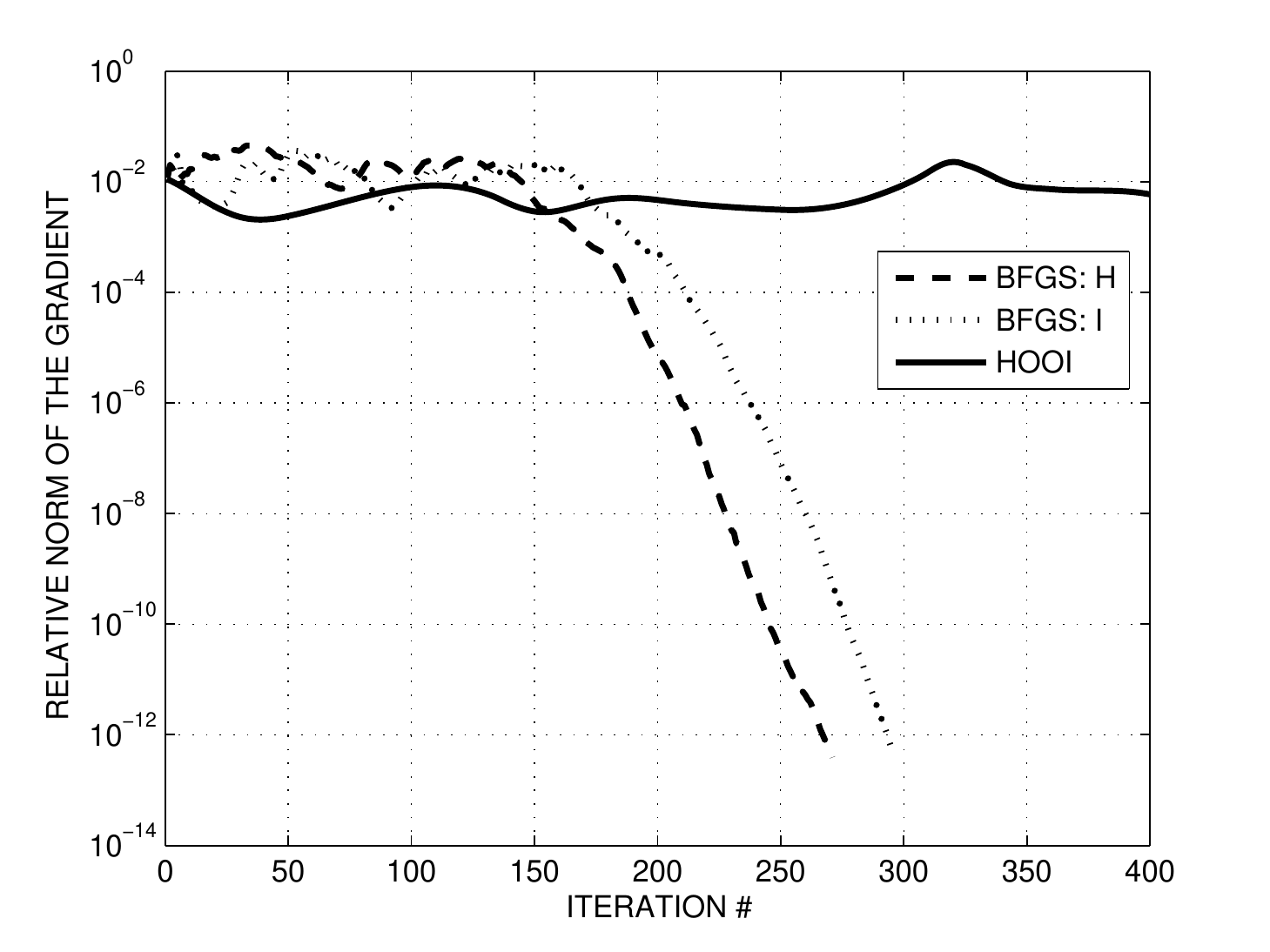}
\includegraphics[width=0.49\textwidth]{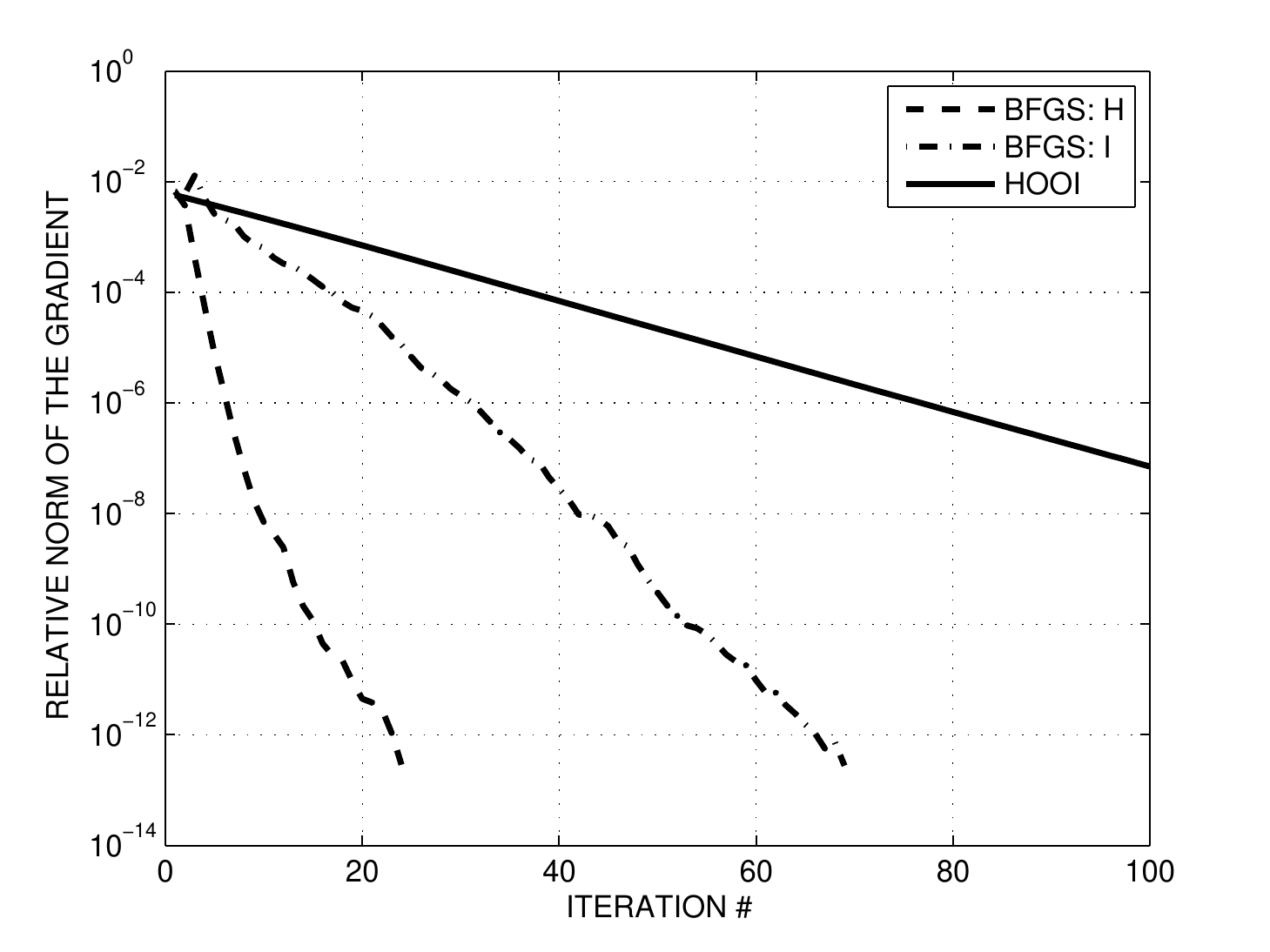}
\caption{\textit{Left:} A $50\times50\times50\times50$ tensor is approximated
by a rank-$(5,5,5,5)$ tensor. \textit{Right:} Here we have a $50\times
50\times50\times50$ symmetric tensor approximated by a rank-$5$ symmetric
tensor.}%
\label{qng:fig:4}%
\end{figure}

In Figure~\ref{qng:fig:4} we show the performance of a local coordinate
implementation of the \textsc{bfgs} algorithm on problems with $4$-tensors.
The first plot shows convergence results for a $50\times50\times50\times50$
tensor approximated by a rank-$(5,5,5,5)$ tensor. The second convergence plot
is for a symmetric $4$-tensor with the same dimensions approximated by a
symmetric rank-$5$ tensor. Again the \textsc{h} and \textsc{i} tags indicate whether the
exact Hessian or a scaled identity is used for initialization.

\begin{figure}[t]
\centering
\includegraphics[width=0.49\textwidth]{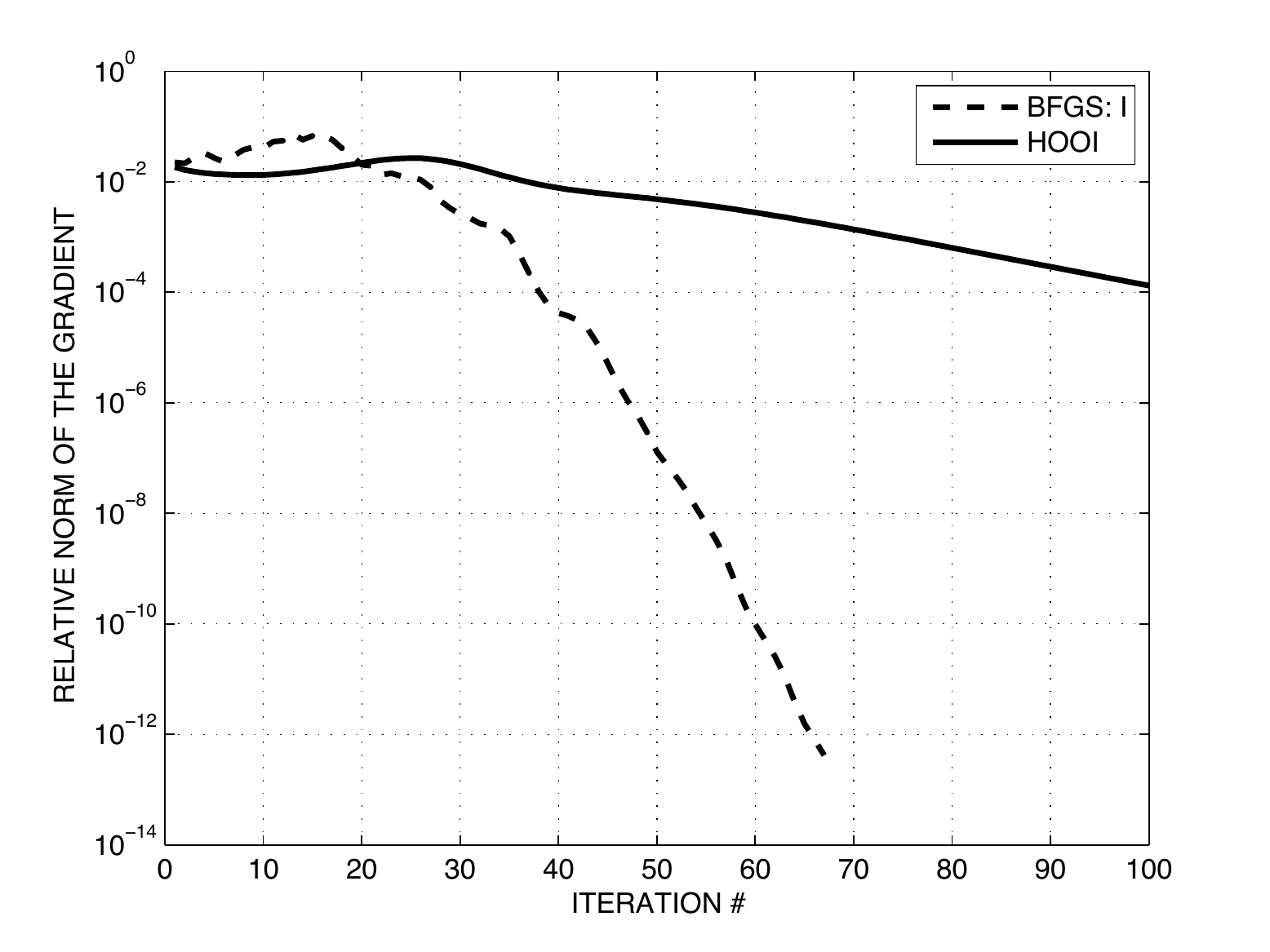}
\includegraphics[width=0.49\textwidth]{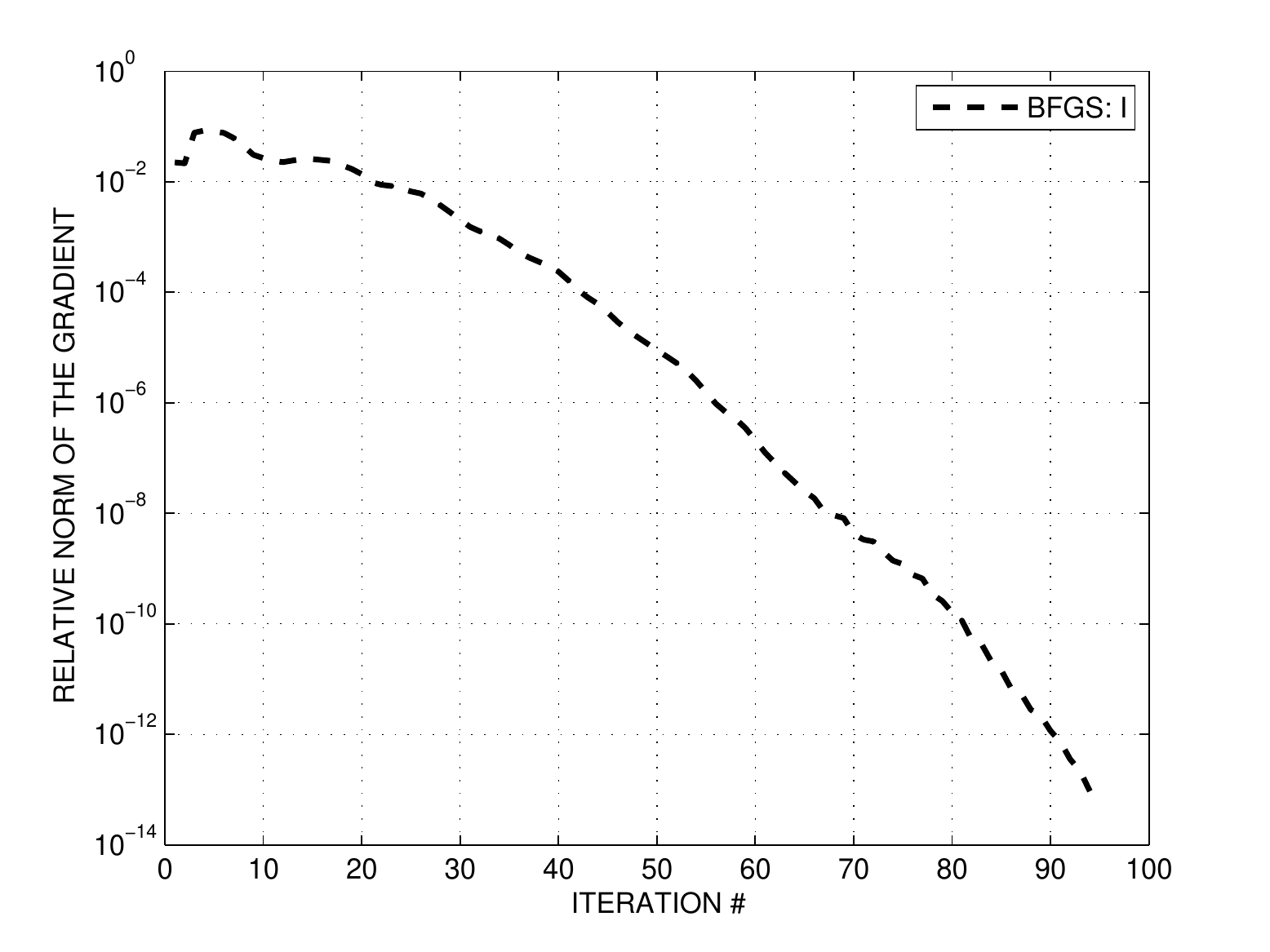}
\caption{\textit{Left:} A $5\times5\times\cdots\times5$ tensor of order-10 is
approximated with a rank-$(2,2,\dots,2)$. \textit{Right:} A `simultaneous'
rank-$5$ approximation of a weighted sum of tensors of orders $2$, $3$, and
$4$.}%
\label{qng:fig:5}%
\end{figure}

We end this section with two unusual examples to illustrate the extent of our
algorithms' applicability: a high order tensor and an objective function that
includes tensors of different orders. The left plot in Figure~\ref{qng:fig:5}
is a high-order example: it shows the convergence of \textsc{bfgs} verses
\textsc{hooi} when approximating an order-$10$ tensor with dimensions
$5\times5\times\cdots\times5$ with a rank-$(2,2,\dots,2)$ tensor. The right
plot in Figure~\ref{qng:fig:5} has an unusual objective function that involves
an order-$2$, an order-$3$, and an order-$4$ tensor,
\[
\Phi(X)=\frac{1}{2!}\lVert S_{2} \cdot (X,X) \rVert_{F}^{2}+\frac{1}{3!}%
\lVert\mathcal{S}_{3}\cdot(X,X,X)\rVert_{F}^{2}+\frac{1}{4!}\lVert
\mathcal{S}_{4}\cdot(X,X,X,X)\rVert_{F}^{2}%
\]
where $S_{2}$ is a $30\times30$ symmetric matrix, $\mathcal{S}_{3}$ is a
$30\times30\times30$ symmetric $3$-tensor, and $\mathcal{S}_{4}$ is a
$30\times30\times30\times30$ symmetric $4$-tensor. Such objective functions have appeared in
\textit{independent component analysis with soft whitening} \cite{delath04b} and in
\textit{principal cumulants components analysis} \cite{ML,LM} where
$S_{2},\mathcal{S}_{3},\mathcal{S}_{4}$ measure the multivariate variance,
kurtosis, skewness respectively (cf.\ Example~\ref{eg:Cum}).
In both examples we observe a fast rate of convergence at the
vicinity of a local minimizer for the \textsc{bfgs} algorithm.

It is evident from the convergence plots here that \textsc{bfgs} and
\textsc{l-bfgs} have faster rate of convergence compared with \textsc{hooi}.
The Newton-Grassmann algorithm takes few iterations but is computationally
more expensive, specifically for larger problems. Our implementation of the
different algorithms in \textsc{matlab} give shortest runtime for the
\textsc{bfgs} and \textsc{l-bfgs} methods. The time for one iteration of
\textsc{bfgs}, \textsc{l-bfgs} and \textsc{hooi} is of the same magnitude for
smaller problems. In larger problems, the \textsc{l-bfgs} performs much faster
tha{n} all other methods.

Our algorithms use the basic arithmetic and data types in the TensorToolbox
\cite{bader06b} for convenience. We use our own object-oriented routines for
operations on Grassmannians and product of Grassmannians, e.g.\ geodesic
movements and {parallel} transports \cite{savas08b}. We note that there are
several different ways to implement \textsc{bfgs} updates \cite{nowr:06}; for
simplicity reasons, we have chosen to update the inverse of the Hessian
approximation. A possibly better alternative will be to update the Cholesky
factors of the approximate Hessians so that one may monitor the approximate
Hessians for indefiniteness during the iterations
\cite{denni96,goldf76,cheng98}.

\subsection{Computational complexity, curse of dimensionality, and
convergence\label{qng:sec:CompComp}}

The Grassmann quasi-Newton methods presented in this report all fit within the
procedural framework given in Algorithm~\ref{qng:alg:frame}.

\begin{algorithm}[t!]
\caption{Algorithmic framework for \textsc{bfgs} and \textsc{l-bfgs} on Grassmannians.}
\label{qng:alg:frame}
\begin{algorithmic}
\STATE Given tensor $\mathcal{A}$ and starting points $(X_0,Y_0,Z_0)\in \operatorname*{Gr}^3$ and an initial Hessian $H_0$
\REPEAT
\STATE \textbf{1} Compute the Grassmann gradient.
\STATE \textbf{2} Parallel transport the Hessian approximation to the new point.
\STATE \textbf{3} Update the Hessian or its compact representation.
\STATE \textbf{4} Solve the quasi-Newton equations to obtain
$\Delta = (\Delta_x, \Delta_y,\Delta_z)$.
\STATE \textbf{5} Move the points $(X_k,Y_k,Z_k)$ along the geodesic curve given by $\Delta$.
\UNTIL{ $\lVert\nabla \widehat{\Phi}\rVert / \Phi < $ TOL}
\end{algorithmic}
\end{algorithm}

\paragraph*{{General case}}

In analyzing computational complexity, we will assume for simplicity that
$\mathcal{A}$ is a general $n\times n\times n$ $3$-tensor being approximated
with a rank-$(r,r,r)$ $3$-tensor. A problem of these dimensions will give rise
to a $3nr\times3nr$ Hessian matrix in global coordinates and a $3(n-r)r\times
3(n-r)r$ Hessian matrix in local coordinates. Table \ref{qng:tab:compCmpl}
gives approximately the amount of computations required in each step of
Algorithm~\ref{qng:alg:frame}. Recall that in \textsc{l-bfgs} $m$ is a small
number, see Section~\ref{qng:sec:lbfgs}. \begin{table}[tbh]
\centering%
\begin{tabular}
[c]{c|c|c|c}\hline\hline
& \textsc{bfgs-gc} & \textsc{bfgs-lc} & \textsc{l-bfgs}\\\hline\hline
\textbf{1} & $6n^{3}r+12n^{2}r^{2}$ & $6n^{3}r+{6}n^{2}r^{2}+6n(n-r)r^{3}$ &
$6n^{3}r+12n^{2}r^{2}$\\
\textbf{2} & $18n^{3}r^{2}$ & --- & $12n^{2}rm$\\
\textbf{3} & $36n^{2}r^{2}$ & $36(n-r)^{2}r^{2}$ & ---\\
\textbf{4} & $18n^{2}r^{2}$ & $18(n-r)^{2}r^{2}$ & $24nrm$\\\hline\hline
\end{tabular}
\caption{Computational complexity of the \textsc{bfgs-gc} (global
coordinates), \textsc{bfgs-lc} (local coordinates) and \textsc{l-bfgs}
algorithms. The numbers in the first column correspond to the steps in
Algorithm~\ref{qng:alg:frame}. }%
\label{qng:tab:compCmpl}%
\end{table}We have omitted terms of lower asymptotic complexity as well as the
cost of point \textbf{5} since that is negligible compared with the costs of
points \textbf{1}--\textbf{4}. For example, the geodesic movement of $X_{k}$
requires the thin \textsc{svd} $U_{x}\Sigma_{x}V_{x}^{\mathsf{T}}=\Delta
_{x}\in\mathbb{R}^{n\times r}$ which takes $6nr^{2}+20r^{3}$ flops (floating
point operations) \cite{golub96}. On the other hand, given the step length $t$
and $U,\Sigma,V$ in \eqref{eq:geodesic}, the actual computation of $X(t)$
amounts to only $4nr^{2}$ flops.

\paragraph*{{Symmetric case}}

The symmetric tensor approximation problem involves the determination of one
$n\times r$ matrix, resulting in an $nr\times nr$ Hessian in global
coordinates and an $(n-r)r\times(n-r)r$ Hessian in local
coordinates. Therefore the complexity of the symmetric problem
differs only by a constant factor from that of the general case.

\paragraph*{{Curse of dimensionality}}

The approximation problem will suffer from the \textit{curse of
dimensionality} when the order of a tensor increases. In general, an
$n\times\dots\times n$ order-$k$ tensor requires the storage of $n^{k}$
entries in memory. The additional memory requirement, mainly for storing the
Hessian, is of order $O(n^{2}r^{2}k^{2})$ for the \textsc{bfgs} methods and
$O(2nrkm)$ for the \textsc{l-bfgs} method, respectively. In the current
approach we assume that the tensor is explicitly given. Our proposed
algorithms are applicable as long as the given tensor fits in memory. There
have been various proposals to deal with the curse of dimensionality using
tensors \cite{khor07,khor09,osel09}. For cases where the tensor is
represented in compact or functional forms our methods can
take direct advantage of these simply by computing
the necessary gradients (and Hessians) using the specific
representations. In fact this was considered in
\cite{morton09} for symmetric tensor approximations.

\paragraph{Convergence}

{There is empirical evidence suggesting that \textsc{als} based algorithms
have fast convergence rate for specific tensors. This was also
pointed out in \cite{delath09}. These are tensors that have inherently low
multilinear rank and the approximating tensor has the correct low ranks, or
tensors that have fast decay in its multilinear singular values \cite{latha00}%
, or a substantial gap in the multilinear singular values at the site of
truncation, e.g. the source tensor is given by a low rank tensor with noise
added. On the other hand not all tensors have gaps or fast decaying
multilinear singular vales. This is specifically true for sparse tensors. It
is still desirable to obtain low rank approximations for these
\textquotedblleft more difficult" tensors. And on these tensors \textsc{als}
performs very poorly, but methods using first and second order derivatives of
the objective function, including the methods presented in this paper perform
good. Among the methods that are currently available, quasi-Newton methods
presented in this paper have the best computational efficiency.}

\section{Related work}\label{sec:relWork}

There are several different approaches to solve the tensor approximation
problem. In this section we will briefly describe them and point out the main
differences with our work. The algorithms most closely related to the
quasi-Newton methods are given in \cite{elsa09,ishteva09,ishteva09b}. All
three references address the best low rank tensor approximation based on the
Grassmannian structure of the problem and use explicit computation of the
Hessian. The obtained Newton equations are solved either fully
\cite{elsa09,ishteva09} or approximately \cite{ishteva09b}. In the latter case
the authors used a truncated conjugate gradient approach to
approximately solve the Newton equations. The iterates are updated using the
more general notion of retractions instead of taking a step along the geodesic
on the manifold. In addition a trust region scheme is incorporated making the
procedure more stable with respect to occasional indefinite Hessians.
The computation of the Hessian is a limiting factor in these algorithms. This is
the case even when the Hessian is not formed explicitly but used implicitly via
its action on a tangent. In our experiments, on moderate-sized problems, e.g.\ $3$-tensors 
of dimensions around $20 \times 20 \times 20$, the \textsc{bfgs} methods noticeably
outperformed Hessian-based methods; and for dimensions around
$100 \times 100 \times 100$, we were unable to get any methods relying 
on Hessians to work despite our best efforts.

There is a different line of algorithms for related tensor
approximation problems based on \textsc{als} and mutltigrid accelerated \textsc{als}
\cite{khor07, khor09}. In our experience, the convergence of \textsc{als}-type
methods depend on the decay of the multilinear singular
values of the given tensor. The exact dependence is unclear but
the relation seems to be that the faster the decay, the
faster the convergence of \textsc{als}. In this regard the
class of functions and operators considered in \cite{khor07,
khor09} appears to possess these favorable properties.

Yet a third approach to obtain low multilinear rank tensor
approximations are the cross methods in \cite{osel08,osel09,flad08}. The
novelty of such methods is that they discard some given
information and retain only a fraction of the original tensor, and as such it
is markedly different from our approach, which uses
all given information to achieve maximal accuracy. In addition, there is an
assumption on the tensor that there exist approximations within pre-specified
bounds and of specific low ranks while we make no such assumptions.

\section{Conclusion\label{qng:sec:conclusion}}

In this paper we studied quasi-Newton algorithms adapted to optimization
problems on Riemannian manifolds. More specifically, we proposed algorithms
with \textsc{bfgs} and \textsc{l-bfgs} updates on a product of Grassmannians that (1) respect the Riemannian metric structure and (2) require
only standard matrix operations in their implementations. Two different
algorithmic implementations are presented: one based on local/intrinsic
coordinates while the other one uses global/embedded coordinates. In
particular, our use of local coordinates is a novelty not previously explored
in other manifold optimization \cite{absil07, absil08, edelm99, gabay82}. We
proved the optimality of our Grassmannian \textsc{bfgs} updates in local
coordinates, showing that the well-known \textsc{bfgs} optimality
\cite{denmor77, densch81} extends to Grassmannian and products of
Grassmannians.

We also applied these algorithms to the problem of determining a best
multilinear rank approximation of a tensor and the analogous (but very
different) problem for a symmetric tensor. While a Newton version of this was
proposed in \cite{elsa09}, here we make substantial improvements with respect
to the Grassmann-Newton algorithm in terms of speed and robustness.
Furthermore, we presented specialized algorithms that take into account the
symmetry in the multilinear approximation of symmetric tensors and related problems.
In addition to the numerical experiments in this paper, we have made our codes freely available
for download \cite{savas08a, savas08b} so that the reader may verify the
speed, accuracy, and robustness of these algorithms for himself.

\appendix

\section{Notation for tensor contractions\label{app:1}}

In this section we define the contracted tensor product notation
used throughout this paper. For given third order tensors
$\mathcal{A}$ and $\mathcal{B}$ we define the following contracted products:
\begin{equation}
\mathcal{C}=\langle\mathcal{A},\mathcal{B}\rangle_{1},\qquad c_{ijkl}%
=\sum_{\lambda}a_{\lambda ij}b_{\lambda kl}.\label{eq:app1}%
\end{equation}
When contracting several indices, with the corresponding
indices of the two arguments being the same, we write
\begin{equation}
\mathcal{C}=\langle\mathcal{A},\mathcal{B}\rangle_{1,2},\qquad c_{ij}%
=\sum_{\lambda,\nu}a_{\lambda\nu i}b_{\lambda\nu j}.\label{eq:app2}%
\end{equation}
The subscript `1' in $\langle\mathcal{A},\mathcal{B}\rangle_{1}$ and
subscripts `1,2' in $\langle\mathcal{A},\mathcal{B}\rangle_{1,2}$
indicate that the contraction is over the first index and
both the first and second indices respectively. If
instead the contraction is to be performed on different indices, we
write
\[
\mathcal{C}=\langle\mathcal{A},\mathcal{B}\rangle_{1;2},\quad c_{ijkl}%
=\sum_{\lambda}a_{\lambda ij}b_{k\lambda l}\quad\text{or}\quad\mathcal{C}%
=\langle\mathcal{A},\mathcal{B}\rangle_{1,3;2,1},\quad c_{ij}=\sum
_{\lambda,\nu}a_{\lambda i\nu}b_{\nu\lambda j}.
\]
The subscripts indicating the indices to be contracted are separated
by a semicolon. It is also convenient to introduce a notation when
contraction is performed in all but one or a few indices. For example
the products in \eqref{eq:app2} and \eqref{eq:app2} may also be written
\[
\langle\mathcal{A},\mathcal{B}\rangle_{1,2}=\langle\mathcal{A},\mathcal{B}%
\rangle_{-3}\quad\text{or}\quad\langle\mathcal{A},\mathcal{B}\rangle
_{1}=\langle\mathcal{A},\mathcal{B}\rangle_{-(2,3)}.
\]

%\bibliography{myBibsENqn}
%\bibliographystyle{siam}

\end{document}